\documentclass[11pt, reqno]{amsart}
\usepackage[margin=1in]{geometry}

\usepackage[T1]{fontenc}
\usepackage[latin1]{inputenc}

\usepackage{amsthm, amsmath, amssymb,booktabs,xcolor,graphicx}
\usepackage[shortlabels]{enumitem}
\usepackage{scalerel,stackengine}
\usepackage[textsize=scriptsize,colorinlistoftodos]{todonotes}
\usepackage{booktabs}
\usepackage[hidelinks, backref=page]{hyperref}

\usepackage[noabbrev,capitalize,nameinlink]{cleveref}
\crefformat{equation}{(#2#1#3)}
\crefrangeformat{equation}{(#3#1#4) to~(#5#2#6)}
\crefmultiformat{equation}{(#2#1#3)}
{ and~(#2#1#3)}{, (#2#1#3)}{ and~(#2#1#3)}

\usepackage{tikz}
\usetikzlibrary{decorations.markings}
\usetikzlibrary{arrows}
\usetikzlibrary{calc}

\newenvironment{proof*}[1][\proofname]{
  
  \begin{proof}[#1]}{\end{proof}}

\newtheorem{theorem}{Theorem}[section]
\newtheorem{proposition}[theorem]{Proposition}
\newtheorem{lemma}[theorem]{Lemma}
\newtheorem{claim}[theorem]{Claim}

\newtheorem{conjecture}[theorem]{Conjecture}

\theoremstyle{definition}
\newtheorem{definition}[theorem]{Definition}
\newtheorem{problem}[theorem]{Problem}
\newtheorem{question}[theorem]{Question}
\newtheorem{example}[theorem]{Example}

\theoremstyle{remark}
\newtheorem{remark}[theorem]{Remark}

\newcommand{\abs}[1]{\left\lvert#1\right\rvert}

\newcommand{\ceil}[1]{\left\lceil #1 \right\rceil}

\newcommand{\sqb}[1]{\left[ #1 \right]}
\newcommand{\set}[1]{\left\{ #1 \right\}}

\renewcommand{\lg}{\text{lg}}
\newcommand{\sm}{\text{sm}}
\newcommand{\eps}{\varepsilon}
\renewcommand{\epsilon}{\varepsilon}

\stackMath
\newcommand\widecheck[1]{
\savestack{\tmpbox}{\stretchto{
  \scaleto{
    \scalerel*[\widthof{\ensuremath{#1}}]{\kern-.6pt\bigwedge\kern-.6pt}
    {\rule[-\textheight/2]{1ex}{\textheight}}
  }{\textheight}
}{0.5ex}}
\stackon[1pt]{#1}{\scalebox{-1}{\tmpbox}}
}

\DeclareMathOperator{\Tr}{Tr}

\DeclareMathOperator{\SAT}{-SAT}
\DeclareMathOperator{\type}{Type}
\DeclareMathOperator{\wt}{wt}

\newcommand{\ol}{\overline}
\newcommand{\ul}{\underline}
\newcommand{\wc}{\widecheck}

\newcommand{\EE}{\mathbb{E}}

\newcommand{\RR}{\mathbb{R}}
\newcommand{\PP}{\mathbb{P}}
\newcommand{\NN}{\mathbb{N}}

\newcommand{\cA}{\mathcal{A}}
\newcommand{\cB}{\mathcal{B}}

\newcommand{\cF}{\mathcal{F}}

\newcommand{\mcG}{\mathcal{G}}

\newcommand{\cH}{\mathcal{H}}
\newcommand{\cI}{\mathcal{I}}
\newcommand{\mcI}{\mathcal{I}}
\newcommand{\mcH}{\mathcal{H}}

\newcommand{\cG}{\mathcal{G}}
\newcommand{\cP}{\mathcal{P}}

\newcommand{\cS}{\mathcal{S}}

\title{Nearly all $k$-SAT functions are unate}
\author[Balogh]{J\'ozsef Balogh}
\author[Dong]{Dingding Dong}
\author[Lidick\'y]{Bernard Lidick\'y}
\author[Mani]{Nitya Mani}
\author[Zhao]{Yufei Zhao}

\thanks{Balogh was supported in part by NSF grants DMS-1764123 and RTG DMS-1937241, FRG DMS-2152488, the Arnold O. Beckman Research Award (UIUC Campus Research Board RB 22000), the Langan Scholar Fund (UIUC)}
\thanks{Lidick\'y was supported in part by NSF grant FRG DMS-2152490 and Scott Hanna professorship.}
\thanks{Mani was supported by the NSF Graduate Research Fellowship Program and a Hertz Graduate Fellowship.}
\thanks{Zhao was supported in part by NSF CAREER award DMS-2044606, a Sloan Research Fellowship, and the MIT Solomon Buchsbaum Fund.}

\address{Balogh: Department of Mathematics, University of Illinois at Urbana-Champaign, Champaign, IL 61801, USA}
\email{jobal@illinois.edu}

\address{Dong: Department of Mathematics, Harvard University, Cambridge, MA 02138, USA}
\email{ddong@math.harvard.edu}

\address{Lidick\'y: Department of Mathematics, Iowa State University, Ames, IA 50010, USA}
\email{lidicky@iastate.edu}

\address{Mani \& Zhao: Department of Mathematics, Massachusetts Institute of Technology, Cambridge, MA 02139, USA}
\email{\{nmani,yufeiz\}@mit.edu}

\begin{document}

\begin{abstract}
We prove that $1-o(1)$ fraction of all $k$-SAT functions on $n$ Boolean variables are unate (i.e., monotone after first negating some variables), for any fixed positive integer $k$ and as $n \to \infty$. This resolves a conjecture by Bollob\'as, Brightwell, and Leader from 2003.
\end{abstract}

\maketitle

\section{Introduction}

\subsection{Background}
We study the following basic question on Boolean functions:
\begin{quote}
How many $k$-SAT functions on $n$ Boolean variables are there? 
What does a typical such function look like?
\end{quote}
This question was first studied by Bollob\'as, Brightwell, and Leader~\cite{BBL03}. 
We focus on the regime where $k$ is fixed and $n \to \infty$.
We count $k$-SAT functions in their disjunctive normal form (DNF). It would be an equivalent problem to enumerate $k$-SAT functions in their conjunctive normal form (CNF) since the negation of a DNF is a CNF. 
For our purpose, a \emph{$k$-SAT function} on $n$ Boolean variables is a function $f \colon \{0,1\}^n \to \{0,1\}$ of the form
\[
f(x_1, \dots, x_n) = C_1 \lor C_2 \lor  \cdots \lor C_m,
\]
where each $C_i$ has the form $z_1 \land \cdots \land z_k$ with $z_1, \dots, z_k \in \{x_1, \ol {x_1}, \dots, x_n, \ol {x_n}\}$.
Here we call $x_1, \dots, x_n$ the \emph{variables}.
Each of $x_i$ and $\ol{x_i}$ is  a \emph{literal} (\emph{positive literal} and \emph{negative literal}, respectively). Each \emph{clause} $C_i$ is a conjunction (``and'') of $k$ literals.
We further restrict that every clause uses $k$ distinct variables (e.g., both $x_1 \land x_2$ and $x_1 \land x_2 \land \ol{x_2}$ are invalid 3-SAT clauses).
This restriction does not lose any generality (the first example can be replaced by $(x_1 \land x_2 \land x_3) \lor (x_1 \land x_2 \land \ol{x_3})$ and the second example is a clause that is never satisfied and so can be deleted).
To simplify notation, we will drop the ``and'' symbol $\land$ when writing a clause.
A \emph{formula} is a set of clauses. 
For instance, the $2$-SAT formula 
$(x_1 \land x_2) \lor  (\ol{x_1}\land x_3) \lor (x_3 \land x_4)$ is written as 
$\{x_1x_2, \ol{x_1} x_3, x_3 x_4\}$. 
Every $k$-SAT function has a $k$-SAT formula, but different $k$-SAT formulae may correspond to the same $k$-SAT function.

Given the importance of $k$-SAT functions, we study the typical structure of such a function.
While the total number of functions $\{0,1\}^n \to \{0,1\}$ is $2^{2^n}$, the number of $k$-SAT functions is significantly smaller. 
As an easy upper bound, since there are $2^k \binom{n}{k}$ possible clauses, the number of $k$-SAT formulae is $2^{2^k \binom{n}{k}}$.
Meanwhile, the number of $k$-SAT functions is at most $2^{2^k \binom{n}{k}}$,
which is significantly smaller than $2^{2^n}$ for a fixed $k$ and large $n$.
The actual number of $k$-SAT functions turns out to be considerably smaller than even this upper bound.

A $k$-SAT formula is \emph{monotone} if it only uses positive literals.
A $k$-SAT function is monotone if it has a monotone $k$-SAT formula.
There are $\binom{n}{k}$ possible monotone clauses, and every monotone $k$-SAT formula produces a unique monotone $k$-SAT function, hence there are $2^{\binom{n}{k}}$ monotone $k$-SAT functions.

A $k$-SAT function or formula is \emph{unate} if it is monotone after replacing some variables with their negations (e.g., $\{\ol{x_1}x_2, \ol{x_1}x_2,x_2x_3\}$ is unate but $\{\ol{x_1}x_2, x_1x_2, x_2x_3\}$ is not).
In other words, a formula is unate if each variable $x_i$ shows up only as a positive literal $x_i$ or only as a negative literal $\ol{x_i}$. A $k$-SAT function is unate if it admits a unate formula.
The number of unate $k$-SAT formulae that use all $n$ variables is at least
\[
2^n \left(2^{\binom{n}{k}} - n 2^{\binom{n-1}{k}} \right) = (1+o(1)) 2^{n + \binom{n}{k}},
\]
for fixed $k$ as $n \to \infty$. Indeed, for each variable $x_i$, there are $2^n$ choices as to whether to use it as a positive literal or a negative literal and there are at least $2^{\binom{n}{k}} - n 2^{\binom{n-1}{k}}$ monotone formulae that use all $n$ variables. 
All unate $k$-SAT formulae represent distinct functions, and thus the number of unate $k$-SAT functions on $n$ variables is at least $(1+o(1)) 2^{n + \binom{n}{k}}$.

Bollob\'as, Brightwell, and Leader~\cite{BBL03} conjectured in 2003 that nearly all $k$-SAT functions are unate. We prove their conjecture.

\begin{theorem} \label{thm:BBL}
Fix $k \ge 2$. 
The number of $k$-SAT functions on $n$ Boolean variables is $(1+o(1)) 2^{n + \binom{n}{k}}$.
Equivalently: a $1-o(1)$ fraction of all $k$-SAT functions on $n$ variables are unate.
\end{theorem}

Bollob\'as, Brightwell, and Leader also proposed a weaker version of this conjecture, namely that the number of $k$-SAT functions on $n$ Boolean variables is $2^{(1+o(1))\binom{n}{k}}$, and established this weaker conjecture for $k=2$.
\cref{thm:BBL} was proved for $k=2$ by Allen \cite{All07} and for $k=3$ by Ilinca and Kahn~\cite{IK12}. The proofs in \cite{BBL03,All07} for $k=2$ used the Szemer\'edi Regularity Lemma~\cite{Szemer},  (\cite{IK09} gave an alternate regularity-free proof for $k=2$), whereas the proof for $k=3$ \cite{IK12}  used the hypergraph regularity method. 

To prove \cref{thm:BBL}, we use hypergraph container and stability arguments to reduce~\cref{thm:BBL} to a Tur\'an density problem for partially directed hypergraphs.
We then solve this Tur\'an problem.
The hypergraph container approach is simpler to carry out than the earlier regularity approaches. One of our insights is distilling the ``right'' Tur\'an problem to be solved.

\subsection*{Acknowledgments} 
This paper integrates results  presented at SODA~`22 by a subset of the authors (Dong, Mani, and Zhao), as well as STOC~`23 by all the authors. The reduction to the Tur\'an problem, encompassing all sections other than \cref{s:turan}, originates from the work presented at SODA~`22. Zhao presented these results at an Oberwolfach workshop in April 2022, sparking the collaboration that culminated in proving the main conjecture by solving the Tur\'an problem (\cref{s:turan}).

Zhao expresses gratitude to Jeff Kahn, who first introduced this problem during Zhao's graduate studies, for his encouragement to pursue this challenge.

\subsection{Strategy}
Our setup is closest to that of Ilinca and Kahn~\cite{IK12} who solved the problem for $3$-SAT functions. Instead of counting functions, we count minimal formulae. In \cite{IK12}, the term ``non-redundant'' is used for what we are calling ``minimal.''

\begin{definition}[Minimal formula] \label{def:minimal}
A formula $G = \{ C_1, C_2, \dots , C_l\}$ on variables $\{x_1, \ldots, x_n\}$ is \emph{minimal} if deleting any clause from $G$ changes the resulting Boolean function. 
Equivalently, $G$ is minimal if for each clause $C_i \in G$, there is a \emph{witness} assignment $\mathsf w \in \{0, 1\}^n$ that satisfies $C_i$ but no other clause $C_j \in G$, $j \neq i$. 
\end{definition}

\begin{example}
The $2$-SAT formula $\{wx,wy,x\ol{z},\ol{y}z\}$ is not minimal since it is impossible to satisfy only $wx$ and no other clause. 
Indeed, if we attempt to only satisfy $wx$, we must assign $w=1$ and $x=1$; 
to avoid satisfying $wy$ and $x\ol{z}$, we must assign $y=0$ and $z=1$, which would then lead to the final clause $\ol{y}z$ being satisfied.
\end{example}

Every $k$-SAT function can be expressed as a minimal formula, but possibly in more than one way. We upper bound the number of minimal $k$-SAT formulae, which in turn upper bounds the number of $k$-SAT functions. Thus, to prove \cref{thm:BBL}, we actually prove the following stronger statement.

\begin{theorem} \label{thm:min-formula}
Let $k \ge 2$ be a fixed positive integer.
The number of minimal $k$-SAT formulae on $n$ Boolean variables is $(1+o(1)) 2^{n + \binom{n}{k}}$.
Equivalently: a $1-o(1)$ fraction of all minimal $k$-SAT formulae on $n$ variables are unate.
\end{theorem}

To upper bound the number of minimal $k$-SAT formulae, we identify a fixed finite set $\cB$ of ``forbidden'' non-minimal formulae.
We then upper bound the number of \emph{$\cB$-free formulae}, which are formulae (not necessarily minimal) not containing any subformula isomorphic to any element of $\cB$.

The problem of counting $\cB$-free formulae is analogous to counting $F$-free graphs on $n$ vertices.
A classic result by Erd\H{o}s, Kleitman, and Rothschild~\cite{EKR76} shows that almost all triangle-free graphs are bipartite. 
Erd\H{o}s, Frankl, and R\"odl~\cite{EFR86} generalized this result and showed that for a fixed graph $F$, the number of $n$-vertex $F$-free graphs is $2^{\operatorname{ex}(n, F) + o(n^2)}$. 
The latter result was initially proved using the Szemer\'edi Regularity Lemma~\cite{Szemer}, i.e., the graph regularity method.
It can also be proved as a quick application of the more recently developed hypergraph container method~\cite{BMS15,ST15} (see also the survey \cite{BMS18}).

Here is a quick sketch of how to use the container method to enumerate $n$-vertex triangle-free graphs. 
By an application of the hypergraph container theorem, there is a collection $\cG$ of $n$-vertex graphs (``containers'')
each with at most  $(1/4 + o(1))n^2$ edges, 
such that $\abs{\cG} = 2^{o(n^2)}$
and every $n$-vertex triangle-free graph is a subgraph of some container $G \in \cG$.
It then follows that the number of $n$-vertex triangle-free graphs is at most $\sum_{G \in \cG} 2^{e(G)} = 2^{(1/4 + o(1)) n^2}$.
See \cref{ss:containerlemma} for more details.

As in the case of counting $F$-free graphs, the enumeration of $\cB$-free formulae reduces to a certain extremal problem. 
Here the situation diverges from the graph theory setting.
For counting triangle-free graphs, the container method reduces the enumeration problem to Mantel's theorem on the maximum number of edges in a triangle-free graph.
It is more intricate to count $k$-SAT functions. 
We identify a structure that we call \emph{partially directed hypergraphs} to cleanly capture the essence of the problem.
Also, unlike counting $F$-free graphs, our extremal problem is not the direct analogue of the maximum size of an $F$-free graph.

We now begin to explain this reduction.

\begin{definition}[Simple formula] \label{def:simple}
A $k$-SAT formula is \emph{simple} if no two clauses use exactly the same set of $k$ variables.
\end{definition}

\begin{example}
The $3$-SAT formula $\{wxy, wx\ol{z},xyz\}$ is simple, whereas $\{xyz,xy\ol{z}\}$ is not.
\end{example}

Every minimal $k$-SAT formula on $n$ variables can be made simple by removing $o(n^k)$ clauses (\cref{prop:nearly-simple}). 
We then observe the following steps:
\begin{align*}
\# \set{\text{$k$-SAT functions}}
&\le 
\# \set{\text{minimal $k$-SAT formulae}}
\\
&\sim  
\# \set{\text{simple minimal $k$-SAT formulae}}
\\
&\le 
\#\set{\text{simple $\cB$-free $k$-SAT formulae}}.
\end{align*}

Recall that $\cB$ is some collection of non-minimal formulae.
Let us focus on enumerating simple $\cB$-free formulae.
The container theorem gives us a collection $\cG$ of $n$-variable ``container'' formulae with $\abs{\cG} = 2^{o(n^k)}$
such that every $n$-variable $\cB$-free formula is a subformula of some container $G \in \cG$. 
Furthermore, every $G \in \cG$ is ``not too large'' in a certain sense that we now describe.
For a container $G \in \cG$, suppose $\alpha \binom{n}{k}$ $k$-subsets of variables support exactly one clause, and $\beta \binom{n}{k}$ $k$-subsets of variables support exactly two clauses.
Let us ignore $k$-subsets of variables supporting more than two clauses since there are $o(n^k)$ of them (\cref{lem:supersat-2}).
For such a container $G$, how many simple subformulae of $G$ are there? For each $k$-subset of variables, if $G$ has exactly $m$ clauses supported on these $k$ variables, then to include at most one of the clauses in our choice of a simple subformula of $G$, there are exactly $m+1$ choices.
So the number of simple subformulae of $G$ is $2^{\alpha\binom{n}{k}} 3^{\beta\binom{n}{k}} 2^{o(n^k)}$, which is at most $2^{(1+o(1))\binom{n}{k}}$ if and only if $\alpha + (\log_2 3) \beta \le 1 +o(1)$. 
So we want $G \in \cG$ to be ``not too large'' in the sense that $\alpha + (\log_2 3) \beta \le 1 +o(1)$.

Even if we wish to only enumerate simple formulae, we need to solve an extremal problem for containers $G$, which are non-simple formulae. 
This is not an inadequacy of our approach, but an essential feature of the problem. (See a related discussion at the end of \cite[Section 2]{BBL03}.)
This is where our situation differs from enumerating $F$-free graphs, which reduces to an extremal problem on $F$-free graphs instead of some more complicated object.

The extremal problem for $\cB$-free formulae appears to be quite complex. 
To solve the problem, we forget just the right amount of  information of a $k$-SAT formula in order to reduce the object to a \emph{partially directed hypergraph}, which we will introduce shortly. 
The reduction leads to a clean formulation of an extremal problem. 
If we had been left with a more complicated extremal problem, we might not have been successful in solving it.

The container approach as described so far leads to an upper bound of the form $2^{(1+o(1))\binom{n}{k}}$. 
Obtaining the tighter asymptotic $(1+o(1))2^{n + \binom{n}{k}}$ of \cref{thm:BBL} requires additional stability arguments.
By analogy, the container method sketched earlier shows that there are $2^{(1/4+o(1))n^2}$ triangle-free graphs on $n$ vertices.
To obtain the Erd\H{o}s--Kleitman--Rothschild result in its full strength, namely that $1-o(1)$ fraction of all triangle-free graphs are bipartite, one needs to apply additional stability arguments on top of the container method, as done in \cite{BBCLMS17}, and later strengthened in \cite{BS20}.

Recall that  \cref{thm:BBL} is equivalent to the statement that a $1-o(1)$ fraction of all $k$-SAT functions are unate.
Far-from-unate formulae can be effectively handled by the container argument.
The final piece is to show is that almost all nearly unate formulae are actually unate. 
We extend the stability argument from \cite[Section 8]{IK12}, though we have to overcome further technical challenges.
This stability argument is quite technical and it occupies the bulk of the latter half of this paper starting at \cref{s:stability}.
The arguments up to the end of \cref{s:firstcount} prove a weaker version of \cref{thm:BBL} that the number of $k$-SAT functions on $n$ Boolean variables is $2^{(1+o(1)) \binom{n}{k}}$, and these sections contain most of the novel ideas of our work.

\subsection{Partially directed hypergraphs}

A \emph{partially directed graph} (also known as a \emph{mixed graph}) is formed by taking a graph and orienting a subset of its edges.
Here to \emph{orient} an edge means to choose one of two directions for the edge.
Every edge is either \emph{directed} or \emph{undirected}. An example is illustrated below.

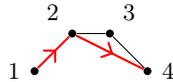
\begin{center}
    \begin{tikzpicture}[
    scale=.5,
    font=\footnotesize,
    v/.style = {circle, fill, inner sep = 0pt, minimum size = 3pt},
    e/.style = {thick, red, postaction={decorate,decoration={
        markings,
        mark=at position .5*\pgfdecoratedpathlength+1.5pt 
        with {\arrow{angle 90}}
      }}}]
    \node[v,label=left:1] (1) at (0,0) {};
    \node[v,label=above left:2] (2) at (1,1) {};
    \node[v,label=above right:3] (3) at (2,1) {};
    \node[v,label=right:4] (4) at (3,0) {};
    \draw (2)--(3)--(4);
    \draw[e] (2)--(4);
    \draw[e] (1)--(2);
     \end{tikzpicture}
\end{center}

\noindent Given a pair of partially directed graphs $\vec H$ and $\vec G$, we say that $\vec H$ is a \emph{subgraph} of $\vec G$ if one can obtain $\vec H$ from $\vec G$ by a combination of (1) removing vertices, (2) removing edges, and (3) removing the orientation of some edges. 

The following partially directed graph plays a special role:
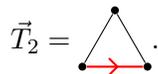
\[
\vec T_2 = 
\begin{tikzpicture}[
    scale=.5,
    baseline={([yshift=-.8ex]current bounding box.center)},
    font=\footnotesize,
    v/.style = {circle, fill, inner sep = 0pt, minimum size = 3pt},
    e/.style = {thick, red, postaction={decorate,decoration={
        markings,
        mark=at position .5*\pgfdecoratedpathlength+1.5pt 
        with {\arrow{angle 90}}
      }}}]
    \node[v] (1) at (90:1) {};
    \node[v] (2) at (210:1) {};
    \node[v] (3) at (-30:1) {};
    \draw (3)--(1)--(2);
    \draw[e] (2)--(3);
    \end{tikzpicture}.
\]
Below, the left graph contains $\vec T_2$ as a subgraph, and the right does not contain $\vec T_2$ as a subgraph.

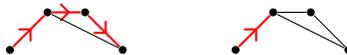
\begin{center}
    \begin{tikzpicture}[
    scale=.5,
    font=\footnotesize,
    v/.style = {circle, fill, inner sep = 0pt, minimum size = 3pt},
    e/.style = {thick, red, postaction={decorate,decoration={
        markings,
        mark=at position .5*\pgfdecoratedpathlength+1.5pt 
        with {\arrow{angle 90}}
      }}}]

    \begin{scope}
        \node[v] (1) at (0,0) {};
        \node[v] (2) at (1,1) {};
        \node[v] (3) at (2,1) {};
        \node[v] (4) at (3,0) {};
        \draw[e] (2)--(3);
        \draw[e] (3)--(4);
        \draw (2)--(4);
        \draw[e] (1)--(2);
    \end{scope}
    
    \begin{scope}[shift={(6,0)}]
        \node[v] (1) at (0,0) {};
        \node[v] (2) at (1,1) {};
        \node[v] (3) at (2,1) {};
        \node[v] (4) at (3,0) {};
        \draw (2)--(3)--(4);
        \draw (2)--(4);
        \draw[e] (1)--(2);
    \end{scope}

    \end{tikzpicture}
\end{center}

The following statement implies \cref{thm:BBL} for $k=2$. (See \cref{thm:strong-count} for the full statement of the implication.)

\begin{theorem}
There is some constant $\theta > \log_2 3$ such that every $n$-vertex partially directed graph with $\alpha \binom{n}{2}$ undirected edges and $\beta \binom{n}{2}$ directed edges and not containing $\vec T_2$ as a subgraph satisfies
\[
\alpha + \theta \beta \le 1 + o_{n\to\infty}(1).
\]
\end{theorem}

Now let us generalize the problem to hypergraphs.
While there are many possible notions of a directed hypergraph, the relevant notion for us is the one where a \emph{directed edge} is an edge with a \emph{pointed} vertex on the edge.
A partially directed hypergraph is formed from a hypergraph by directing some of its edges and leaving others intact.

\begin{definition}[$k$-PDG]
A \emph{partially directed $k$-graph}, abbreviated \emph{$k$-PDG}, is given by a set $V$ of $n$ vertices, and, for each unordered $k$-element subset $\{v_1, \dots, v_k\}$, one of the following possibilities:
\begin{enumerate}
    \item[(a)] No edge with vertices $\{v_1, \dots, v_k\}$, or
    \item[(b)] An \emph{undirected edge} with vertices  $\{v_1, \dots, v_k\}$, or
    \item[(c)] A \emph{directed edge} using vertices  $\{v_1, \dots, v_k\}$ along with a choice of some $v_i \in \{v_1, \dots, v_k\}$;  we say that the edge is \emph{directed} (or \emph{pointed}) toward $v_i$. We notate such a directed edge by $v_1 \cdots \wc v_i \cdots v_k$.
\end{enumerate}
In particular, no two edges (whether directed or undirected) of a $k$-PDG can use the exact same set of $k$ vertices.
We usually denote $k$-PDGs by simply writing their edge sets.
\end{definition}

\begin{example}
This is the edge set of a $3$-PDG: $\{123,12\wc4, 134,235,\wc245\}$.
Non-examples of edge sets of $3$-PDGs include $\{123,12\wc 3\}$ (two edges using the same triple of vertices), $\{1\wc2\wc3\}$ (not a valid edge), and $\{123,12\}$ (not 3-uniform).
\end{example}

The following $k$-PDG plays a special role.

\begin{definition}[$\vec T_k$] \label{def:Tk}
For each $k \ge 3$, define $\vec T_k$ to be the $k$-PDG obtained by starting with $\vec T_2$ defined earlier and then adding $k-2$ common vertices to all three edges.
\end{definition}

There are some examples:
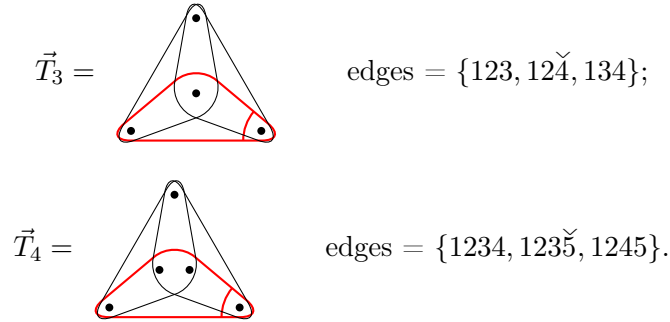
\[\vec T_3
=
    \begin{tikzpicture}[
    font=\footnotesize,
    baseline={([yshift=-.8ex]current bounding box.center)},
    v/.style = {circle, fill, inner sep = 0pt, minimum size = 3pt},
    e/.style = {thick, red, postaction={decorate,decoration={
        markings,
        mark=at position .5*\pgfdecoratedpathlength+1.5pt 
        with {\arrow{angle 90}}
      }}}]
	\node[v] (0) at (0:0) {};
    \node[v] (1) at (90:1) {};
    \node[v] (2) at (210:1) {};
    \node[v] (3) at (-30:1) {};
	
	\draw[red, thick, rounded corners=10pt] (90:.4) -- +(-40:1.6) coordinate (a) -- +(-140:1.6)--cycle;
	\draw[red,  thick] ($ (a) + (-180:.6) $) arc (180:140:.6);

	\draw[rounded corners=10pt,rotate=120] (90:.4) -- +(-40:1.6) coordinate (a) -- +(-140:1.6)--cycle;
	\draw[rounded corners=10pt,rotate=-120] (90:.4) -- +(-40:1.6) coordinate (a) -- +(-140:1.6)--cycle;
    \end{tikzpicture}
    \qquad 
    \text{edges = } 
    \{ 123, 12\wc 4, 134 \};\]
\[
\vec T_4
=
    \begin{tikzpicture}[
    font=\footnotesize,
    baseline={([yshift=-.8ex]current bounding box.center)},
    v/.style = {circle, fill, inner sep = 0pt, minimum size = 3pt},
    e/.style = {thick, red, postaction={decorate,decoration={
        markings,
        mark=at position .5*\pgfdecoratedpathlength+1.5pt 
        with {\arrow{angle 90}}
      }}}]
	\node[v] (0) at (0:.2) {};
    \node[v] (0') at (180:.2) {};
    
    \node[v] (1) at (90:1) {};
    \node[v] (2) at (210:1) {};
    \node[v] (3) at (-30:1) {};
	
	\draw[red, thick, rounded corners=10pt] (90:.4) -- +(-40:1.6) coordinate (a) -- +(-140:1.6)--cycle;
	\draw[red,  thick] ($ (a) + (-180:.6) $) arc (180:140:.6);

	\draw[rounded corners=10pt,rotate=120] (90:.4) -- +(-40:1.6) coordinate (a) -- +(-140:1.6)--cycle;
	\draw[rounded corners=10pt,rotate=-120] (90:.4) -- +(-40:1.6) coordinate (a) -- +(-140:1.6)--cycle;
    \end{tikzpicture}
    \qquad 
    \text{edges = } 
    \{ 1234, 123\wc 5, 1245 \}.
\]

\noindent Given a $k$-PDG $\vec{H}$, we write
\begin{align*}
e_u(\vec H) &= \text{the number of undirected edges in $\vec H$, and}
\\
e_d(\vec H) &= \text{the number of directed edges in $\vec H$.}	
\end{align*}
Denote the edge densities by
\[
\alpha(\vec{H}) = \frac{e_u(\vec H)}{\binom{n}{k}}
\quad \text{and} \quad 
\beta(\vec{H}) = \frac{e_d(\vec H)}{\binom{n}{k}}.
\]

\begin{definition}[Subgraph of $k$-PDG] \label{def:subgraph}
We say that a $k$-PDG $\vec{F}$ is a \emph{subgraph} of another $k$-PDG $\vec{H}$ if we can obtain a $k$-PDG isomorphic to $\vec{F}$ starting from $\vec{H}$ by some sequence of actions of the following types: (a) deleting an edge, (b) deleting a vertex, and (c) forgetting the direction of an edge.
\end{definition}

\begin{definition}[$\vec F$-free $k$-PDG]
Given $k$-PDG $\vec F$, we say that a $k$-PDG $\vec{H}$ is \emph{$\vec F$-free} if $\vec{H}$ does not contain $\vec F$ as a subgraph.
\end{definition}

\begin{definition}[Tur\'an density for $k$-PDGs]\label{d:densitypdg}
Given a $k$-PDG $\vec F$, let $\pi(\vec F, \theta)$ be the smallest real number such that every $n$-vertex $\vec F$-free $k$-PDG $\vec{H}$ satisfies
\[
\alpha(\vec H) +\theta \beta(\vec H)  \le \pi(\vec F, \theta) + o_{n\to\infty}(1).
\]
\end{definition}

The Tur\'an problem for $\vec T_k$-free $k$-PDGs will play a central role, as $\vec T_k$ encodes non-minimal formulae that turn out to be critical for the $k$-SAT enumeration problem. 
Note that for every $k$ and $\theta$ we always have
\[
\pi(\vec T_k, \theta) \ge 1,
\]
since the complete undirected $k$-PDG is $\vec T_k$-free. (This is related to the fact that every monotone formula is minimal.) 
Also, note the trivial inequality $\pi(\vec T_k, \theta) \le \pi(\vec T_k,\theta')$ for $\theta < \theta'$.

\subsection{Results}  

We prove the following two results, which together imply \cref{thm:BBL}.

\begin{theorem}[Reduction to Tur\'an problem] \label{thm:strong-count}
Let $k \ge 2$. If $\pi(\vec T_k, \theta) = 1$ for some $\theta > \log_2 3$, then~\cref{thm:BBL} and \cref{thm:min-formula} are true for this $k$.
\end{theorem}

\begin{theorem}[Tur\'an problem] \label{t:pi}
For every $k \ge 2$, there is some constant $\theta > \log_2 3$ such that
\[
\pi(\vec T_k, \theta) = 1.
\]
\end{theorem}

The proof of \cref{t:pi} is  in \cref{s:turan}. The proof of \cref{thm:strong-count} occupies the rest of the paper starting in \cref{s:hyptools}.
On our way to proving \cref{thm:strong-count}, by the end of \cref{s:firstcount}, we establish the following weaker conclusion.

\begin{theorem}[Weak reduction] \label{thm:weak-count}
    For $k \ge 2$, the number of $n$-variable $k$-SAT functions is at most 
    \[ 2^{(\pi(\vec T_k, \log_2 3)+o(1))\binom{n}{k}}.
    \]
\end{theorem}

\cref{t:pi,thm:weak-count} already show that the number of $k$-SAT functions is $2^{(1+o(1))\binom{n}{k}}$ for all $k \ge 2$. It is easier to establish this weaker conclusion since it does not require the technical stability arguments that occupy the second half of the paper starting with \cref{s:stability}. Most of the new ideas in our paper are contained in the proof of this weaker conclusion.

We prove \cref{t:pi} by reducing it to a simple though mysterious looking inequality. 
While this inequality is easy to check by hand, its discovery required a computer search, drawing inspiration from Razborov's innovative flag algebras method \cite{Raz07}. Notably, our approach diverges from many other flag algebra applications commonly seen in extremal graph and hypergraph problems, as our proof certificate is remarkably concise and readily verifiable without computer assistance.

Despite the extensive study of hypergraph Tur\'an density problems, definitive solutions have remained elusive, with only particular cases being resolved. For further insights into hypergraph Tur\'an problems, we recommend Keevash's comprehensive survey \cite{Kee11}.

\begin{remark}[Earlier version] \label{rem:earlier}
An earlier version of this work, available at \texttt{arXiv:2107.09233v3}, contains additional arguments and results that do not appear in this version since they are no longer necessary for proving \cref{thm:BBL}. Here we summarize some these elements that might have potential future value.

First, we identify a family of forbidden $k$-PDGs that completely characterize the asymptotic enumeration of minimal $k$-SAT formulae. 
This means that the enumeration of $k$-SAT functions is equivalent to the Tur\'an problem with this larger family of forbidden subgraphs. 
It turns out, as we now know from~\cref{t:pi}, that forbidding $\vec T_k$ alone is sufficient, although this was not clear a priori.
It is possible that the additional forbidden configurations (that are no longer mentioned in this paper) may play a critical role in further extensions of this problem discussed next.

We also deleted a proof of~\cref{t:pi} for $k=4$ via an interesting argument that applied a recent result by F\"uredi and Maleki~\cite{FM17}, as well as an appendix (by Nitya Mani and Edward Yu) that gives a computer verification of the $k=5$ case of~\cref{t:pi}.	
\end{remark}

\subsection{Further directions} 

As suggested by Bollob\'as, Brightwell, and Leader~\cite{BBL03}, this investigation opens doors to a \emph{theory of random $k$-SAT functions}. 
For example, \cref{thm:min-formula} implies that a typical $k$-SAT function admits a unique minimal $k$-SAT formula, and furthermore this formula has $(1/2+o(1))\binom{n}{k}$ clauses. 
Note that our model is very different from that of \emph{random $k$-SAT formulae} where clauses are added at random (e.g., the recent breakthrough on the satisfiability conjecture~\cite{DSS22}).
Rather, \cref{thm:min-formula} concerns a random $k$-SAT formula \emph{conditioned} on mini\-mality. 
In this light, our results are analogous to the theory of dense Erd\H{o}s--R\'enyi random graphs.
In random graph theory, understanding the typical structure of a dense triangle-free graph has led to many fruitful research directions concerning typical structures and thresholds, and inspired important combinatorial techniques such as the hypergraph container method.
One such direction concerns the typical structure of sparser triangle-free graphs. 
Another direction concerns the ty\-pical structure of a graph without a $t$-clique, when $t$ increases with $n$.
We refer to \cite{BMSW16,BS20} and their references for discussions of the extensive work on these and related problems.
It is natural to ask similar questions about $k$-SAT functions.

\begin{question}
Fix $k$.
For a given $m = m(n)$, what is the typical structure of a $k$-SAT function on $n$ Boolean variables that can be expressed as a formula with $m$ clauses?
What is the typical structure of a minimal $k$-SAT formula with $n$ Boolean variables and $m$ clauses? What is the threshold $m_0$ above which a typical such formula is unate?
\end{question}

Another interesting further direction concerns large $k$.
Bollob\'as and Brightwell~\cite{BB03} studied the $k$-SAT function enumeration problem for large $k = k(n)$. They conjectured that as long as $k \le (1/2-c) n$ for some constant $c > 0$, the number of $k$-SAT functions on $n$-variables is $2^{(1+o(1))\binom{n}{k}}$. 
In that paper \cite{BB03}, they proved bounds on the number of $k$-SAT functions on $n$ variables in the regime $k \ge n/2$, where a completely different asymptotic behavior arises.

\begin{question}
	For a given $k = k(n)$, how many $k$-SAT functions on $n$ Boolean variables are there?
	What about minimal $k$-SAT formula?
	What is the threshold $k_0(n)$ below which it is typically unate?
\end{question}

A significant bottleneck here is the quantitative bounds of the container lemma. 
Even with a more efficient container lemma~\cite{BS20}, it seems difficult to analyze what happens when $k$ grows faster than $\log n$, let alone linear in $n$.

\subsection{Outline}

The rest of the paper can be divided into three parts.

\smallskip

\emph{Part I. Tur\'an problem.}

In~\cref{s:turan}, we prove~\cref{t:pi} (the Tur\'an density problem).

\smallskip

\emph{Part II. Exponential asymptotics.}

In~\cref{s:hyptools}, we recall some tools for hypergraphs, including the hypergraph container theorem. 
In \cref{s:firstcount}, we prove \cref{thm:weak-count} giving an upper bound on the number of $k$-SAT functions. 

\smallskip

\emph{Part III. Stability.}

Finally, we prove \cref{thm:strong-count} via a technical stability argument.
In~\cref{s:stability}, by a more careful analysis of the containers, we reduce the problem to showing that there are negligibly many non-unate but nearly unate minimal formulae.
This final claim is then established in \cref{s:istar} by extending the arguments in \cite[Section 8]{IK12}.

\section{A Tur\'an density problem}\label{s:turan}
In this section, we prove \cref{t:pi}.
We quote the following lemma.

\begin{lemma}[F\"uredi {\cite[Lemma 2.1]{Fur92}}] \label{lem:furedi}
Given a graph $G = (V, E)$, let $G^2$ be the graph on vertex set $V$ with $xy \in E(G^2)$ if and only if there is some $z \in V$ with $xz, yz \in E$. Then, for any graph $G$,
$$e(G^2) \ge e(G) - \lfloor n/2 \rfloor.$$
\end{lemma}

\begin{theorem} \label{thm:T2}
	$\pi(\vec T_2, 2) = 1$. 
\end{theorem}

\begin{proof} 
Let $\vec H$ be $\vec T_2$-free 2-PDG on $n$ vertices. 
Let $\alpha = \alpha(\vec H)$ and $\beta = \beta(\vec H)$.
Let $G$ be the underlying simple graph of $\vec{H}$ (the graph obtained from $\vec H$ by replacing directed edges with simple edges). Then $e(G) = (\alpha + \beta){n \choose 2}$. 
We say that edge $xy \in E(G)$ is a \emph{triangular edge} if there exists some $z \in V(G)$ so that $xz,yz \in E(G)$.    
Only non-triangular edges in $G$ can arise as directed edges in $\vec{H}$, i.e., $|E(G)-E(G^2)|\ge \beta {n\choose 2}$. We further note that all triangular edges in $G$ are in $G^2$, i.e., a pair of vertices either forms an edge in $G^2$, or in $G-G^2$, or in neither of them.  Thus by~\cref{lem:furedi},
$${n\choose 2} \ge e(G^2) + \beta{n \choose 2} \ge e(G) + \beta{n \choose 2} - O(n)
= (\alpha + 2\beta){n \choose 2} - O(n),$$
which implies that $\alpha + 2 \beta \le 1 + o(1)$.  
\end{proof}

\begin{theorem}
	$\pi(\vec T_3, 5/3) = 1$.
\end{theorem}

\begin{proof}
Let $\vec H$ be a $\vec T_3$-free $3$-PDG on vertex set $V$.
Let $\alpha = \alpha(\vec H)$ and $\beta = \beta(\vec H)$.
For each vertex $v \in V$, 
the \emph{link of $v$}, denoted $\vec{H}_v$, is the $2$-PDG with vertex set $V \setminus \{v\}$ and with edges obtained by taking all edges of $\vec H$ containing $v$ and then deleting the vertex $v$ from each edge. We preserve all the direction data, except that every edge directed toward $v$ in $\vec H$ becomes an undirected edge in $\vec H_v$.
Write
$\alpha_v = \alpha(\vec H_v)$
and
$\beta_v = \beta(\vec H_v)$.
By linearity of expectation over a uniform random vertex $v \in V$, $\EE_v \alpha_v = \alpha + \beta/3$ and 
$\EE_v \beta_v = 2\beta/3$.
For every $v \in V$, $\vec{H}_v$ is $\vec T_2$-free (for otherwise combining with $v$ forms a $\vec T_3$), 
so $\alpha_v + 2 \beta_v \le 1+o(1)$ by \cref{thm:T2}.
Hence
$
\alpha + 5\beta/3 = 
\EE_v[\alpha_v + 2\beta_v]
\le 1 + o(1)$.
\end{proof}

We only consider $k\ge 4$ from now on.
Here is the key lemma.

\begin{lemma} \label{lem:flag-ineq}
	Let $k \ge 4$ be a positive integer.
	There exist $\theta > \log_23$ and $a,b \in \RR$ (depending on $k$) such that  the following holds.
	Suppose $\vec F$ is a $\vec T_k$-free $k$-PDG on $k+1$ vertices.
	Let $x_1, \dots, x_{k-1}, y, z$ be a permutation of vertices of $\vec F$ chosen uniformly at random. 
 
	Define $x := x_1\ldots x_{k-1}$, let $x y$ denote the event that $x_1 \ldots x_{k-1} y$ forms an undirected edge in $\vec F$, let
	$x \wc y$ denote the event that $x_1 \ldots x_{k-1} \wc y$ is a directed edge in $\vec F$,
	and let $\ul{x y}$ denote the event that there is no edge with vertices $x_1, \ldots, x_{k-1}, y$ in $\vec F$.
	Then, we have that
	\[
	\PP(x y) + k\theta \PP(x\wc y) + a^2 \PP(\ul{x y} \wedge \ul{x z}) - 2ab \PP(\ul{x y} \wedge x\wc z) + b^2 \PP(x \wc y \wedge x \wc z) \le 1.
	\]
\end{lemma}
\begin{remark}
The proof gives $\theta = 1 + \frac{1}{\sqrt{2}} \ge 1.707 > 1.585 > \log_2 3$, with $
a = \frac{1}{\sqrt{2}},
\,\,
b = \frac{k(\theta - 1) - 1}{\sqrt{2}}=\frac{k}{2}- \frac{1}{\sqrt{2}}.$ 
\end{remark}

\begin{theorem}There exists $\theta>\log _23$ such that
    $\pi(\vec T_k,\theta)=1$ for all $k\geq 4$.
\end{theorem}

\begin{proof}Let $\vec H$ be an $n$-vertex $\vec T_k$-free $k$-PDG with $\alpha \binom{n}{k}$ undirected edges and $\beta \binom{n}{k}$ directed edges.

Let $x_1,\dots,x_{k-1},y,z$ be vertices of $\vec H$ chosen without replacement uniformly at random. 
Applying \cref{lem:flag-ineq} (by first conditioning on the set of $k+1$ randomly selected vertices) and using the notation of the lemma,
\[
\PP(x y) + k\theta \PP(x\wc y) + a^2 \PP(\ul{x y} \wedge \ul{x z}) - 2ab \PP(\ul{x y} \wedge x\wc z) + b^2 \PP(x \wc y \wedge x \wc z) \le 1.
\]
Note that $\PP(x y)=\alpha$ and $k \PP(x\wc y)=\beta$.
It remains to show that
\[
a^2 \PP(\ul{x y} \wedge \ul{x z}) - 2ab \PP(\ul{x y} \wedge x\wc z) + b^2 \PP(x \wc y \wedge x \wc z) \ge -o(1).
\]
We will show that this inequality holds for every fixed choice of $x$.
Conditioned on $x$, vertices $y$ and $z$ are uniformly chosen  without replacement in $V(\vec F) \setminus \set{x_1, \dots, x_k}$. 
When $n$ is large, this is not much different from choosing  with replacement, in which case $y$ and $z$ would be conditionally independent, given $x$.
In particular, 
$\PP(\ul{x y} \wedge \ul{x z} | x) = \PP(\ul{x y}| x)^2 + o(1)$, 
and $\PP(\ul{x y} \wedge x\wc z | x) = \PP(\ul{xy} | x) \PP(x \wc y | x) + o(1)$,
and
$\PP(x \wc y \wedge x \wc z | x) = \PP(x \wc y | x)^2 + o(1)$.
Then we can prove the above displayed inequality, conditioned on any $x$, by observing that
\begin{align*}
 &a^2 \PP(\ul{x y} \wedge \ul{x z} | x) - 2ab \PP(\ul{x y} \wedge x\wc z | x) + b^2 \PP(x \wc y, x \wc z | x) 
\\
&= a^2 \PP(\ul{x y}| x)^2  - 2ab  \PP(\ul{xy} | x) \PP(x \wc y | x) +  b^2 \PP(x \wc y | x)^2 - o(1)
\\
&= (a \PP(\ul{x y}| x) - b \PP(x \wc y | x))^2 - o(1) \ge -o(1). \qedhere 
\end{align*}
\end{proof}

Finally, it remains to prove \cref{lem:flag-ineq}.

\begin{proof}[Proof of \cref{lem:flag-ineq}]
Let $\vec H$ be a $\vec T_k$-free $k$-PDG on $k+1$ vertices.
Construct the following associated digraph $D$ on the same vertex set as $\vec H$:
\begin{enumerate}
	\item For every directed edge in $\vec H$ missing vertex $i$ and directed towards vertex $j$, add the directed edge $i \to j$ in $D$.
	\item For every undirected edge in $\vec H$ missing vertex $i$, add the loop $i \to i$ in $D$.
\end{enumerate}
Notice that every vertex in $D$ has out-degree at most $1$. Moreover, since $\vec H$ does not contain $\vec T_k$ as a subgraph, $D$ is free of the following forbidden pattern:
\begin{itemize}
	\item \emph{Forbidden pattern:} $i_1 \to i_2 \to *$ and $i_3 \to *$ for three distinct vertices $i_1,i_2,i_3$, where $*$ can be any vertex (the two $*$'s do not have to be the same).
\end{itemize}
Indeed, if $D$ contains the above pattern, then $\vec H$ has an edge missing $i_2$, an edge missing $i_3$, and an edge containing both $i_2$ and $i_3$ that is directed at $i_2$. These three edges contain $\vec T_k$ as a subgraph.

It is easy to deduce the following exhaustive classification of all digraphs $D$ on $k+1$ vertices  with out-degree at most $1$ at every vertex and without the above forbidden pattern.
\begin{enumerate}
	\item [(A)] $i \to j \to i$ for distinct $i,j$. No other edges.
	\item [(B)] $i \to j \to j$ for distinct $i,j$. No other edges.
	\item [(C)] $i_1 \to i_2 \to i_3$ for distinct $i_1,i_2,i_3$. No other edges.
	\item [(D)] For every edge of the form $i\to j$ with $i,j$ distinct, $j$ has zero out-degree. (In other words, $D$ is vertex-disjoint union of self-loops and stars in which all edges are directed toward the central vertex.)
\end{enumerate}

Write $i\rightarrow\varnothing$ to denote  that $i$ has out-degree 0 in $D$. The inequality in \cref{lem:flag-ineq} then translates to
\begin{multline}
	 \label{eq:d-ineq}
	\PP(z \to z) 
	+ k \theta \cdot \PP(z \to y) + a^2 \PP(y \to \varnothing \wedge z \to \varnothing)
	\\
	- 2ab \PP(y \to z \wedge z \to \varnothing)
	+ b^2 \PP(y \to z \wedge z \to y)
	\le 1.
\end{multline}

Let $u$ and $d$ denote the number of undirected and directed edges in $\vec H$, respectively. We calculate the terms in  inequality \cref{eq:d-ineq}:
\begin{align*}
\PP(z \to z) &= \frac{u}{k+1}, \quad 
\PP(z \to y) = \frac{d}{k(k+1)},
\\
\PP(y \to \varnothing \wedge z \to \varnothing) &= \frac{(k+1-u-d)(k-u-d)}{(k+1)k},
\\
\PP(y \to z \wedge z \to \varnothing)
&= \begin{cases}
 	0 & \text{in cases (A), (B)}, \\
 	\frac{1}{(k+1)k} & \text{in case (C)}, \\
 	\frac{d}{(k+1)k} & \text{in case (D)},
 \end{cases}
\\
\PP(y \to z \wedge z \to y) &= \begin{cases}
 \frac{2}{(k+1)k} & \text{in case (A)}, \\
 0 & \text{in cases (B), (C), (D)}. \\ 	
 \end{cases}
\end{align*}
Then \cref{eq:d-ineq} reduces to simpler inequalities in each of the four cases:
\begin{itemize}
	\item [(A)]  When our digraph is of the form $i \to j \to i$ for distinct $i,j$ with no other edges,~\cref{eq:d-ineq} simplifies to 
	\begin{equation}\label{eq:a}
	\frac{2\theta}{k+1} + \frac{(k-1)(k-2)}{(k+1)k} a^2 + \frac{2}{(k+1)k} b^2 \le 1.
	\end{equation}
	\item [(B)] Since $\theta \ge 1$, the associated inequality in this case is implied by (A); we note the inequality below, but do not need to consider this case separately:
	\[
	\frac{1}{k+1} + \frac{\theta}{k+1} + \frac{(k-1)(k-2)}{(k+1)k} a^2 \le 1.
	\]
	\item [(C)] The inequality in this case is also implied by (A) (we will end up choosing $a,b > 0$) so we also need not consider it separately:
	\[
	\frac{2\theta}{k+1}
	+ \frac{(k-1)(k-2)}{k(k+1)}a^2
	- \frac{2}{k(k+1)}ab\le 1.
	\]
	\item[(D)] In this case, if $d = 0$ then $u \le k+1$, and if $d \ge 1$ then $u + d \le k$. The inequality reduces to
	\begin{equation}\label{eq:c}
	\frac{u}{k+1} + \frac{\theta d}{k+1} + \frac{(k+1-u-d)(k-u-d)}{(k+1)k} a^2  - \frac{2d}{(k+1)k} ab \le 1.
	\end{equation}
\end{itemize}
Thus it remains to show that we can choose parameters $\theta > \log_2 3$ and $a, b > 0$ to make the inequalities implied by (A) and (D) true. We choose
\[
\theta = 1 + \frac{1}{\sqrt{2}} > 1.707,
\, \quad \text{ satisfying }\quad  \ 2\theta^2 - 4\theta + 1 = 0,\]
\[
a = \frac{1}{\sqrt{2}},
\,\,\quad \quad
b = \frac{k(\theta - 1) - 1}{\sqrt{2}}=\frac{k}{2}- \frac{1}{\sqrt{2}}.
\]

We can verify~\cref{eq:a}, establishing the desired inequality for case (A) (and thus cases (B) and (C)) via a direct substitution:

\begin{align*}
 &\frac{2\theta}{k+1} + \frac{(k-1)(k-2)}{(k+1)k} a^2 + \frac{2}{(k+1)k} b^2 - 1
\\
&= 	
\frac{2\theta}{k+1} + \frac{(k-1)(k-2)}{2(k+1)k} + \frac{(k(\theta-1)-1)^2}{(k+1)k} - 1
\\&=
\frac{4k\theta + (k-1)(k-2) + 2(k(\theta-1)-1)^2 - 2(k+1)k}{2(k+1)k}
\\
&=
\frac{4 - k + k^2 ( 2\theta^2 - 4\theta + 1)}{2(k+1)k}
=
\frac{4 - k}{2(k+1)k} \le 0.
\end{align*}

To verify~\cref{eq:c} for (D), 
observe that the left-hand side of \cref{eq:c} is non-decreasing in $u$. Indeed, replacing $u$ by $u+1$ increases the expression by at least 
\[
\frac{1}{k+1} - \frac{(k+1)k - k(k-1)}{(k+1)k} a^2	
= \frac{1}{k+1} - \frac{2}{k+1} a^2 = 0.
\]
Consequently, it remains to verify case (D) when $u$ is as large as possible, meaning $u+d \in \set{k,k+1}$, which makes the third term in~\cref{eq:c} zero. The left-hand side of~\cref{eq:c} becomes
\begin{align*}
\frac{u}{k+1} + \frac{\theta d}{k+1}  - \frac{2d}{(k+1)k} ab.
\end{align*}
If $d = 0$, then $ u = k+1$, and the inequality clearly holds.
Otherwise, $d + u = k$.
Since the above expression is linear in $u$ (or equivalently, linear in $d$), the maximum is attained at one of endpoints $(u,d) = (k,0)$ or $(0,k)$. The only nontrivial situation to check is $(u,d) = (0,k)$, in which case the above expression is
\[
\frac{k\theta}{k+1} - \frac{2}{k+1} ab 
= \frac{k\theta}{k+1} - \frac{k(\theta - 1) -1}{k+1} = 1.
\]
Therefore, \cref{eq:d-ineq} always holds.
\end{proof}

Our work also motivates the following natural Tur\'an problem for $k$-PDGs.

\begin{problem}
For each $k \ge 2$, what is the largest $\theta$ such that $\pi(\vec T_k, \theta) = 1$?
\end{problem}

The following $k$-PDG is $\vec T_k$-free: 
partitioning the $n$ vertices as $A \cup B$ with $\abs{A} = n/k$ and $\abs{B} = (k-1)n/k$, and add all directed edges with $k-1$ vertices in $B$ and directed toward a vertex in $A$.
This $k$-PDG has $\abs{A}\binom{\abs{B}}{k-1} = (n/k) \binom{(1-1/k)n}{k-1} = (1+o(1))(1-1/k)^{k-1}\binom{n}{k}$ edges. So $\pi(\vec T_k, \theta) > 1$ whenever $\theta >  (1-1/k)^{-k+1}$.
We conjecture that this construction is optimal in the following sense.

\begin{conjecture} \label{conj:pi-opt}
$\pi(\vec T_k, (1-1/k)^{-k+1}) = 1$ for every $k \ge 2$.
\end{conjecture}

\cref{thm:T2} confirms this conjecture for $k=2$. 
The conjectured optimal construction is reminiscent of the solution to the hypergraph Tur\'an problem corresponding to $3$-graphs with independent  neighborhoods by F\"uredi, Pikhurko, and Simonovits~\cite{FPS05}.

\section{Hypergraph tools}\label{s:hyptools}

\emph{Remarks on asymptotic notation.} We view $k\ge 2$ as a fixed constant and do not explicitly mention dependence on $k$ in the hidden constants. We use standard asymptotic notation including $O(\cdot)$ and $\Omega(\cdot)$. We use subscripts to denote that the hidden constant factors may depend on the the subscripted parameters (along with $k$, which is omitted). For example, $f(n) = \Omega_B(g(n))$ for some $g> 0$ means that there is some constant $c = c(B, k) > 0$ such that $f(n) \ge cg(n)$ for all $n$.

\subsection{Densities of blowups} \label{sec:blowup}

Given a $k$-SAT formula $G$, its \emph{$b$-blowup}, denoted $G[b]$, is obtained from $G$ by replacing each variable by $b$ identical duplicates. For example, for the $3$-SAT formula $G = \{x_1x_2x_3,x_1\ol x_3x_4\}$, we have
\begin{multline*}    
G[2]=\left\{x_1x_2x_3, x_1x_2x_3',x_1x_2'x_3,x_1x_2'x_3',
x_1'x_2x_3, x_1'x_2x_3',x_1'x_2'x_3,x_1'x_2'x_3',\right. \\
\left.x_1\ol {x_3}x_4,x_1\ol {x_3}x_4', x_1\ol {x_3'}x_4,x_1\ol {x_3'}x_4',
x_1'\ol {x_3}x_4,x_1'\ol {x_3}x_4', x_1'\ol {x_3'}x_4,x_1'\ol {x_3'}x_4'\right\}.
\end{multline*}

By a standard Cauchy--Schwarz argument, we have the following.

\begin{lemma}\label{lem:formula-blowup}
For every $\epsilon > 0$ and $k$-SAT formula $F$, there is some $\delta > 0$ and $n_0$ such that for all $n > n_0$, if a $k$-SAT formula on $n$ variables has at least $\epsilon n^{v(F)}$ copies of $F$, then it has at least $\delta n^{2v(F)}$ copies of $F[2]$.
\end{lemma}

\begin{proof}
    Suppose $G$ is a $k$-SAT formula on variables $x_1,\dots,x_n$ that has at least $\epsilon n^{v(F)}$ copies of $F$. Let $\ell=v(F)$. Define 
    $Z=\{x_1,\dots,x_n,\ol{x_1},\dots,\ol{x_n}\}$ to be the set of $2n$ literals, and $V$ to be the disjoint union of $\ell$ copies of $Z$, i.e., $V=Z_1\cup\dots\cup Z_\ell$ where every $Z_j=\{x_1^{(j)},\dots,x_n^{(j)},\ol{x_1}^{(j)},\dots,\ol{x_n}^{(j)}\}$.
    
    Let $W=\{w_1,\dots,w_\ell\}$ be the set of variables that appeared in $F$. Construct the associated $\ell$-graph $A_{F}$ on $V$ as follows: for every subformula $F'\subseteq G$ isomorphic to $F$ with isomorphism $\phi:W\to Z$, add edge $\phi(w_1)^{(1)}\cdots\phi(w_\ell)^{(\ell)} $ to $A_{F}$. Then $A_{F}$ is an $\ell$-partite $\ell$-graph with at least $\epsilon n^{\ell}$ edges.

    For $U_1\subseteq Z_1,\dots,U_\ell\subseteq Z_\ell$, we say that $A_F[U_1,\dots,U_\ell]$ is complete if $u_1\cdots u_\ell\in E(A_F)$ for all $u_1\in U_1,\dots, u_\ell\in U_\ell$. For every $i\in[\ell]$, let $v_i,v_i'$ be two independent random literals in $Z_i$.
    For $j=0,\dots,\ell$, define the probability $$p_j=\PP[A_F[\{v_1,v_1'\},\dots, \{v_j,v_j'\}, \{v_{j+1}\},\dots, \{v_\ell\}]\text{ is complete}].$$
    By the Cauchy--Schwarz inequality, for each $j=0, \dots, \ell-1$,
    \begin{align*}
        p_{j+1}&=\EE_{\substack{v_1,v_1',\dots,v_{j+1},v_{j+1}'\\v_{j+2},\dots,v_\ell}}\prod_{u_1\in\{v_1,v_1'\}}\cdots \prod_{u_{j+1}\in\{v_{j+1},v_{j+1}'\}} \mathbf 1[u_1\cdots u_{j+1}v_{j+2}\cdots v_\ell\in E(A_F)]\\
        &=\EE_{\substack{v_1,v_1',\dots,v_{j},v_{j}'\\v_{j+2},\dots,v_\ell}}(\EE_{v_{j+1}}\prod_{u_1\in\{v_1,v_1'\}}\cdots \prod_{u_{j}\in\{v_{j},v_{j}'\}}\mathbf 1[u_1\cdots u_{j}v_{j+1}\cdots v_\ell\in E(A_F)])^2\\
        &\geq (\EE_{\substack{v_1,v_1',\dots,v_{j},v_{j}'\\v_{j+1},\dots,v_\ell}}\prod_{u_1\in\{v_1,v_1'\}}\cdots \prod_{u_{j}\in\{v_{j},v_{j}'\}}\mathbf 1[u_1\cdots u_{j}v_{j+1}\cdots v_\ell\in E(A_F)])^2=p_j^2.
    \end{align*}
    Since $p_0=|E(A_F)|/n^\ell \ge \eps$,
    we have $p_{\ell}\ge \eps^{2^\ell}$.
    Therefore, the number of copies of $F[2]$ in $G$ is at least $p_{\ell} n^{2\ell} - O_\ell(n^{2\ell-1})$ (the error term accounts for non-injective maps), which is at least $\eps^{2^\ell} n^{2\ell}/2$ for sufficiently large $n$.
\end{proof}

\subsection{Kruskal--Katona theorem}

We need the following special case of the Kruskal--Katona theorem~\cite{Kru63,Kat68}.
Here a \emph{simplex} in a $k$-graph is a clique on $k+1$ vertices.

\begin{theorem} \label{thm:kruskalkatona}
A $k$-uniform hypergraph on $n$-vertices with at most $\beta n^k/k!$ edges contains at most $\beta^{\frac{k+1}{k}} n^{k+1}/(k+1)!$ simplices.
\end{theorem}

Here is a quick proof of this statement using Shearer's entropy inequality~\cite{CGFS86} (see also \cite{Kah01}).

\begin{proof}
Let $(X_1, \dots, X_{k+1})$ be the vertices of a uniformly chosen simplex in the $k$-graph, with the $k+1$ vertices permutated uniformly at random.
By Shearer's inequality, letting $\vec X_{-i}$ denote $(X_1, \dots, X_{k+1})$ with the $i$-th coordinate removed, we have
\begin{align*}
k \log((k+1)!\text{\#simplices}) 
&= k H(X_1, \dots, X_{k+1}) 
\\
&\le H(\vec X_{-1}) + \cdots +  H(\vec X_{-(k+1)}) 
\\
&= (k+1) H(X_1, \dots, X_{k})
\le (k+1) \log (k!\text{\#edges}). 
\end{align*}
Since $k!\text{\#edges} \le \beta n^k$, we have $(k+1)!\text{\#simplices} \le \beta^{\frac{k+1}{k}}n^{k+1}$.
\end{proof}

\subsection{Hypergraph containers}\label{ss:containerlemma}

For motivation, let us first recall a container theorem for triangle-free graphs \cite[Theorem 2.1]{BMS18}.

\begin{theorem}[Containers for triangle-free graphs] \label{thm:container-triangle-free}
For every $\epsilon > 0$, there exists $C>0$ such that for every $n$, there exists a set $\cG$ of graphs on $n$ labeled vertices, such that
\begin{enumerate}[(a)]
	\item every triangle-free graph on $n$ labeled vertices is contained in some $G \in \cG$, and
	\item each $G \in \cG$ contains at most $\epsilon n^3$ triangles, and
	\item $\abs{\cG} \le n^{C n^{3/2}}$.
\end{enumerate}
\end{theorem}
This theorem outputs a collection of $2^{o(n^2)}$ \emph{containers} $\mcG$, such that each triangle-free graph is contained within some container $G \in \mcG$ and each container $G$ is \emph{almost triangle-free.}

We now state and prove an analogous result for formulae avoiding a given finite set of subformulae. For any formula $B$, we let $v(B)$ denote the number of variables that appear in some clause in $B$.
Recall that we view a formula as a set of clauses, so $\abs{B}$ denotes the number of clauses. We formally define the notions of a subformula and a formula being ``free'' of subformulae below.
\begin{definition}[Subformula]
A \emph{subformula} of a formula $G$ is a subset of clauses of $G$. 
We say that two formulae are \emph{isomorphic} if one can be obtained from the other by relabeling variables. 
Given another formula $F$, a \emph{copy} of $F$ in $G$ is a subformula of $G$ that is isomorphic to $F$.
Given a set $\cB$ of $k$-SAT formulae, we say that $G$ is \emph{$\cB$-free} if $G$ has no copy of $B$ for every $B \in \cB$.
\end{definition}

We write
\[
    m(\cB)=\max_{B\in\cB}|B|.
\]

We shall apply the hypergraph container theorem, proved independently by 
Balogh, Morris, and Samotij \cite{BMS15} 
and 
Saxton and Thomason \cite{ST15}
to show the following.

\begin{theorem}\label{thm:gencontub}
Let $\cB$ be a finite set of simple $k$-SAT formulae.
For every $\delta > 0$, there exists $C = C(\cB, \delta) > 0$ such that for every $n$, there exists a collection $\cG=\cG(\cB,\delta)$ of formulae on variable set $X = \{x_1, \dots, x_n\}$ such that:
\begin{enumerate}
    \item[(a)] Every $\cB$-free $k$-SAT formula with variables in $X$ is a subformula of some $G \in \cG$, and
    \item[(b)] For every $G \in \cG$ and $B \in \cB$, $G$ has at most $\delta n^{v(B)}$ copies of $B$, and
    \item[(c)] $\abs{\cG} \le n^{Cn^{k-1/(m(\cB)-1)}}$.
\end{enumerate}
\end{theorem}

In the rest of this section, we deduce the above claim from the more general hypergraph container theorem. To state the general result, we introduce some notation for an  $\ell$-uniform hypergraph $\cH$. Let $\mathcal I(\cH)$ denote the collection of independent vertex sets in $\cH$. For every $j\in\{1,\dots,\ell\}$, let $\Delta_j(\cH)$ denote the maximum $j$-codegree of vertices in $\cH$, i.e., the maximum number of edges containing a given $j$-vertex subset in $\cH$. 
Let $\cP(V)$ denote the family of all subsets of $V$.
Let ${V\choose \leq s}$ denote the family of those subsets of $V$ having size at most $s$.
For $\cA \subseteq \cP(V)$, write $\ol\cA = \cP(V) \setminus \cA$.

We apply the following version of the
hypergraph container theorem as stated in \cite[Theorem 2.2]{BMS15}.
A similar result was proved in \cite{ST15}.

\begin{theorem}[Hypergraph container theorem]\label{thm:hc}
For every $\ell\in\NN$ and $c,\epsilon>0$, there exists $C_0=C_0(\ell,c,\epsilon)>0$ such that the following holds. Let $\cH$ be an $\ell$-uniform hypergraph and $\cA\subseteq \cP(V(\cH))$ be an increasing family of sets such that for all $A\in\cA$, we have  $|A|\geq \epsilon v(\cH)$ and $e(\mcH[A])\geq \epsilon e(\mcH)$. Suppose  there exists $p\in(0,1)$ such that, for every $j\in\{1,\dots,\ell\}$, we have
$$\Delta_j(\cH)\leq c\cdot p^{j-1}\cdot\frac{e(\cH)}{v(\cH)}.$$ 
Then there exists a family $\mathcal S\subseteq {V(\cH)\choose \leq C_0\cdot p\cdot v(\cH)}$ and functions $f:\mathcal S\rightarrow\ol{\cA}$, $g:\cI(\cH)\rightarrow\cS$ such that for every $I\in \mathcal I(\mcH)$, we have:
\begin{enumerate}
    \item[(i)] $g(I)\subseteq I$;
    \item[(ii)] $I\setminus g(I)\subseteq f(g(I))$.
\end{enumerate}
\end{theorem}

As motivation for our deduction of~\cref{thm:gencontub} using \cref{thm:hc}, we first review why~\cref{thm:container-triangle-free} follows from \cref{thm:hc}.
To count triangle-free graphs, we define a $3$-uniform hypergraph $\mcH$ encoding triangles, i.e., where $V(\mcH) = {[n] \choose 2}$ and $E(\mcH)$ comprise those triples of edges that form triangles. The triangle-free graphs on $n$ vertices then are in one-to-one correspondence to independent sets on $\mcH$ and after verifying the codegree conditions with $p = n^{-1/2}$,~\cref{thm:hc} implies~\cref{thm:container-triangle-free}.

\cref{thm:gencontub} follows from \cref{thm:hc} in a similar way. For a finite set $\cB$ of simple $k$-SAT formulae, we define a hypergraph $\cH$ encoding formulae in $\cB$, so that $\cB$-free formulae become independent sets in $\cH$. We then use \cref{thm:hc} to deduce \cref{thm:gencontub} after verifying the required codegree conditions.

With this preparation, we now prove \cref{thm:gencontub}.

\begin{proof}[Proof of~\cref{thm:gencontub} using~\cref{thm:hc}]
Fix a finite set $\cB$ of simple $k$-SAT formulae. 
Let $m = m(\cB) = \max_{B \in \cB} |B|$ and $v = \max_{B\in\mathcal B}v(B)$.

Let $\mcH$ be an $m$-uniform hypergraph, where $V(\mcH)$ is the set of all $k$-SAT clauses on variable set $\{x_1,\dots,x_n\}$ (so $
v(\mcH) = 2^k{n\choose k}$), and
$$E(\mcH)=\left\{G\in \binom{V(\mcH)}{m}\mid G\text{ is simple; $v(G)\leq v$; $G$ contains a copy of some }B\in \cB\right\}.$$
Observe that, for any formula $B\in\cB$, there are at least $\binom{n}{v(B)}$ different formulae supported on $x_1,\dots,x_n$ isomorphic to $B$. Thus, with $v=\max_{B\in\cB}v(B)$, we have the lower bound $e(\mathcal H)\geq \binom{n}{v}$. On the other hand, since $e(\mathcal H)$ is at most the number of simple formulae on  $\leq v$ vertices, we have $e(\mathcal H)= O(n^v)$.

Using the above bounds on $e(\mcH)$ and $v(\mcH)$, we obtain the following codegree bounds for $\mcH$.
\begin{itemize}
    \item Since every clause has $k$ variables, we have
    $$\Delta_1(\mcH) = O_{\cB}(n^{v - k}) =O_\cB\left( e(\mcH)/v(\mcH)\right).$$
    \item Take $j \in \{2, 3, \ldots, m\}$. Note that any $G \in E(\mcH)$ arises from a \emph{simple} formula, and thus any $j \ge 2$ clauses from such $G$ must include at least $k + 1$ variables. This implies that 
    $$\Delta_j(\mcH) = O_{\cB} (n^{v - k - 1})=O_\cB\left( p^{j - 1}e(\mcH)/v(\mcH)\right)$$
    with $p = n^{-1/(m-1)}$.
\end{itemize}

Fix $\delta>0$ and choose $\epsilon=\epsilon(\cB,\delta)>0$ sufficiently smaller than $\delta$ to be determined later. Define the family
\[
\cA=\{A\subseteq V(\cH):e(\cH[A])\geq \epsilon e(\cH)\}.
\]
For every $A\in\cA$, since $\epsilon e(\cH)\leq e(\cH[A])\leq \Delta_1(\cH)|A|$ with  $\Delta_1(\cH)=O_\cB(n^{v-k})$, we have $|A|\geq  \epsilon e(\cH)/\Delta_1(\cH)=\Omega_\cB(\epsilon)  v(\cH)$.
 
Applying~\cref{thm:hc}, we know that there exists $C=C(\cB,\delta)$ and a family
$$\mathcal S\subseteq {V(\mcH) \choose \leq C p v(\mcH)}$$
with functions $f:\mathcal S\rightarrow\ol{\cA}$, $g:\cI(\cH)\rightarrow\cS$ such that for every $I\in \mathcal I(\mcH)$, we have:
\begin{enumerate}
    \item[(i)] $g(I)\subseteq I$;
    \item[(ii)] $I\setminus g(I)\subseteq f(g(I))$.
\end{enumerate}

Let $\cG=\cG(\cB,\delta)=\{S \cup f(S):S\in\cS\}$. We show that $\cG$ meets all the conditions described in~\cref{thm:gencontub}.

(a) \emph{Every $\cB$-free formula with $n$ variables is a subformula of some $G \in \cG$.} 

Every $\cB$-free formula is an independent set $I$ in $\mcH$. Since $I\setminus g(I)\subseteq f(g(I))$ and $g(I)\in\cS$, we have $I\subseteq g(I)\cup f(g(I)) \in\cG$.

(b) \emph{For every $G \in \cG$ and $B \in \cB$, $G$ has at most $\delta n^{v(B)}$ copies of $B$.} 

Every $G \in \cG$ is of the form $G=S \cup f(S)$ for some $S\in\cS$. Since $|S|\leq  C p v(\mcH)$, the number of edges incident to $S\setminus f(S)$ is at most (recalling $\Delta_1(\cH)=O_\cB(e(\cH)/v(\cH))$ and $p = n^{-1/(m-1)}$)
$$
    C p v(\mcH) \Delta_1(\cH)
    =O_{\cB,\delta}( p e(\mathcal H))
    =O_{\cB,\delta}( n^{v-1/(m-1)}).
$$
Meanwhile, since $f(S)\in\ol{\cA}$, we have $e(\cH[f(S)])<\epsilon  e(\cH)$.
Hence for every $G=S\cup f(S) \in \cG$, we have
$$e(\cH[G])\leq \epsilon e(\cH)+O_{\cB,\delta}( n^{v-1/(m-1)})\leq 2\epsilon e(\cH).$$
Let $\epsilon=\epsilon(\cB,\delta)$ be sufficiently small such that, for any $B\in\cB$, 
if $G$ has more than $\delta n^{v(B)}$ copies of $B$, then $e(\mcH[G])>2\epsilon  e(\cH)$. Notice that such an $\epsilon$ always exists: if $v(B)=v$, then the statement is straightforward; if $v(B)<v$ and $G$ has more than $\delta n^{v(B)}$ copies of $B$, we get $\Omega(\delta)n^v$ edges in $\cH[G]$ by adding $v-v(B)$ free variables to every copy of $B$. Therefore, for every $G\in\cG$, since $e(\cH[G])\leq 2\epsilon  e(\cH)$, we know that  $G$ has at most $\delta n^{v(B)}$ copies of $B$ for every $B\in\cB$.

(c) \emph{$\abs{\cG} \le n^{C'n^{k-1/(m-1)}}$ for some $C'=C'(\cB,\delta)$.} 

Since $\abs{\cG}\leq \abs{\cS}$, we have (recalling $p = n^{-1/(m-1)}$)
\[
\abs{\cG} 
\le
\binom{v(\mcH)}{\le C pv(\mcH)} 
\le n^{O_{\cB,\delta}(n^{k-1/(m-1)})}. 
\qedhere 
\]
\end{proof}

\section{Weak upper bound on the number of $k$-SAT functions}\label{s:firstcount}
In this section, we prove~\cref{thm:weak-count}. 
Recall that a $k$-SAT formula is \emph{simple} (\cref{def:simple}) if it does not contain a pair of clauses with the same associated variable set.
We begin by observing that every minimal formula is nearly simple.
It remains to upper bound the number of \emph{simple minimal formulae}.
We will select some fixed finite set $\cB$ of non-minimal formulae (the forbidden subformulae).
Then every simple minimal formula is $\cB$-free.
To upper bound the number of simple $\cB$-formula, we proceed via the following steps.

\begin{enumerate}[(a)]
    \item \textbf{Container formulae.} \cref{thm:gencontub} gives a collection of $2^{o(n^k)}$ \emph{container formulae} such that (i) each simple $\cB$-free formulae is a subformula of some container, and (ii) every container has few copies of elements in $\cB$.
    \item \textbf{$k$-PDGs.} Each container formula is nearly ``semisimple'' (\cref{def:semisimple}). We convert each container formula $G$ to a $k$-PDG $\vec{F}$.
    If $\vec{F}$ had many copies of $\vec{T}_k$, then $G$ would have many copies of some element of $\cB$,
    so $G$ must be nearly $\vec{T}_k$-free.
    The Tur\'an bound, \cref{t:pi}, then implies an upper bound on the number of simple subformulae of $G$.
\end{enumerate}

\subsection{Minimal formulae are nearly simple}\label{ss:minimalsimple}

\begin{lemma}\label{lem:pairclausesnonmin}
The $2$-blowup of a pair of clauses on the same set of $k$ variables is non-minimal.
\end{lemma}

For example, the lemma tells us that the $2$-blowup of $\{xyz,xy\ol z\}$ is non-minimal.

\begin{proof}
Up to relabeling and/or negating some of the variables, we can write the pair of clauses as $v_1 v_2 \cdots v_k$ and $\overline{v_1} \cdots \ol {v_j}v_{j+1}\cdots v_k$ for some $1\leq j\leq k$. 
Then its $2$-blowup has clauses (among others)
$$v_1 v_2 \cdots v_k,v_1' v_2 v_3 \cdots v_k, v_1 v_2' v_3 \cdots v_k, \ldots, v_1 \cdots v_{j-1}v_j'v_{j+1}\cdots v_k, \overline{v_1'} \cdots \ol {v_j'}v_{j+1}\cdots v_k.$$ 
It is impossible to satisfy the first clause only, since it would involve setting $v_1 = \cdots = v_k = 1$ and $v'_1 = \cdots = v'_j = 0$, which would then satisfy the last clause,
 so the $2$-blowup is non-minimal.
\end{proof}

In the above proof, if $j \ge 2$, then the displayed sequence of clauses is simple. We record this fact here as it will be used later~\cref{ss:supersaturation}.
\begin{lemma}\label{lem:parallel}
Given a pair of clauses on the  set of $k$ variables such that the clauses differ by at least two variable negations, then its $2$-blowup contains a simple non-minimal subformula. \qed 
\end{lemma}

\begin{proposition}\label{prop:nearly-simple}
Every minimal $k$-SAT formula on $n$ variables can be made simple by deleting $o(n^k)$ clauses.
\end{proposition}

\begin{proof}
Let $G$ be a minimal $k$-SAT formula on $n$ variables.
By \cref{lem:pairclausesnonmin}, $G$ does not contain any 2-blowup of a pair of clauses on the same $k$ variables. 
It follows from \cref{lem:formula-blowup} that the number of $k$-subsets of variables that support at least two clauses in $G$ is $o(n^k)$. Removing all such clauses yields a simple subformula of $G$.
\end{proof}

\subsection{From $k$-SAT to $k$-PDGs}

As hinted by \cref{lem:parallel}, for the container argument, we focus our attention on the following special type of formulae.

\begin{definition}[Semisimple formula] \label{def:semisimple}
A $k$-SAT formula is \emph{semisimple} if for every $k$-subset $S$ of variables, one of the following is true:
\begin{enumerate}
    \item[(a)] There is no clause on $S$, or
    \item[(b)] There is exactly one clause on $S$, or 
    \item[(c)] There are exactly two clauses on $S$ and they differ by exactly one variable negation.
\end{enumerate}
\end{definition}

\begin{example}
The 3-SAT formula $\{ x_1\overline{x}_2x_3, x_1\ol{x}_2\ol{x}_3,x_1x_2\overline{x}_4\}$ is semisimple, whereas the formula $\{ x_1\overline{x}_2x_3, x_1 x_2 \ol{x}_3\}$ is not semisimple.
\end{example}

We define a forgetful map from semisimple formulae to $k$-PDGs.
\begin{definition}[Type map]\label{d:type}
Define
\[
\type\colon \{\text{semisimple $k$-SAT formula on $n$ variables} \} \to \{\text{$k$-PDGs on $n$ vertices}\}
\]
as follows.
Given a semisimple formula $G$, let $\type(G)$ be the $k$-PDG $\vec{F}$ formed as follows:
\begin{itemize}
    \item the vertex set of $\vec{F}$ are the variables of $G$;
    \item for each $k$-subset of variables supporting exactly one clause in $G$, we add an associated undirected edge to $\vec{F}$;
    \item for every $k$-subset of variables supporting exactly two clauses in $G$, we add a directed edge to $\vec{F}$, directed towards the variable where these two clauses differ by negation.
\end{itemize}
\end{definition}

\begin{example}
The semisimple formula
\[
G = \{abc,ab\ol{c}, \ol{a}b \ol{d},bd\ol{e},\ol{b}d\ol{e}\}
\]
has $\type(G)$ the $3$-PDG with vertex set $\{a, b, c, d, e\}$ and edge set $\{ab\wc{c}, abd, \wc{b}de\}$.
\end{example}

\subsection{Forbidden subgraphs and non-minimal formulae}

\begin{proposition}\label{prop:fktypenonmin}
Let $G$ be a semisimple $k$-SAT formula with $\type(G) = \vec T_k$.
Then $G[2]$ has a simple non-minimal subformula.
\end{proposition} 

\begin{example}
We consider a $3\SAT$ example.
An example of a semisimple formula with Type $\vec{H}= (\{1,2,3,4\}, \{123, 13\wc 4, 234\})$ is the collection of $4$-clauses on $4$-variables $G =\{abc, a\ol bd, a\ol b\ol d, \ol bc \ol d\} $. 
The $2$-blowup of $G$ has the following simple non-minimal formula: $B = \{abc, ab'c, a\ol {b'}d, \ol{b'}c\ol{d}\}$. We can observe that $B$ is non-minimal by noting that to only satisfy $abc$, we must assign $a,b,c \mapsto 1$ which forces us to assign $b' \mapsto 0$ to avoid satisfying the second  clause in $B$. However then we must assign $d \mapsto 0$ to avoid satisfying the third clause in $B$ and this combination of assignments satisfies the fourth clause in $B$.
\end{example}

\begin{proof}[Proof of \cref{prop:fktypenonmin}]
Let $G$ be a semisimple formula with $$\type(G) = \vec T_k = (\{1,2, \ldots, k+1\},\{\underbrace{12 3 \cdots k}_{e_a},\, \underbrace{2 3 \cdots k (k+1)}_{e_b},\, \underbrace{1 2 \cdots (k-1) (\widecheck{k+1}}_{e_c})\} ).$$

We suppose that variables $x_1,\dots,x_{k+1}$ in $G$ corresponds to vertices $1,\dots,k+1$ in $\vec T_k$ under the Type map. Without loss of generality, suppose that clause $C_a = x_1 x_2 \cdots x_k$ lies in $G$, which gets mapped to $e_a=1\cdots k$ under the Type map. We also suppose that the blowup $G[2]$ has variable set  $\{x_1, \ldots, x_{k+1}, x_1', \ldots, x_{k+1}'\}$. 

We will construct a simple non-minimal subformula of $G[2]$, that we term $H$. We begin by including the following clauses from $G[2]$ into $H$.
\begin{align*}
C_a &= x_1 x_2 x_3 \cdots x_k, \\
C_a^{(1)} &= x_1' x_2 x_3 \cdots x_k, \\
C_a^{(2)} &= x_1 x_2' x_3 \cdots x_k ,\\
C_a^{(3)} &= x_1 x_2 x_3' \cdots x_k, \\
&\vdots \\
C_a^{(k)} &= x_1 x_2 x_3 \cdots x_k'.
\end{align*}
If $H$ was minimal, then we should be able to find an assignment $\phi$ of variables $\{x_1, \ldots, x_{k+1}, x_1', \ldots, x_{k+1}'\}$ that \emph{only} satisfies $C_a$ among all clauses in $H$. To satisfy $C_a$, this assignment $\phi$ must satisfy $\phi(x_i) = 1$ for all $i \in [k]$. Then, in order to \emph{avoid} satisfying  $C_a^{(i)}$ for any $i \in [k]$,  we must assign $\phi(x_i') = 0$ for all $i \in [k]$.

Suppose clause $C_b\in G$ gets mapped to $e_b$, and clauses $C_c,C_c'\in G$ get mapped to $e_c$ under the Type map. Write
\begin{align*}
    C_b&=z_2\cdots z_kz_{k+1},\\
    C_c&=w_1\cdots w_{k-1}z_{k+1},\\
    C_c'&=w_1\cdots w_{k-1}\ol{z_{k+1}},
\end{align*}
with literals $z_i,w_i\in\{x_i,\ol{x_i}\}$ for every $i$. If $\phi(x_1)=\dots=\phi(x_k)=1$ and $\phi(x_1')=\dots=\phi(x_k')=0$, then no matter which literals the $z_i,w_i$ are equal to, we can always pick clauses $D_b\in C_b[2]$, $D_c\in C_c[2]\cup C_c'[2]$ such that $\phi(D_b)=\phi(z_{k+1})$ and $\phi(D_c)=\phi(\ol{z_{k+1}})$. Here is an illustrating example: if $$C_b=\ol{x_2}x_3x_4,\quad C_c=x_1\ol{x_2}x_4, \quad C_c'=x_1\ol{x_2}\ol{x_4},$$ then we can pick $D_b=\ol{x_2'}x_3x_4\in C_b[2]$ and $D_c=x_1\ol{x_2'}\ol{x_4}$ such that $\phi(x_1)=\phi(\ol{x_2'})=\phi(x_3)=1$, which gives $\phi(D_b)=\phi(x_4)$ and $\phi(D_c)=\phi(\ol {x_4})$. Therefore, $\phi$ must satisfy one of $D_b$ and $D_c$. 

Let $H=\{C_a,C_a^{(1)},\dots,C_a^{(k)},D_b,D_c\}\subseteq G[2]$. Per above, we know that there is no assignment $\phi$ of variables $\{x_1, \ldots, x_{k+1}, x_1', \ldots, x_{k+1}'\}$ that \emph{only} satisfies $C_a$ among all clauses in $H$, so $H$ is non-minimal.
\end{proof}

\subsection{Supersaturation}\label{ss:supersaturation}
For a $k$-SAT formula $G$ on $n$ variables, define $\alpha_i(G)$ by
\[
\alpha_i(G) {n \choose k} = \left| \left\{S \in {V \choose k} \mid  G \text{ has exactly } i \text{ clauses with variables } S \right\} \right|.
\]
For example, $\alpha_2(G){n \choose k}$ is the number of $k$-subsets of variables $S$ such that $G$ has \emph{exactly} $2$ clauses supported on $S$.

Furthermore, we let $\alpha_2'(G){n \choose k}$ be the number of $k$ subsets of variables $S$ where $G$ has exactly $2$ clauses supported on $S$ \emph{and} this pair of clauses supported on $S$ differs in exactly one variable (which is negated in one clause and positive in the other). We also let 
$\alpha_2''(G) := \alpha_2(G) - \alpha_2'(G)$.

\begin{lemma} \label{lem:supersat-1}
For every $\epsilon, \theta > 0$ there exists $\delta > 0$ such that for all sufficiently large $n$, if an $n$-variable $k$-SAT formula $G$ satisfies $\alpha_1(G)+\theta \alpha_2'(G)\geq \pi(\vec T_k,\theta)+\epsilon$, 
then there is some simple non-minimal subformula $B$ with at most $2(k+1)$ variables, such that $G$ contains at least $\delta n^{v(B)}$ copies of $B$ as subformulae.
\end{lemma}

\begin{proof}Fix $\epsilon > 0$ and $\theta>0$. Consider an $n$-variable $k$-SAT formula $G$ such that $\alpha_1(G)+\theta \alpha_2'(G)\geq \pi(\vec T_k,\theta)+\epsilon$. After deleting all clauses of $G$ not accounted for by $\alpha_1(G)$ and $\alpha_2'(G)$, we arrive at a semisimple subformula. We restrict our attention to this subformula and thus assume hereafter without loss of generality that $G$ is semisimple.

Let $\vec H = \type(G)$. Then $\vec H$ has $\alpha_1(G) \binom{n}{k}$ undirected edges and $\alpha_2'(G)\binom{n}{k}$ directed edges.
It follows by the definition of Tur\'an density $\pi(\vec T_k,\theta)$ and a standard supersaturation argument that $\vec H$ contains at least $\Omega_{\eps, \theta}(n^{k+1})$ copies of $\vec T_k$.
So there is some $k$-SAT semisimple formula $G'$, with $\type(G')=\vec T_k$, such that $G$ contains $\Omega_{\eps, \theta}(n^{k+1})$ copies of $G'$ as subformulae.
Then, by \cref{lem:formula-blowup}, $G$ contains $\Omega_{\epsilon, \theta}(n^{2(k+1)})$ copies of $G'[2]$ as subformulae.
By \cref{prop:fktypenonmin}, each copy of $G'[2]$ contains some simple non-minimal subformula.
Thus there is some simple non-minimal subformula $B$ on at most $2(k+1)$ variables that appears $\Omega_{\epsilon, \theta}(n^{v(B)})$ times.
\end{proof}

\begin{lemma} \label{lem:supersat-2}
Let $\epsilon > 0$.
There exists $\delta = \delta(k, \epsilon)>0$ such that for all sufficiently large $n$, if an $n$-variable $k$-SAT formula has at least $\alpha_i(G) \ge \epsilon$ for some $i \ge 3$ or $\alpha_2''(G) \ge \epsilon$, then there is some simple non-minimal subformula $B$ on at most $2k$ vertices such that $G$ contains at least $\delta n^{v(B)}$ copies of $B$ as subformulae.
\end{lemma}

\begin{proof}
Consider an $n$-variable $k$-SAT formula $G$ such that $\alpha_i(G) \ge \epsilon$ for some $i \ge 3$, or $\alpha_2''(G) \ge \epsilon$. By the pigeonhole principle, there is some $2 \le j \le k$ such that $G$ has $\Omega_{\epsilon}(\epsilon n^k)$ 
subformulae of the form $G' = \{v_1 v_2 \cdots v_k$ and $\overline{v_1} \cdots \ol {v_j}v_{j+1}\cdots v_k\}$ where each $v_i$ is some literal in $\{x_i, \ol x_i\}$. In other words, $G$ has $\Omega_{\epsilon}(\epsilon n^k)$ pairs of clauses that are supported on the same variable set but differ by the negation on at least $k$ literals. 
By \cref{lem:formula-blowup}, $G$ contains $\Omega_{\eps}(n^{2k})$ copies of the $2$-blowups of these clause pairs as subformulae. By \cref{lem:parallel}, each copy of $G'[2]$ has a non-minimal simple subformula. The conclusion then follows as earlier.
\end{proof}

\subsection{Applying the container theorem}\label{ss:apply-containers}

Recall the definition of the $\alpha_i$'s from the beginning of \cref{ss:supersaturation}.

\begin{definition} \label{def:weight}
Define the \emph{weight} of a $k$-SAT formula $G$ by
\[
\wt(G)=\alpha_1(G)+\log_23 \cdot \alpha_2(G)+\log_24 \cdot \alpha_3(G)+\dots+\log_2(2^k+1) \cdot \alpha_{2^k}(G).
\]
\end{definition}

Given a $k$-SAT formula $G$ on $n$ variables, the number of simple $k$-SAT subformulae of $G$ is
\begin{equation}\label{eq:wt-count}
\prod_{j =1}^{2^k} (j+1)^{\alpha_j(G) {n \choose k}} = 2^{\wt(G) {n \choose k}},
\end{equation}
since there are $j+1$ ways to choose at most one clause from  $j$ clauses. This motivates our definition of weight.

We are now ready to prove \cref{thm:weak-count} that the number of minimal $k$-SAT formulae on $n$ variables is at most $2^{(\pi(\vec T_k, \log_2 3) + o(1)) \binom{n}{k}}$.

\begin{proof}[Proof of \cref{thm:weak-count}]
Let $\cB$ be the set of all non-minimal simple formulae on at most $2(k+1)$ variables, one for each isomorphism class (so that $\cB$ is finite).
Fix any $\epsilon > 0$.
Pick $\delta > 0$ to be smaller than the minimum of the $\delta$'s chosen in \cref{lem:supersat-1} (with $\theta = \log_2 3$) and \cref{lem:supersat-2}.
Applying the container result, \cref{thm:gencontub}, we find a collection $\cG$ of $n$-variable formulae on $X=\{x_1,\dots,x_n\}$ such that
\begin{enumerate}
    \item [(a)] Every $\cB$-free formula on $X$ is a subformula of some $G \in \cG$, so in particular, every simple minimal formula on $X$ is a subformula of some $G \in \cG$; and
    \item[(b)] For every $G \in \cG$ and $B \in \cB$, $G$ contains at most $\delta n^{v(B)}$ copies of $B$, which implies, by \cref{lem:supersat-1,lem:supersat-2}, that $\alpha_1(G) + \log_2 3 \cdot \alpha_2'(G)< \pi(\vec T_k, \log_2 3) + \epsilon$, $\alpha_2''(G) < \epsilon$ and $\alpha_i(G)< \epsilon$ for all $i \ge 3$; and 
    \item[(c)] $\abs{\cG} \le 2^{o(n^k)}$.
\end{enumerate}

The number of simple minimal $k$-SAT subformula of each $G \in \cG$ is thus at most
\[
2^{\wt (G) \binom{n}{k}} \le 2^{(\pi(\vec T_k, \log_2 3) + O(\epsilon))\binom{n}{k}}.
\]
By taking a union bound over all $G \in \cG$, of which there are at most $2^{o(n^k)}$, and noting that $\epsilon$ can be taken to be arbitrarily small, the number of simple minimal $k$-SAT formula on $n$ variables is $2^{(\pi(\vec T_k, \log_2 3) + o(1))\binom{n}{k}}$. 

Finally, by \cref{prop:nearly-simple},  we can obtain any minimal $k$-SAT formulae by adding $o(n^k)$ clauses to a simple minimal formula. This adds a negligible factor $2^{o(n^k)}$. So the total number of minimal $k$-SAT formulae on $n$ variables is $2^{(\pi(\vec T_k, \log_2 3) + o(1))\binom{n}{k}}$.
\end{proof}

\section{Stability I: near vs.\ far from unate}\label{s:stability}

We prove~\cref{thm:strong-count} in the following two sections. 
Let us recall some definitions. The \emph{variables} are $x_1, \dots, x_n$. Each variable can appear as either a \emph{positive literal} $x_i$ or a \emph{negative literal} $\ol{x_i}$. We say that a clause \emph{uses} a variable $x$ if it contains one of the literals $x$ and $\ol x$.

We usually omit mentioning dependencies on $k$ as we consider it fixed throughout. Recall that a clause is \emph{monotone} if it only uses positive literals.
An assignment $\mathsf w$ of variables is a \emph{witness} for a clause $C$ in a formula if $\mathsf w$ satisfies $C$ and no other clauses. A formula is \emph{minimal} if every clause has a witness.

We next define several useful properties a $k$-SAT formula might have.
\color{black}

\begin{definition}
\label{d:cistar}

Let $\Phi$ be an $n$-variable $k$-SAT formula and let $\zeta, \zeta' > 0$.
\begin{itemize}
    \item $\Phi$ is \emph{$\zeta$-nearly monotone} if it has at most $\zeta n^k$ non-monotone clauses.
    \item $\Phi$ is \emph{$\zeta$-sparsely minimal} if it is minimal and every clause has a witness assignment with $< \zeta n$ variables assigned to $1$. 
    \item $\Phi$ is \textit{positive-dominant}  if for all variables $x$,  $m(x) \ge m(\overline{x})$, where for a literal $v$, $m(v)$ is the number of occurrences of $v$ in $\Phi$.
    \item $\Phi$ is $\zeta'$-\textit{dense} if every variable is used by at least $\zeta' {n - 1 \choose k - 1}$ clauses.
\end{itemize}
Let $\cI(n)$ be the set of minimal $k$-SAT formulae on $n$ variables, and let $\cI^*(n, \zeta, \zeta')$ be the set of minimal $k$-SAT formula that enjoy the above properties with parameters $\zeta, \zeta'$.
\end{definition}

Here is the main result of this section.

\begin{proposition}\label{prop:reducetoboundik}
Suppose $\theta > \log_2 3$ satisfies $\pi(\vec T_k,\theta)=1$.
For every $\zeta>0$, there exists $\epsilon > 0$ (depending only on $k, \theta, \zeta$) such that  for all sufficiently large $n$,
\[
\abs{\cI(n)} 
\le
2^{(1-\epsilon)\binom{n}{k}} + 2^{\frac12{\binom{n-1}{k-1}}}|\cI(n-1)| +2^{n} \abs{\cI^*\left(n, \zeta,\frac{1}{10k}\right)}.
\]
\end{proposition}

\begin{remark}[Notation for setting sufficiently small constants]
Throughout the next two sections, our arguments often involve picking a sequence of constants. We use $\eps \lll \delta$ to mean that $\eps > 0$ is a sufficiently small constant depending on $\delta > 0$. Here we omit from the notation the dependencies on $k$ and $\theta > \log_2 3$, which we view as fixed throughout.    
\end{remark}

To prove \cref{prop:reducetoboundik}, we refine the proof of \cref{thm:weak-count} from \cref{s:firstcount} by separating containers that are nearly unate from those that far from unate.

\begin{definition}
For every $\rho > 0$, we say that a $k$-SAT formula $G$ is $\rho$-\emph{nearly unate} if it can be made unate by removing up to $\rho n^k$ clauses. 
Otherwise we say that $G$ is $\rho$-\emph{far from unate}.
\end{definition}

\subsection{Nearly complete far from unate simple formulae have many forbidden subformulae}

We call a simple $k$-SAT formula with $k+1$ variables and $k+1$ clauses a \emph{simplex}. 
An example of a simplex with $k=3$ is $\set{wxy,w\ol xz,\ol w y\ol z,xy\ol z}$.

\begin{lemma}\label{lem:blowupsofk+1cliquesnonminimal}
The $2$-blowup of every non-unate simplex is non-minimal.
\end{lemma}

\begin{proof}
The proof is analogous to that of~\cref{prop:fktypenonmin}.
Let $x_1, \ldots, x_{k+1}$ be the variables of the simplex.
Denote the vertices of the $2$-blowup by $x_1, x_1', \ldots, x_{k+1}, x_{k+1}'$.
Without loss of generality (by negating or relabeling variables), the non-unate simplex contains the following three clauses:
\begin{align*}
    C_1 &= x_1 x_2 x_3 \cdots x_{k-1} x_{k}, \\
    C_2 &= \phantom{x_1} y_2 y_3 \cdots y_{k-1} y_{k} x_{k+1}, \\
    C_3 &= z_1 \phantom{x_2} z_3 \cdots z_{k-1} z_{k} \ol{x_{k+1}},
\end{align*}
where $y_i, z_i \in \{x_i, \ol{x_i}\}$ for each $i$.
Let us attempt to find a witness assignment $\mathsf w$ of the variables that satisfies only the clause $x_1\cdots x_k$ in the $2$-blowup.
For each $i = 1, \dots, k$, we have $\mathsf w(x_i) = 1$. Furthermore, since $\mathsf w$ does not satisfy the clause $x_1\cdots x_i' \cdots x_k$  from the $2$-blowup of $C_1$, 
we have
$\mathsf w(x_i') = 0$.
If $\mathsf w(x_{k+1}) = 1$, then there is a clause in the $2$-blowup of $C_2$ that is satisfied under $\mathsf{w}$, 
whereas if $\mathsf w(x_{k+1}) = 0$, then there is clause in the $2$-blowup of $C_3$ that is satisfied. 
So $\mathsf w$ cannot be a witness. Thus the $2$-blowup is non-minimal.
\end{proof}

\begin{lemma}\label{lem:far-from-unate-non-unate-simplex}
For every $\rho > 0$, there exists $\epsilon > 0$ such that for all sufficiently large $n$, every simple $\rho$-far from unate $n$-variable $k$-SAT formula with at least $(1 - \eps)\binom{n}{k}$ clauses contains $\Omega_\rho(n^{k+1})$ non-unate simplices.
\end{lemma}

\begin{proof}
Let $G$ be a simple $\rho$-far from unate $n$-variable $k$-SAT formula with at least $(1 - \eps)\binom{n}{k}$ clauses.
For each variable $x$, let $G_x$ denote the subformula consisting of clauses containing $x$ (as a positive literal), and let $G_{\ol{x}}$ the subformula consisting of clauses containing $\ol x$ (as a negative literal). 

Since $G$ is $\rho$-far from unate, there are at least $\rho n/2$ variables $x$ with 
\begin{equation} \label{eq:far-from-unate-non-unate-simplex-1}
    \min\set{\abs{G_x}, \abs{G_{\ol x}}} \ge \rho n^{k-1}/2.
\end{equation}
Since $\abs{G} \ge (1-\eps) \binom{n}{k}$, for all but at most $k\sqrt{\eps} n$ variables $x$,
\begin{equation} \label{eq:far-from-unate-non-unate-simplex-2}
   \abs{G_x \cup G_{\ol x}} \ge (1-\sqrt{\eps}) \binom{n-1}{k-1}.
\end{equation}
By choosing $\eps$ small enough so that $k\sqrt{\eps} \le \rho/4$, there is some set $S$ of variables with $\abs{S} \ge \rho n/4$ 
so that \cref{eq:far-from-unate-non-unate-simplex-1} and \cref{eq:far-from-unate-non-unate-simplex-2} both hold for every $x \in S$.

For every $x \in S$, by the Kruskal--Katona theorem (\cref{thm:kruskalkatona}),
\[
\# \set{\text{simplices in $G_x$}} \le \frac{((k-1)! \abs{G_x})^{k/(k-1)}}{k!},
\]
and
\[
\# \set{\text{simplices in $G_{\ol x}$}} \le \frac{((k-1)! \abs{G_{\ol x}})^{k/(k-1)}}{k!}.
\]
Adding the two inequalities, and noting that $\abs{G_x} + \abs{G_{\ol x}} \le \binom{n-1}{k-1}$, $\abs{G_x} = \Omega(\rho n^{k-1}), \abs{G_{\ol x}} = \Omega(\rho n^{k-1})$, and the convexity of $t \mapsto t^{k/(k-1)}$, we have
\begin{equation} \label{eq:far-from-unate-non-unate-simplex-3}
\# \set{\text{simplices in $G_x$}}
+ \# \set{\text{simplices in $G_{\ol x}$}} \le (1 - \Omega(\rho))\frac{n^k}{k!}.
\end{equation}

For each $x \in S$, due to \cref{eq:far-from-unate-non-unate-simplex-2}, there are $(1 - O(\sqrt{\eps}))\binom{n-1}{k}$ simplices in $G_x \cup G_{\ol x}$.
Combining with \cref{eq:far-from-unate-non-unate-simplex-3}, we see that there are $\Omega(\rho n^k)$ simplices that simultaneously contain both $x$ and $\ol x$ as literals. Such simplices are non-unate. Since $\abs{S} = \Omega(\rho n)$, there are $\Omega(\rho^2 n^{k+1})$ non-unate simplices.
\end{proof}

\begin{lemma}\label{lem:manyblowupsofk+1cliques}
For every $\rho > 0$, there exists $\eps > 0$ such that for all sufficiently large $n$, every simple $\rho$-far from unate $n$-variable $k$-SAT formula with at least $(1 - \eps){n\choose k}$ clauses contains at least $\Omega_\rho(n^{v(B)})$ copies of some non-minimal simple subformula $B$ on at most $2(k+1)$ variables. 
\end{lemma}

\begin{proof}
Let $G$ be such a formula.  
By \cref{lem:far-from-unate-non-unate-simplex}, $G$ contains $\Omega_\rho(n^{k+1})$ non-unate simplices.
By \cref{lem:formula-blowup}, $G$ contains $\Omega_\rho(n^{2k+2})$ copies of 2-blowups of non-unate simplices, and each such copy, by \cref{lem:blowupsofk+1cliquesnonminimal}, contains a non-minimal subformulae on $2k+2$ variables. The result then follows.
\end{proof}

\subsection{Small and large containers}
Recall from \cref{ss:supersaturation} that the definition of $\alpha_i(G)$, namely that $\alpha_i(G)\binom{n}{k}$ is the number of $k$-tuples of variables that support exactly $i$ clauses. 
Also the weight (\cref{def:weight}) is defined by $\wt(G) = \sum_{i=1}^{2^k} (\log_2 i) \alpha_i(G)$.
We consider a container ``large'' if its weight is close to $1$.

\begin{lemma}[Large containers are nearly unate] \label{lem:lg-nearly-unate}
Suppose $\theta > \log_2 3$ satisfies $\pi(\vec T_k,\theta)=1$.
For every $\rho > 0$, there exists $\delta > 0$ such that if
$G$ is an $n$-variable $k$-SAT formulae with $\wt(G) \ge 1 - \delta$ and at most $\delta n^{v(B)}$ copies of each non-minimal simple formula $B$ on at most $2(k+1)$ variables,
then $G$ is $\rho$-nearly unate.
\end{lemma}

\begin{proof}
Select constants $0 < \delta \lll \eps \lll \rho$.
Since $G$ has at most $\delta n^{v(B)}$ copies of each non-minimal simple formula $B$ on at most $2(k+1)$ variables,
we have $\alpha_1(G) + \theta \alpha'_2(G) \le 1 + \eps$ by \cref{lem:supersat-1},
and $\alpha''_2(G) \le \eps$ and $\alpha_i(G)\le \eps$ for each $i \ge 3$ by \cref{lem:supersat-2}.
Thus
\begin{align*}
1 - \delta < \wt(G) 
&= \alpha_1(G) + (\log_2 3) \alpha'_2(G) + O(\eps) 
\\
&= \alpha_1(G) + \theta \alpha'_2(G) - (\theta - \log_2 3) \alpha'_2(G) + O(\eps)
\\
&\le 1 - (\theta - \log_2 3) \alpha'_2(G) + O(\eps).
\end{align*}
Thus $\alpha'_2(G) = O(\eps)$ (recall that we treat $\theta > \log_2 3$ as a fixed constant whose dependence is omitted),
and $\alpha_1(G) = 1 - O(\eps)$. Hence, by \cref{lem:manyblowupsofk+1cliques}, $G$ is $\rho$-nearly unate. 
\end{proof}

\begin{lemma}\label{lem:dense-witness}For every $\zeta>0$, there exists $\rho>0$ such that the following holds. Suppose  $(G,H)$ is a pair of formulae such that $G$ is $\rho$-nearly unate,   $H$  is minimal and $H\subseteq G$. Then, at least one of the following holds for $(G,H)$:
    \begin{enumerate}
        \item $H$ is not $\frac{1}{10k}$-dense.
        \item There exists a variable negation of $G$ into $G'$ such that $G'$ has $\leq \rho n^k$ non-monotone clauses, but the resulting $H'\subseteq G'$ is not $\zeta/2$-sparsely minimal.
        \item We can negate a subset of variables to obtain $H^\uparrow\in\mathcal I^*(n,\zeta,\frac{1}{10k})$ from $H$.
    \end{enumerate}
\end{lemma}

\begin{proof}
    Suppose $(G,H)$ does not satisfy conditions (1)(2) in lemma. Since $G$ is $\rho$-nearly unate, we can negate a subset of variables to obtain some formula $G'$ that has $\leq \rho n^k$ non-monotone clauses. Consequently, the formula $H'\subseteq G'$ also has $\leq \rho n^k$ non-monotone clauses.
    
    Let $X_{H'}$ be the set of variables with $m(x)<m(\overline{x})$ in $H'$. Observe that 
    $|X_{H'}|\leq 100k^{k+1}\rho n$,
    as otherwise, since $H'$ is $\frac{1}{10k}$-dense, there would be more than
    \[
    100k^{k+1}\rho n\cdot \frac{1}{20k}\binom{n-1}{k-1}\cdot \frac{1}{k}>\rho n^k
    \]
    non-monotone clauses in $H'$. Therefore, by negating those variables in $X_{H'}$, we can obtain a positive-dominant formula $H^\uparrow$ from $H'$ that has at most (suppose $\rho\lll\zeta$)
    \[
    100k^{k+1}\rho n\binom{n-1}{k-1}+\rho n^k\leq \zeta n^k
    \]
    non-monotone clauses. It is clear that $H^\uparrow$ is $\frac{1}{10k}$-dense. Since every clause in $H'$ has a witness assignment with fewer than $\zeta n/2$ variables assigned to 1,  every clause in $H^\uparrow$ has a witness assignment with fewer than $\zeta n/2+100k^{k+1}\rho n \leq \zeta n$ variables assigned to 1. Thus, $H^\uparrow$ is $\zeta$-sparsely minimal. Combining the above gives that $H^\uparrow\in\mathcal I^*(n,\zeta,\frac{1}{10k})$.
\end{proof}

\begin{lemma}\label{lem:dense-not-witness}
For every $\zeta>0$, there exists $\rho>0$ such that the following holds for all sufficiently large $n$:
if $\cG$ is a collection of $\rho$-nearly unate $n$-variable $k$-SAT formulae with $\abs{\cG} \le 2^{\rho n^k}$, then 
\begin{align*}
    &|\{(G,H):\text{$G\in\cG$, $H\subseteq G$, $H$ minimal and $(G,H)$ meets condition (2) in \cref{lem:dense-witness}}\}|\leq 2^{(1-\rho)\binom{n}{k}}.
\end{align*}
\end{lemma}

\begin{proof}
Let $\cG$ be a collection of $\rho$-nearly unate $n$-variable $k$-SAT formulae, with $\abs{\cG} \le 2^{\rho n^k}$.
Per~\cref{lem:dense-witness}, we will count pairs $(G,H)$ such that
\begin{enumerate}
    \item  $G\in\cG$, $H\subseteq G$, $H$ is minimal;
    \item there exists a variable negation of $G$ into $G'$ such that $G'$ has $\leq \rho n^k$ non-monotone clauses, but the resulting $H'\subseteq G'$ is not $\zeta/2$-sparsely minimal.
\end{enumerate}
Every such $(G,H)$ is uniquely determined by the following data:
\begin{enumerate}[(a)]
    \item $G\in\cG$ ($\leq 2^{\rho n^k}$ possibilities);
    \item The variable negation of $G$ into $G'$ ($\leq 2^n$ possibilities);
    \item The non-monotone clauses of $H'\subseteq G'$ ($\le 2^{\rho n^k}$ possibilities since $G'$ is $\rho$-nearly monotone);
    \item The $\ceil{\zeta n/2}$ lexicographically smallest variables $x_i$ where $\mathsf{w}_{H'}(x_i) = 1$ for some witness $\mathsf{w}_{H'}$ to clause $C_{H'}$ ($\le \binom{n}{\ceil{\zeta n/2}}$ possibilities);
    \item All monotone clause(s) supported on the above $\ceil{\zeta n/2}$ variables (which must be either uniquely $C_{H'}$ or none due to $\mathsf w$ being a witness for $C_{H'}$; $\le \binom{\ceil{\zeta n/2}}{k} + 1$ possibilities);
    \item The remaining monotone clauses of $H'$ not covered by (e) ($\le 2^{\binom{n}{k} - \binom{\ceil{\zeta n/2}}{k}}$ possibilities).
\end{enumerate}
Combining the bounds, we see that with $\rho\lll\zeta$, the number of such pairs is at most
\begin{align*}
&\le 
\underbrace{2^{\rho n^k}}_{(a)}\cdot \underbrace{2^{n}}_{(b)}\cdot \underbrace{2^{\rho n^k}}_{(c)} \cdot 
\underbrace{\binom{n}{\ceil{\zeta n/2}}}_{(d)} \cdot 
\underbrace{\left({\ceil{\zeta n/2} \choose k} + 1\right)}_{(e)} \cdot 
\underbrace{2^{\binom{n}{k} - { \ceil{\zeta n/2}  \choose k}}}_{(f)} \\
&\le 2^{(1 - \rho)\binom{n}{k}}. \qedhere 
\end{align*}
\end{proof}

We will use $H(\cdot)$ to denote the \emph{binary entropy function} $H(p)=-p\log_2p-(1-p)\log_2(1-p)$, which satisfies the standard inequality
\begin{equation}\label{eq:bin-ent}
{n\choose \le \lfloor pn \rfloor }\leq 2^{H(p)n} \quad \text{ for all } 0<p\le 1/2 \text{ and } n \in \NN.
\end{equation}

\begin{proof}[Proof of \cref{prop:reducetoboundik}]

We first count the number of formulae in $\cI(n)$ that are not $\frac{1}{10k}$-dense. Every formula  $G\in\cI(n)$ that is not $\frac{1}{10k}$-dense is uniquely determined by the following information:
\begin{enumerate}
    \item The lexicographically smallest variable $x\in X$ that is used by less than $\frac{1}{10 k}{n-1\choose k-1}$ clauses in $G$;
    \item The clauses in $G$ that use the variable $x$;
    \item The clauses in $G$ that do not use the variable $x$.
\end{enumerate}
We get the following bound by upper bounding the number of choices for each item:
\begin{align*}
    |\{G\in \cI(n): G\text{ is not $\frac{1}{10k}$-dense}\}| &\leq \underbrace{n}_{(1)}\underbrace{\binom{2^k\binom{n-1}{k-1}}{\leq \frac{1}{10 k}{n-1\choose k-1}}}_{(2)}\underbrace{|\cI(n-1)|}_{(3)}
      < 2^{\frac12{n-1\choose k-1}}|\cI(n-1)|,
\end{align*}
where the last inequality uses~\cref{eq:bin-ent}, since $H\left(\frac{1}{ 10k\cdot 2^k}\right)2^k < \frac13$ for $k \ge 3$ and $n$ is sufficiently large.

It remains to upper bound $|\{G\in \cI(n): G\text{ is $\frac{1}{10k}$-dense}\}|$. To do so, we will choose a sequence of constants
\[
0< \delta \lll \rho \lll \zeta.
\]
As in \cref{s:firstcount}, let $\cB$ be the set of all simple non-minimal $k\SAT$ formulae on at most $2(k+1)$ vertices. Applying the container theorem, \cref{thm:gencontub}, we obtain a collection $\cG$ of formulae satisfying:
    \begin{itemize}
        \item Every $\cB$-free $n$-variable $k$-SAT formula is a subformula of some $G \in \cG$;
        \item For every $G \in \cG$ and $B \in \cB$, $G$ has at most $\delta n^{v(B)}$ copies of $B$;
        \item $\abs{\cG} = 2^{o(n^k)}$.
    \end{itemize}

We partition $\cG =\cG^{\sm} \cup \cG^{\lg}$ where
\[
    \cG^\sm=\left\{G\in \cG \mid \wt (G)\leq 1- \delta \right\} 
    \quad\text{and}\quad
    \cG^\lg=\left\{G\in \cG \mid \wt (G)> 1- \delta \right\}.
\]
We know from \cref{prop:nearly-simple} that every minimal formulae can be made simple by deleting $o(n^k)$ clauses. 
The number of simple subformulae of each $G \in \cG^\sm$ is $2^{\wt(G)\binom{n}{k}} \le 2^{(1-\delta)\binom{n}{k}}$, and so the number of minimal subformulae of $G$ is at most $2^{(1-\delta + o(1))\binom{n}{k}}$. 
Taking a union bound over all $G \in \cG^\sm$, as $\abs{\cG^\sm} = 2^{o(n^k)}$, 
\[
\abs{\bigcup_{G\in \cG^\sm}\set{\text{minimal } H\subseteq G}} \le 2^{(1-\delta + o(1))\binom{n}{k}}.
\]
By \cref{lem:lg-nearly-unate}, every $G \in \cG^\lg$ is $\rho$-nearly unate. 
By \cref{lem:dense-witness} and \cref{lem:dense-not-witness}, we have
\[
\abs{\bigcup_{G\in \cG^\lg}\set{\text{minimal, $\frac{1}{10k}$-dense }H\subseteq G}} 
\le
2^{(1-\rho)\binom{n}{k}} + 2^n\abs{\cI^*(n, \zeta, \frac{1}{10k})}.
\]
Since every minimal formula is a subformula of some element of $\cG =\cG^{\sm} \cup \cG^{\lg}$, combining these two bounds, we get that

\[
\left|\left\{G\in \cI(n): G\text{ is  $\frac{1}{10k}$-dense}\right\}\right|\leq 2^{(1-\eps)\binom{n}{k}}+2^n\left|\cI^*\left(n,\zeta,\frac{1}{10k}\right)\right|. \qedhere
\]
\end{proof}

\section{Stability II: nearly monotone and sparsely minimal  formulae}\label{s:istar}

Recall that in the previous section, we proved~\cref{prop:reducetoboundik}, which reduced the problem of bounding the number of minimal $k$-SAT formula to counting formulae in

 $$\cI_1^*(n, \zeta) := \cI^*\left(n, \zeta, \frac{1}{10k}\right)$$
(see \cref{d:cistar}). 
In this section, we complete the proof of~\cref{thm:strong-count} by establishing a recursive inequality for $|\cI_1^*(n, \zeta)|$. 
The inequality is analogous to equation (20) in Ilinca--Kahn \cite{IK12}.

\begin{theorem}\label{t:boundcistar}
For $k \ge 3$, there exist $\zeta>0$ and $c>0$ such that for all $n$,
\begin{multline*}
\left|\cI_1^*\left(n, \zeta\right)\right| \le 2^{n\choose k}+\exp_2\sqb{(1-c){n\choose k}}+\sum_{i=1}^{k-2}\exp_2\sqb{i(1-c){n\choose k-1}}|\cI(n-i)|
\\+\exp_2\sqb{\frac{k}{2}(1-c){n\choose k-1}} |\cI(n-k)|+\exp_2\sqb{{n\choose k}-cn}.
\end{multline*}
When $k = 2$, we have 
\[\left|\cI_1^*\left(n, \zeta\right)\right| \le 2^{n\choose k}+2\exp_2\sqb{(1-c){n\choose k}}+\exp_2\sqb{\frac{k}{2}(1-c){n\choose k-1}} |\cI(n-k)|.\]
\end{theorem}

We first prove \cref{thm:strong-count} assuming \cref{t:boundcistar}.

\begin{proof}[Proof of \cref{thm:strong-count}] Choose  $\theta>\log_23$ with $\pi(\vec T_k,\theta)= 1$ and  $\zeta>0$ as in \cref{t:boundcistar}. We then choose $\epsilon>0$ as in \cref{prop:reducetoboundik}, so that for $n$ sufficiently large, we have
$$|\cI(n)|\leq  2^{\left(1-\epsilon\right)\binom{n}{k}} + 2^{\frac12{n-1 \choose k-1}} |\mcI(n-1)| + 2^{n} \abs{\cI_1^*\left(n, \zeta\right)}.$$
We show by induction that there exist $c_1, c_0>0$ such that for all $n$,
\begin{equation}\label{eq:strongub}
    \abs{\cI(n)}\leq (1+ 2^{c_1-c_0n}) 2^{{n\choose k}+n}.
\end{equation}

We first choose $0<c<1$ sufficiently small so that \cref{t:boundcistar} holds. We then choose $n_0\in\NN$ such that  for all $n\geq n_0$, the following conditions hold:
\begin{enumerate}[(i)]
\item $\binom{n}{k}-\binom{n-i}{k}    >i\left(1-\frac{c}{2}\right)\binom{n}{k-1}$ for all $i\in[k]$; 
\item if $k\geq 3$, we require that $\frac{c}{2}\binom{n}{k-1}\geq 2n$;
\item $2^{-c\binom{n}{k}+1}+k2^{-n+1}+2^{-cn+1}\leq 2^{-\frac{c}{k}n}$.
\end{enumerate}

Choose $c_0>0$ such that for all $n>n_0$, we have
 \begin{equation}\label{eq:c0}
   2^{-n} + 2^{-\frac12 {n - 1 \choose k-1}}  + 2^{ -\frac{c}{k} n}\leq 2^{-c_0n}.
\end{equation}
Then, take $c_1>0$ such that~\cref{eq:strongub} holds for all $n \leq n_0$.
It remains to show that ~\cref{eq:strongub} holds for all $n > n_0$. Fix one such $n$ and assume by induction that~\cref{eq:strongub} holds for all smaller values of $n.$ By the inductive hypothesis and condition (i) on $c$, for every $i\in[k]$, we have
\begin{align*}
    \exp_2\sqb{i(1-c){n\choose k-1}}|\cI(n-i)|&\leq (1 + 2^{c_1 - c_0 (n-i)}) 2^{{n-i \choose k} + n  - i}\cdot 2^{i(1-c){n\choose k-1}}  \\
    &\leq (1 + 2^{c_1 - c_0 (n-i)}) 2^{{n \choose k} + n - \frac{c}{2}{n \choose k-1} - i}.
\end{align*}
When $k \ge 3$, substituting the above into~\cref{t:boundcistar} gives that 
\begin{align*}
\left|\cI_1^*\left(n, \zeta\right)\right| &\le 2^{n\choose k}+\exp_2\sqb{(1-c){n\choose k}}+\sum_{i=1}^{k}\exp_2\sqb{i(1-c){n\choose k-1}}|\cI(n-i)|
 +\exp_2\sqb{{n\choose k}-cn} \\
&\leq 2^{n\choose k}+\exp_2\sqb{(1-c){n\choose k}}+\sum_{i=1}^{k}(1 + 2^{c_1 - c_0 (n-i)}) 2^{{n \choose k} + n - \frac{c}{2}{n \choose k-1} - i} 
+\exp_2\sqb{{n\choose k}-cn}.
\end{align*}
Let 
\[
J_{n,k}:= \sum_{i=1}^{k}(1 + 2^{c_1 - c_0 (n-i)}) 2^{{n \choose k} + n - \frac{c}{2}{n \choose k-1} - i}.
\]
We know from  condition (ii) that $\frac{c}{2}\binom{n}{k-1}\geq 2n$.  Thus, we have 
\[
J_{n,k}\leq k(1+2^{c_1})2^{\binom{n}{k}-n},
\]
which gives
\begin{align*}
\left|\cI_1^*\left(n, \zeta\right)\right|&\leq 2^{{n \choose k}}+2^{(1-c){n \choose k}}+k(1+2^{c_1})2^{\binom{n}{k}-n}+2^{\binom{n}{k}-cn}\\
& \leq (1+2^{c_1}(2^{-c\binom{n}{k}}+k2^{-n+1}+2^{-cn}))2^{{n \choose k}} \leq (1 + 2^{c_1 - \frac{c}{k} n}) 2^{{n \choose k}},
\end{align*}
where the final inequality follows from condition (iii).

When $k=2$, we apply the simpler recurrence of~\cref{t:boundcistar} and obtain:
\begin{align*}
    \mcI^*_1(n, \zeta) &\leq 2^{{n \choose 2}} +2 \exp_2\sqb{(1-c){n\choose 2}} + 2^{(1-c)n}(1+2^{c_1})2^{\binom{n-1}{2}+n-1} \\
    &
    \le 2^{{n \choose 2}} +2 \exp_2\sqb{(1-c){n\choose 2}} + (1+2^{c_1})2^{\binom{n}{2}-cn+1}\\
    &\leq (1+2^{c_1}(2^{-c\binom{n}{2}+1}+2^{-cn+1}))2^{\binom{n}{2}}
    \leq (1+2^{c_1-\frac{c}{2}n})2^{\binom{n}{2}},
\end{align*}
where the final inequality follows from condition (iii).

\cref{prop:reducetoboundik} then implies
\begin{align*}
|\cI(n)|&\leq 2^{\left(1-\epsilon\right)\binom{n}{k}}+  2^{\frac12{n-1 \choose k-1}} |\mcI(n-1)| + 2^n|\cI_1^*\left(n, \zeta\right)|\\
&\leq (2^{-n} + 2^{-\frac12 {n - 1 \choose k-1}+c_1} + 1 + 2^{c_1 -\frac{c}{k} n} ) 2^{{n\choose k}+n} \\
&\le (1 + 2^{c_1-c_0 n}) 2^{n + {n \choose k}},
\end{align*}
where the final inequality is a consequence of our choice of $c_0$ satisfying~\cref{eq:c0}.
This proves the inductive hypothesis.
\end{proof}

Our proof of \cref{t:boundcistar} extends the arguments in \cite[Section 8]{IK12}. 
The following steps~\cref{s:startik} to~\cref{s:endik} roughly correspond to Steps 1--5 in~\cite[Section 8]{IK12}.

\subsection{Basic properties of $\cI_1^*$}\label{s:startik}

It remains for us to count formulae in $\cI_1^*(n, \zeta)$.
We recall that such formulae enjoy the following properties:
\begin{itemize}
  \item $m(x) \ge m(\ol x)$ for every variable $x \in X$; 
    \item every such formula has at most $\zeta n^k$ non-monotone clauses;
    \item every clause has a witness $\mathsf w$ supported on fewer than $< \zeta n$ variables;
    \item every variable $x \in X$ is used by at least $\frac{1}{10k} {n - 1 \choose k-1}$ clauses.
\end{itemize}

We will show that in addition, most clauses are monotone in every formula in $\cI^*_1(n, \zeta)$. The following  observations will be useful in our subsequent steps.

\begin{lemma}\label{lem:ci_1starproperties}
For every $G \in \cI_1^*(n, \zeta)$, the following are true.
\begin{enumerate}[label=(\alph*)]
    \item For every $k-1$ literals $z_1,\dots,z_{k-1}$,
    there are at most $\zeta n$ variables $w$ such that $z_1\cdots z_{k-1}\ol w\in G$. 
    \item For every variable $v$ and $j\in\{1,\dots,k-1\}$, there are at most $\zeta n^{k-1}$ clauses in $G$ containing the positive literal $v$ and exactly $j$ negative literals.
    \item For every variable $v$ and $j\in\{1,\dots,k-1\}$, there are at most $\zeta n^{k-1}$ clauses in $G$ containing the negative literal $\ol v$ and exactly $j$ negative literals other than $\ol v$.
\end{enumerate}
\end{lemma}
\begin{proof}Consider any formula $G\in \cI_1^*(n, \zeta)$. We verify properties (a)--(c):
\begin{enumerate}[label=(\alph*)] 
    \item Fix $k-1$ literals $z_1,\dots,z_{k-1}$. Suppose there are more than $\zeta n$ variables $w$ such that $z_1\cdots z_{k-1}\ol w\in G$. Let $w_0$ be one of these variables, and $\mathsf w\in\{0,1\}^n$ be the witness that satisfies $z_1\cdots z_{k-1}\ol {w_0}$ but no other clauses in $G$. Then for any other $w\neq w_0$ such that  $z_1\cdots z_{k-1}\overline{w} \in G$, $\mathsf w$ must assign $w$ to $1$. But then $\mathsf w$ assigns at least $\zeta n$ variables to $1$, contradicting the definition of $\cI^*(n, \zeta, \zeta')$ (\cref{d:cistar}).
    \item Fix a variable $v$. By (a), for any $j-1$ negative literals $\ol {v_1}, \dots, \ol {v_{j-1}}$ and $k-1-j$ positive literals $v_j,\ldots,v_{k-2}$,  there are at most $\zeta n$ variables $w$ such that $G$ contain $v\ol v_1 \cdots \ol v_{j-1} v_j\cdots v_{k-2}\ol w$ as a clause. Summing over all possible choices $\ol {v_1}, \dots, \ol {v_{j-1}}$ and $v_j,\ldots,v_{k-2}$ gives the result.
    \item Fix a variable $v$. By (a), for any $j-1$ negative literals $\ol {v_1}, \dots, \ol {v_{j-1}}$ and $k-1-j$ positive literals $v_j,\ldots,v_{k-2}$, there are at most $\zeta n$ variables $w$ such that $G$ contain $\ol v\ol v_1 \cdots \ol v_{j-1} v_j\cdots v_{k-2}\ol w$ as a clause. Summing over all possible choices $\ol {v_1}, \dots, \ol {v_{j-1}}$ and $v_j,\ldots,v_{k-2}$ gives the result. \qedhere
\end{enumerate}
\end{proof}

\subsection{Most minimal formulae have few clauses with negated literals} 
We will eventually show that most formulae in $\cI_1^*(n, \zeta)$ have  few non-monotone clauses. In this subsection, we give an upper bound on the number of formulae in $\cI_1^*(n, \zeta)$ that have a large numbers of clauses containing $i$ negative literals, for any $i \in [k-2].$

\begin{definition}\label{def:I2}
For every $\beta_1>0$, define 
\begin{multline*}
    \cI^*_{2,1}(n, \zeta,\beta_1)=\left\{G\in \cI_{1}^*(\zeta):\text{for all $u\in X$, $G$ contains at most $\beta_1n^{k-1}$ clauses}\right. \\\left.\text{of the form } \ol uv_1\cdots v_{k-1}\right\}.
\end{multline*}
For $k \ge 4$, for every $i\in\{2,\dots,k-2\}$ and $\beta_i>0$, define
\begin{multline*}
    \cI^*_{2,i}(n, \zeta,\beta_1,\beta_i)=\left\{G\in \cI_{2,1}^*(n, \zeta,\beta_1):\text{for all $\{u_1,\dots,u_i\}\in{X\choose i}$, }\right. \\\left.\text{$G$ contains at most $\beta_in^{k-i}$ clauses of the form } \ol u_1\cdots\ol u_iv_1\cdots v_{k-i}\right\}.
\end{multline*}
For $k \ge 3$ and every $\vec \beta > 0$, define
\begin{align*}
    \cI_2^*(n, \zeta,\vec\beta)&=\bigcap_{i=2}^{k-2}\cI^*_{2,i}(n, \zeta,\beta_1,\beta_i)\\
    &=\bigl\{G\in\cI^*(\zeta):\text{for all distinct variables } u_1,\dots,u_i \in X \text{ where } i\in [k-2],\\
    &\hspace{3cm} \text{$G$ contains at most $\beta_{i}n^{k-i}$ clauses of the form $\ol{u_1}\dots\ol{u_i}v_1\dots v_{k-i}$}\bigr\}.
\end{align*}
For $k = 2$, we define $\mcI_2^*(n, \zeta, \vec\beta) = \cI^*_{2,1}(n, \zeta,\beta_1)$.
\end{definition}
We will first show that most formulae in $\cI_1^*(n, \zeta)$ lie in $\cI_{2,1}^*(n, \zeta,\beta_1)$.

\subsubsection{Formulae with lots of clauses $\overline{u}v_1 \cdots v_{k-1}$}

We begin by studying the $k = 2$ case.
\begin{lemma}\label{l:step2.1-k=2}
Suppose $k=2$. For every $\beta_1<1/2$, if $\zeta$ is sufficiently small relative to $\beta_1$, there exists $c > 0$ such that
\begin{align*}
    |\cI^*_{1}(n, \zeta)\setminus \cI_{2,1}^*(n, \zeta,\beta_1)|\leq \exp_2\sqb{(1-c){n\choose k}}.
\end{align*}
\end{lemma}
\begin{proof}
    For every formula $G\in \cI_1^*(n,\zeta)\setminus\cI_{2,1}^*(n,\zeta)$, there exists some variable $u\in X$ such that $G$ contains more than $\beta_1n$ clauses of the form $\ol uv$. From the fact that (1) $u$ is used by at least $\frac{n}{20}$ clauses in $G$,  (2) $m(u)\ge m(\ol u)$ and (3) $G$ contains at most $\zeta n$ clauses of the form $u\ol v$, we know that $G$ contains at least $(1/40-\zeta)n$ clauses of the form $uv$. Define
    \begin{align*}
        N_u&=\{v\in X_{-u}:uv\in G\},\\
        N_{\ol u}&=\{v\in X_{-u}:u\ol v\in G\}.
    \end{align*}
Observe that, to avoid the non-minimal formula $\{\ol {u}v_2, uv_1,v_1v_2\}$, all monotone clauses $v_1v_2$ with $v_1\in N_u$, $v_2\in N_{\ol u}$ cannot lie in $G$. Thus, given $N_u$ and $N_{\ol u}$, we know at least $|N_u||N_{\ol u}|-n$ monotone clauses that cannot lie in $G$. 

Each formula $G\in \cI_1^*(n,\zeta)\setminus\cI_{2,1}^*(n,\zeta,\beta_1)$ is then uniquely determined by specifying the following in order:
\begin{enumerate}
    \item $u\in X$;
    \item $N_u, N_{\ol u}\subseteq X$;
    \item $\{v\in X:\ol uv\in G\}$ and $\{v\in X:\ol u\ol v\in G\}$;
    \item the non-monotone clauses in $G$ that do not use the variable $u$.
    \item the monotone clauses in $G$ that do not use the variable $u$.
\end{enumerate}
We then have the following upper bound on $|\cI^*_{1}(n, \zeta)\setminus \cI_{2,1}^*(n, \zeta,\beta_1)|$, in terms of the number of choices for each item.
\[
|\cI^*_{1}(n, \zeta)\setminus \cI_{2,1}^*(n, \zeta,\beta_1)|
\le
\underbrace{n}_{(1)}\underbrace{5^n}_{(2)(3)}\underbrace{\binom{(2^2-1)\binom{n}{2}}{\leq \zeta n^k}}_{(4)}\underbrace{2^{\binom{n}{2}-\beta_1(1/40-\zeta)n^2-n}}_{(5)}\leq \exp_2\sqb{(1-c)\binom{n}{2}}. 
\qedhere 
\]
\end{proof}

We now turn to the $k \ge 3$ case.
\begin{lemma}\label{l:step2.1}
For $k \ge 3$ and any positive $\beta_1 < 1/k$, if $\zeta$ is sufficiently small relative to $\beta_1$, there exists $c > 0$
such that
$$|\cI^*_{1}(n, \zeta)\setminus \cI_{2,1}^*(n, \zeta,\beta_1)|\leq \exp_2\sqb{(1-c){n\choose k}}+\exp_2\sqb{(1-c){n\choose k-1}}|\cI (n-1)|.$$
\end{lemma}

\begin{proof}
Choose constants $0<c\lll\zeta\lll\theta_1\lll\beta_1 \lll 1/k$.

For every formula $G\in \cI^*_{1}(n,\zeta)\setminus \cI_{2,1}^*(n, \zeta,\beta_1)$, by definition, there is some variable $u\in X$ such that $G$ contains more than $\beta_1n^{k-1}$ clauses of the form $\ol uv_1\cdots v_{k-1}$. Fix this variable $u$ and let $X_{-u}=X\setminus \{u\}$. Define the \emph{monotone links} of literals $u$ and $\ol u$ to be the $(k-1)$-graphs with edge sets given by
\begin{align*}
    N_u&=\left\{\{v_1,\dots,v_{k-1}\}\subseteq X_{-u}:uv_1\cdots v_{k-1}\in G\right\},\\
    N_{\ol u}&=\left\{\{v_1,\dots,v_{k-1}\}\subseteq X_{-u}:\ol uv_1\cdots v_{k-1}\in G\right\}.
\end{align*}
Now consider the collections of $(k-2)$-subsets of $X_{-u}$ that have small co-degree (less than $\theta_1 n$) in the hypergraphs $N_u, N_{-u}$. Precisely,
\begin{align*}
    S_u&=\left\{\{v_1,\dots,v_{k-2}\}\subseteq X_{-u} \text{ such that }|\{w:\{w,v_1,\dots, v_{k-2}\}\in N_{u}\}|<\theta_1n\right\},\\
    S_{\ol u}&=\left\{\{v_1,\dots,v_{k-2}\}\subseteq X_{-u}\text{ such that }|\{w:\{w,v_1,\dots, v_{k-2}\}\in N_{\ol u}\}|<\theta_1n\right\}.
\end{align*}
Finally, let $T_{-u}$ be the $(k-2)$-subsets of $X_{-u}$ that fall in neither $S_u$ nor $S_{\ol u}$, i.e. that have large co-degree into both $N_u$ and $N_{\ol u}$, so that 
$
T_{-u}={X_{-u}\choose k-2}\setminus(S_u\cup S_{\ol u}).
$
Since $|S_u|, |S_{\ol u}| < {n \choose k-2}$, we see that
\begin{align*}
|N_{\ol u}|&\leq \theta_1 n|S_{\ol u}|+n\left(|S_{ u}|+|T_{-u}|\right)
    \leq n(|S_{u}|+|T_{-u}|)+\frac{\theta_1 n^{k-1}}{(k-2)!}.\\
\end{align*}
We assumed that $G \not \in \cI_{2,1}^*(n, \zeta, \beta_1)$, so $|N_{\ol u}|\geq \beta_1n^{k-1}$. 

Combining this lower bound with the above upper bound on $|N_{\ol u}|$ means that at least one of the following must hold:
\begin{enumerate}
    \item $|T_{-u}|>\theta_1 n^{k-2}$; 
    \item $|T_{-u}|\leq \theta_1 n^{k-2}$; in this case, from the above inequalities, we can get a lower bound for $|S_u|$: $$|S_u|\geq \beta_1n^{k-2}-\frac{\theta_1 n^{k-2}}{(k-2)!} -|T_{-u}| \geq \beta_1n^{k-2}-\frac{\theta_1n^{k-2}}{(k-2)!}-\theta_1n^{k-2}>\theta_1 n^{k-2},$$ as $\theta_1 \lll\beta_1$.
\end{enumerate}

We upper bound the number of formulae in $\cI^*_{1}(n,\zeta)\setminus \cI_{2,1}^*(n, \zeta,\beta_1)$ in each of these two cases.
\begin{claim}\label{claim:I21-case1}
There are  at most $\exp_2\left[(1-c){n\choose k}\right]$ formulae $G\in \cI^*_{1}(\zeta)\setminus \cI_{2,1}^*(n, \zeta,\beta_1)$ such that  for some variable $u\in X$, we have $|T_{-u}| > \theta_1 n^{k-2}$.
\end{claim}
\begin{proof*}
Each such formula $G$ is uniquely determined by specifying the following information in order:
\begin{enumerate}
    \item $u\in X$;
    \item $T_{-u}\subseteq\binom{X}{k-2}$;
    \item $\{ \{v_1,\dots,v_{k-1}\}\in N_u: \text{there exists $i\in[k-1]$ such that $\{v_1,\dots,v_{k-1}\}\setminus \{v_i\}\in T_{-u}$}\}$;
    \item $\{ \{v_1,\dots,v_{k-1}\}\in N_{\ol u}: \text{there exists $i\in[k-1]$ such that $\{v_1,\dots,v_{k-1}\}\setminus \{v_i\}\in T_{-u}$}\}$;
    \item The non-monotone clauses in $G$;
    \item The monotone clauses in $G$.
\end{enumerate}
We prove the claim by giving an upper bound on the number of choices for each piece of information:
\begin{align*}
    \underbrace{n}_{(1)} \cdot \underbrace{2^{\binom{n}{k-2}}}_{(2)} \cdot \underbrace{2^{|T_{-u}|n}}_{(3)} \cdot \underbrace{2^{|T_{-u}|n}}_{(4)} \cdot \underbrace{\binom{(2^k-1)
\binom{n}{k})}{\leq \zeta n^k}}_{(5)}\underbrace{\exp_2\left[{n\choose k}- \Omega(n^k)\right]}_{(6)}<\exp_2\left[(1-c){n\choose k}\right]
\end{align*}
(recall $c\lll\zeta\lll\theta_1$). In particular, bounds (5) and (6) come from the following observations:
\begin{enumerate}
    \item[(5)] Since $G \in \cI_1^*(n, \zeta)$ has at most $\zeta n^k$ non-monotone clauses.
    \item[(6)] Suppose we have specified (1)--(4). Recall from the definition of $T_{-u}$ that for every $\{v_1,\dots, $ $v_{k-2}\}\in T_{-u}$, there are at least $\theta_1n$ variables $w\in X$ such that $v_1\cdots v_{k-2}wu\in G$, and  at least $\theta_1n$ variables $w'\in X$ such that $v_1\cdots v_{k-2}w'\ol u\in G$. For every such pair $w,w'$, the monotone clause $v_1\cdots v_{k-2}ww'$ is forbidden from $G$, else $G$ contains the non-minimal formula
$$\left\{v_1\cdots v_{k-2}ww',\quad v_1\cdots v_{k-2}wu,\quad v_1\cdots v_{k-2}w'\ol u\right\}$$
(as it is impossible to satisfy only the first clause). Therefore, given (1)--(4) and the hypothesis that $|T_{-u}|>\theta_1n^{k-2}$, there are at least $|T_{-u}|\cdot (\theta_1n)^2/\binom{k}{2} = \Omega(n^k)$ monotone clauses that cannot lie in $G$. \qedhere
\end{enumerate}
\end{proof*}

\begin{claim}\label{claim:I21-case2}
There are  at most $\exp_2\left[(1-c){n\choose k-1}\right]|\cI(n-1)|$ formulae $G\in \cI^*_{1}(\zeta)\setminus \cI_{2,1}^*(n, \zeta,\beta_1)$ such that  for some variable $u\in X$, we have $|T_{-u}|\leq \theta_1 n^{k-2}$ and $|S_u|>\theta_1 n^{k-2}$.

\end{claim}
\begin{proof*}
We follow a similarly inspired approach to the previous claim; each such formula $G$ is uniquely determined by specifying the following information in order:
\begin{enumerate}
    \item $u\in X$;
    \item $S_u,S_{\ol u},T_{-u}\subseteq\binom{X_{-u}}{k-2}$;
    \item $\{ \{v_1,\dots,v_{k-1}\}\in N_u: \text{there exists $i\in[k-1]$ such that $\{v_1,\dots,v_{k-1}\}\setminus \{v_i\}\in S_u\cup T_{-u}$}\}$;
    \item $\{ \{v_1,\dots,v_{k-1}\}\in N_{\ol u}: \text{there exists $i\in[k-1]$ such that $\{v_1,\dots,v_{k-1}\}\setminus \{v_i\}\in S_{\ol u}\cup T_{-u}$}\}$;
    \item $\{ \{v_1,\dots,v_{k-1}\}\in N_u: \text{for all $i\in[k-1]$, $\{v_1,\dots,v_{k-1}\}\setminus \{v_i\}\in  S_{\ol u}\setminus S_u$}\}$;
    \item $\{ \{v_1,\dots,v_{k-1}\}\in N_{\ol u}: \text{for all $i\in[k-1]$, $\{v_1,\dots,v_{k-1}\}\setminus \{v_i\}\in S_u\setminus S_{\ol u}$}\}$;
    \item The monotone clauses in $G$ containing $u$;
    \item The clauses in $G$ containing $\ol u$ and $k-1$ positive literals;
    \item The clauses in $G$ that use the variable $u$ and contain some negative literal other than $\ol u$;
    \item The clauses in $G$ that do not use the variable $u$.
\end{enumerate}

Note that (7) and (8) are uniquely determined by (3)--(6). We prove the claim by giving an upper bound on the number of choices for each item:
\begin{align*}
\underbrace{n}_{(1)}&\underbrace{4^{\binom{n}{k-2}}}_{(2)}\underbrace{\exp_2[ (2\theta_1+H(\theta_1))n^{k-1}]}_{(3)(4)} \underbrace{\exp_2 \sqb{\left(\theta_1^2 + \frac{1}{(k-1)!} \right) n^{k-1}}}_{(5)(6)} \underbrace{\exp_2\left[\frac{1}{k} {n \choose k-1}\right]}_{(9)}\underbrace{|\cI(n-1)|}_{(10)}\\
    &<\exp_2\left[(1-c){n\choose k-1}\right]|\cI(n-1)|
\end{align*}
(recall $c\lll\zeta\lll\theta_1\lll 1/k$). 
These bounds come from the following observations:
\begin{enumerate}
    \item[(2)] Since $T_{-u}$ is uniquely determined by $S_u,S_{-u}$, there are at most $4^{\binom{n}{k-2}}$ possible triples $S_u,S_{\ol u},T_{-u}$.
    \item[(3)(4)] For every $(k-2)$-set in $T_{-u}$, there are at most $2^n$ ways to choose the $(k-1)$-sets in $N_u$  containing it. For every $(k-2)$-set  in $S_{u}$, the number of $(k-1)$-sets in $N_u$ containing it is at most $\theta_1n$, so there are at most $\binom{n}{\leq\theta_1n}$ ways to choose these sets. Hence the number of possible choices for (3) is at most $ 2^{|T_{-u}|n}{n\choose \le \theta_1n  }^{|S_u|}$. By a similar reasoning, the number of possible choices for (4) is at most $ 2^{|T_{-u}|n}{n\choose \le \theta_1n  }^{|S_{\ol u}|}$. Since $|T_{-u}|\leq\theta_1n^{k-2}$, the number of choices for (3)(4) is at most
    \begin{align*}
        2^{2|T_{-u}|n}{n\choose \le \theta_1n  }^{|S_u+S_{\ol u}|}&\leq \exp_2[2\theta_1n^{k-1}+2H(\theta_1)n^{k-1}].
    \end{align*}

\item[(5)(6)] The candidate elements for (5) are those $(k-1)$ subsets of $X_{-u}$  whose $(k-2)$-subsets all lie in $S_{\ol u}\setminus S_u$. By the Kruskal--Katona theorem (\cref{thm:kruskalkatona}), the number of such candidates is at most $((k-2)!|S_u\setminus S_{\ol u}|)^\frac{k-1}{k-2}/(k-1)!$. By a similar reasoning, the number of candidates for (6) is at most $((k-2)!|S_{\ol u}\setminus S_u|)^\frac{k-1}{k-2}/(k-1)!$. We know from the claim hypothesis that $|S_u| \geq \theta_1 n^{k-2}$. 
Moreover, since $G$ is $\frac{1}{10k}$-dense and positive-dominant, there are at least $\frac{1}{20k}\binom{n-1}{k-1}$ clauses in $G$ that contains the positive literal $u$. By \cref{l:step2.1}(b), at most $(k-1)\zeta n^{k-1}$ of these clauses are non-monotone, so we have $|N_u|\geq \frac{1}{20k}\binom{n-1}{k-1}-(k-1)\zeta n^{k-1}$. Since 
$$|N_u|\leq \frac{1}{k-1}\left(\theta_1n|S_u|+n\left(\binom{n-1}{k-2}-|S_u|\right)\right)$$
with $\zeta\lll\theta_1\lll1/k$, we also have the upper bound $|S_u|\leq \binom{n-1}{k-2}-\theta_1 n^{k-2}$. Thus,
we have the upper bound
\begin{align*}
 &\frac{((k-2)!|S_u\setminus S_{\ol u}|)^\frac{k-1}{k-2}}{(k-1)!}+\frac{((k-2)!|S_{\ol u}\setminus S_u|)^\frac{k-1}{k-2}}{(k-1)!} \\
&=  \frac{((k-2)!)^\frac{k-1}{k-2}}{(k-1)!}\left(|S_u\setminus S_{\ol u}|^\frac{k-1}{k-2}+|S_{\ol u}\setminus S_u|^\frac{k-1}{k-2} \right)\\
&\overset{(*)}\le \frac{((k-2)!)^\frac{k-1}{k-2}}{(k-1)!} \left(\theta_1^\frac{k-1}{k-2}+\left(\frac{1}{(k-2)!}-\theta_1\right)^\frac{k-1}{k-2}\right)n^{k-1},
\end{align*}
where (*) follows from the convexity of $x\mapsto x^{(k-1)/(k-2)}$ and the fact that $\theta_1 n^{k-2}\leq |S_u|\leq  \binom{n-1}{k-2}-\theta_1 n^{k-2}$:
\begin{align*}
    |S_u\setminus S_{\ol u}|^\frac{k-1}{k-2}+|S_{\ol u}\setminus S_u|^\frac{k-1}{k-2}&\leq |S_u|^\frac{k-1}{k-2}+\left({n\choose k-2}-|S_u|\right)^\frac{k-1}{k-2}\\
    &\leq \left(\theta_1^\frac{k-1}{k-2}+\left(\frac{1}{(k-2)!}-\theta_1\right)^\frac{k-1}{k-2}\right)n^{k-1}.
\end{align*}
Thus the number of choices for (5)(6) is at most
\begin{align*}
\exp_2\sqb{\frac{((k-2)!)^\frac{k-1}{k-2}}{(k-1)!} \left(\theta_1^\frac{k-1}{k-2}+\left(\frac{1}{(k-2)!}-\theta_1\right)^\frac{k-1}{k-2}\right)n^{k-1}} &\le \exp_2\left[\left(\theta_1^{\frac{k-1}{k-2}} + \frac{1}{(k-1)!} \right) n^{k-1} \right],
\end{align*}
since $\frac{((k-2)!)^\frac{k-1}{k-2}}{(k-1)!} \le \frac12$ for $k \ge 3$.

    \item[(9)] By \cref{lem:ci_1starproperties}(b), for every $j\in\{1,\dots,k-1\}$, there are at most $\zeta n^{k-1}$ clauses in $G$ containing the positive literal $u$ and exactly $j$ negative literals; hence there are at most
    $\binom{{2n-2\choose k-1}}{ \le \zeta n^{k-1}}^{k-1}$ ways to determine the non-monotone clauses in $G$ that contain the positive literal $u$. Similarly, by \cref{lem:ci_1starproperties}(c), there are at most
\begin{align*}
    {{2n-2\choose k-1}\choose \le \zeta n^{k-1}}^{k-1} &\le \exp_2\left[{2n-2\choose k-1} H\left(\frac{\zeta n^{k-1}}{{2n-2\choose k-1}}\right)(k-1) \right] \tag{Equation~\ref{eq:bin-ent}} \\
    &\le \exp_2\left[\frac{1}{k} {n \choose k-1}\right] \tag{$\zeta \lll 1/k$}
\end{align*}
 ways to determine the clauses in $G$ that contain $\ol u$ and some other negative literal. \qedhere
\end{enumerate}
\end{proof*}

Combining \cref{claim:I21-case1} and \cref{claim:I21-case2}, we get that
\[
|\cI^*_{1}(n,\zeta)\setminus \cI_{2,1}^*(n, \zeta,\beta_1)|\leq \exp_2\sqb{(1-c){n\choose k}}+\exp_2\sqb{(1-c){n\choose k-1}}|\cI (n-1)|.
\qedhere \] 
\end{proof}

\subsubsection{Formulae with lots of clauses $\overline{u_1} \cdots \ol u_i v_1 \cdots v_{k-i}$}
Throughout this subsection, we assume $k \ge 4$.
Next, we upper bound the size of $\cI_{2,1}^*(n, \zeta,\beta_1)\setminus \cI_{2}^*(n, \zeta,\vec\beta)$. Recall from \cref{def:I2}  that 
for every $i\in\{2,\dots,k-2\}$ and $\beta_i>0$, we have
\begin{multline*}
    \cI^*_{2,i}(n, \zeta,\beta_1,\beta_i)=\bigl\{G\in \cI_{2,1}^*(n, \zeta,\beta_1):\text{for all distinct variables } u_1,\dots,u_i \in X,  \\ \text{$G$ contains at most $\beta_in^{k-i}$ clauses of the form } \ol u_1\cdots\ol u_iv_1\cdots v_{k-i}\bigr\},
\end{multline*}
and $\cI_{2}^*(n, \zeta,\vec\beta)=\bigcap_{k=2}^{k-2} \cI^*_{2,i}(n, \zeta,\beta_1,\beta_i)$.

We give upper bound on the size of  $\cI^*_{2,1}(n, \zeta,\beta_1)\setminus \cI_{2,i}^*(n, \zeta,\beta_1,\beta_i)$ for every $i\in\{2,\dots,k-2\}$.

\begin{lemma}\label{l:step2.2}
For $k \ge 4$ and every $i\in\{2,\dots,k-2\}$ and $\beta_i>0$, if positive $\zeta$ and $\beta_1$ are chosen sufficiently small relative to $\beta_i$, there exists $c > 0$ such that
$$|\cI^*_{2,1}(n, \zeta,\beta_1)\setminus \cI_{2,i}^*(n, \zeta,\beta_1,\beta_i)|\leq \exp_2\left[(1-c){n\choose k}\right]+\exp_2\left[i(1-c){n\choose k-1}\right]|\cI(n-i)|.$$
\end{lemma}
\begin{proof}
For every $i\in\{2,\dots,k-2\}$ and $\beta_i>0$, choose constants $0<c\lll\zeta\lll\theta_i\lll\beta_i \lll 1/k$ and $ \beta_1\lll\beta_i$. We will follow a similar proof strategy to that used in~\cref{l:step2.1}.

For every $G\in \cI_{2,1}^*(n,\zeta,\beta_1)\setminus \cI_{2,i}^*(n,\zeta,\beta_1,\beta_i)$, by definition, there exist some $u_1,\dots,u_i\in X$ such that there are more than $\beta_in^{k-i}$ clauses of the form $\ol u_1\cdots\ol u_iv_1\cdots v_{k-i}$. Let $X_{-\mathbf u}=X\setminus\{u_1,\dots,u_i\}$ and define the associated $(k-i)$-graph 
\begin{align*}
        N_{\ol{\mathbf u}}&=\left\{\{v_1,\dots,v_{k-i}\}\subseteq X_{-\mathbf u}:\ol u_1\cdots\ol u_iv_1\cdots v_{k-i}\in G\right\},
\end{align*}
so that by our assumption on $G$, $|N_{\ol{\mathbf u}}|>\beta_in^{k-i}$.
For every $\mathbf v=\{v_1,\dots,v_{k-i}\}\subseteq X_{-\mathbf u}$, let $N_{\mathbf u,\mathbf v}$ be the $(i-1)$-graph comprising the common co-neighborhood of sets $\{u_j, v_1, \ldots, v_{k-i}\}$ over all $j \in [i].$ Formally, let
\begin{align*}
    N_{\mathbf u,\mathbf v}&=\left\{\{w_1,\dots,w_{i-1}\}\subseteq X_{-\mathbf u}:u_jv_1\cdots v_{k-i}w_1\cdots w_{i-1}\in G\text{ for all }j\in[i]\right\},\\
    A_{\mathbf u}&=\left\{\mathbf v\in N_{ \ol {\mathbf u}}:\left|N_{\mathbf u,\mathbf v}\right|\geq\theta_in^{i-1}\right\}, \\
    B_{\mathbf u}&=N_{\ol{\mathbf u}}\setminus A_{\mathbf u}.
\end{align*}
Since $|N_{\ol{\mathbf u}}|>\beta_in^{k-i}$, at least one of the following holds:
\begin{enumerate}
    \item $|A_{\mathbf u}|>\theta_in^{k-i}$;
    \item $|B_{\mathbf u}|>(\beta_i-\theta_i)n^{k-i}$.
\end{enumerate}

We upper bound the number of formulae in $\cI_{2,1}^*(n,\zeta,\beta_1)\setminus \cI_{2,i}^*(n,\zeta,\beta_1,\beta_i)$ under each of these two cases.
\begin{claim}\label{claim:I2i-case1}
There are at most $\exp_2\left[(1-c){n\choose k}\right]$ formulae $G \in \cI_{2,1}^*(n,\zeta,\beta_1)\setminus \cI_{2,i}^*(n,\zeta,\beta_1,\beta_i)$ such that the following holds for some $u_1,\dots,u_k\in X$: (1) $G$ contains more than $\beta_in^{k-i}$ clauses of the form $\ol u_1\cdots\ol u_iv_1\cdots v_{k-i}$; (2) $|A_{\mathbf u}|>\theta_in^{k-i}$.

\end{claim}
\begin{proof*}
Every such $G$ is uniquely determined by sequentially specifying the following information:
\begin{enumerate}
    \item $u_1,\dots,u_i\in X$;
    \item $A_{\mathbf u}\in\binom{X_{-\mathbf u}}{k-i}$;
    \item $N_{\mathbf u,\mathbf v}$ for every $\mathbf v\in A_{\mathbf u}$;
    \item The non-monotone clauses in $G$;
    \item The monotone clauses in $G$.
\end{enumerate}
We prove the claim by giving an upper bound on the number of choices for each item:

\begin{align*}
    &\underbrace{\binom{n}{i}}_{(1)}
    \underbrace{2^{\binom{n}{k-i}}}_{(2)}
    \underbrace{2^{n^{k-1}}}_{(3)}
    \underbrace{\binom{2^k\binom{n}{k}}{\leq \zeta n^k}}_{(4)}    \underbrace{\exp_2\sqb{{n\choose k}-\Omega(n^k)}}_{(5)} \leq \exp_2\left[(1-c){n\choose k}\right]
\end{align*}
(recall that $c\lll \zeta\lll\theta_i$). The bounds come from the following observations:
\begin{enumerate}
    \item[(3)] For every fixed $\mathbf v\in A_{\mathbf u}$, there are $2^{{n\choose {i-1}}}$ possible choices for $N_{\mathbf u,\mathbf v}$. Since $|A_{\mathbf u}|\leq\binom{n}{k-i}$, the number of possible choices for (3) is at most $2^{\binom{n}{i-1}\binom{n}{k-i}}\leq 2^{n^{k-1}}$.
    \item[(4)] Since $G \in \cI_1^*(n, \zeta)$, it has at most $\zeta n^k$ non-monotone clauses.
    \item[(5)]  For every $\mathbf v=\{v_1,\dots,v_{k-i}\}\in A_{\ol {\mathbf u}}$ and $\{w_1,\dots,w_{i-1}\}\in N_{\mathbf u,\mathbf v}$, to avoid the non-minimal subformula 
\begin{align*}
    &\{v_1\cdots v_{k-i}w_1\cdots w_{i-1}z, \\
    &\quad v_1\cdots v_{k-i}w_1\cdots w_{i-1}u_1,\,\,\ldots,\,\,v_1\cdots v_{k-i}w_1\cdots w_{i-1}u_i,\\
    &\hspace{6cm} \quad \ol u_1\cdots\ol u_iv_1\cdots v_{k-i}\}
\end{align*}
(as it is impossible to satisfy only the first clause), we cannot have $v_1\cdots v_{k-i}w_1\cdots w_{i-1}z\in G$ for any $z\in X$. Therefore, (1)--(3) have determined at least
$$\frac{1}{i{k\choose i}}\cdot |A_{\mathbf u}|\theta_in^{i-1}(n-(k-1))\geq \frac{1}{i{k\choose i}}\cdot \theta_in^{k-i}\theta_in^{i-1} (n-(k-1))\geq \frac{\theta_i^2 n^{k}}{i{k\choose i}} = \Omega(n^k)$$
monotone clauses that cannot belong to $G$. \qedhere
\end{enumerate}
\end{proof*}

\begin{claim}\label{claim:I2i-case2}
There are at most $\exp_2[i(1-c){n\choose k-1}]|\cI(n-i)|$ formulae $G \in \cI_{2,1}^*(n,\zeta,\beta_1)\setminus \cI_{2,i}^*(n,\zeta,\beta_1,\beta_i)$ such that the following holds for some $u_1,\dots,u_k\in X$: (1) $G$ contains more than $\beta_in^{k-i}$ clauses of the form $\ol u_1\cdots\ol u_iv_1\cdots v_{k-i}$; (2) $|B_{\mathbf u}|>(\beta_i-\theta_i)n^{k-i}$.

\end{claim}
\begin{proof*}
Every such $G$ is uniquely determined by sequentially specifying the following information:
\begin{enumerate}
    \item $u_1,\dots,u_i\in X$;
    \item $B_{\mathbf u}\in\binom{X_{-\mathbf u}}{k-i}$;
    \item $N_{\mathbf u,\mathbf v}$ for every $\mathbf v\in B_{\mathbf u}$;
    \item The monotone clauses in $G$ that use at least one of the variables $u_1,\dots,u_i$;
    \item The clauses in $G$ that contain one of the negative literals $\ol {u_1},\dots,\ol{u_i}$ and $k-1$ other positive literals;
    \item The clauses in $G$ that  (i) use at least one of the variables $u_1,\dots,u_i$, and (ii) contain a negative literal that is different from the negation of this variable;
    \item The clauses in $G$ that do not use any of the variables $u_1,\dots,u_i$.
\end{enumerate}
We prove the claim by giving an upper bound on the number of choices for each item:
\begin{align*}
    &\underbrace{n}_{(1)}\underbrace{2^{\binom{n}{k-i}}}_{(2)}\underbrace{\binom{\binom{n}{i-1}}{\leq \theta_in^{i-1}}^{n^{k-i}}}_{(3)}\underbrace{\exp_2\left[i{n-1\choose k-1}-(i-\log(2^i-1))\Omega(\beta_i)n^{k-1}\right]}_{(4)}\\
    &\qquad \underbrace{\binom{\binom{n-1}{k-1}}{\leq \beta_1n^{k-1}}^i}_{(5)}\underbrace{\binom{2^{k-1}\binom{n-1}{k-1}}{\leq k\zeta n^{k-1}}^{2i}}_{(6)}\underbrace{|\cI(n-i)|}_{(7)}\\&< \exp_2\sqb{i(1-c){n\choose k-1}}|\cI(n-i)|
\end{align*}
(recall that $c\lll \zeta\lll\theta_i\lll\beta_i$ and $\beta_1\lll\beta_i$). The bounds come from the following observations:
\begin{enumerate}
    \item [(3)]For every fixed $\mathbf v\in B_{\mathbf u}$, since $|N_{\mathbf u,\mathbf v}|<\theta_in^{i-1}$, there are at most $\binom{\binom{n}{i-1}}{\leq \theta_in^{i-1}}^{n^{k-i}}$ possible choices for $N_{\mathbf u,\mathbf v}$. 
    \item[(4)] Let
    \begin{align*}
        &Q_{\mathbf u}=\{\{v_1\cdots v_{k-i}w_1\cdots w_{i-1}\}\in\binom{X_{-\mathbf u}}{k-1}:\\
        &\hspace{5cm}\mathbf v=\{v_1,\dots,v_{k-i}\}\in B_{\mathbf u},\, \{w_1,\dots,w_{i-1}\}\notin N_{\mathbf u,\mathbf v}\}.
    \end{align*}
    By definition of $N_{\mathbf u,\mathbf v}$, for all $\{v_1\cdots v_{k-i}w_1\cdots w_{i-1}\}\in Q_{\mathbf u}$, there exists some $j\in[i]$ such that $u_jv_1\cdots v_{k-i}w_1\cdots w_{i-1}\notin G$.
    Since $|B_{\mathbf u}|>(\beta_i-\theta_i)n^{k-i}$ and $|N_{\mathbf u,\mathbf v}|<\theta_in^{i-1}$ for every $\mathbf v\in B_{\mathbf u}$, (1)--(3) have specified at least $${k-1\choose i-1}^{-1} (\beta_i-\theta_i)n^{k-i} \left({n-1-(k-i)\choose i-1}-\theta_in^{i-1}\right)\geq \Omega(\beta_i)n^{k-1} $$
    $(k-1)$-subsets of $X$ that lie in $Q_{\mathbf u}$ (recall that $\theta_i\lll \beta_i$) . For every $\{v_1\cdots v_{k-i}w_1\cdots w_{i-1}\}\in Q_{\mathbf u}$, there are at most $2^i-1$ ways to choose which of $u_1v_1\cdots v_{k-i}w_1\cdots w_{i-1}$, \dots,\\ $u_iv_1\cdots v_{k-i}w_1\cdots w_{i-1}$ lie in $G$. Therefore, there are at most 
\begin{align*}
    \exp_2\left[i{n-1\choose k-1}-(i-\log(2^i-1))\Omega(\beta_i)n^{k-1}\right]
\end{align*}
choices for the monotone clauses in $G$ that contain some of $u_1,\dots,u_i$.
\item[(5)] For every $j\in[i]$, since $G\in \cI_{2,1}^*(n, \zeta,\beta_1)$, there are at most $\beta_1n^{k-1}$ clauses in $G$ that contain $\ol{u_j}$ and $k-1$ positive literals.
\item[(6)] For every variable $u_j$, by \cref{lem:ci_1starproperties}(b)(c), there are at most $k\zeta n^{k-1}$ clauses in $G$ that contain the positive literal $u_j$ and some negative literal; there are at most $k\zeta n^{k-1}$ clauses in $G$ that contain the negative literal $\ol{u_j}$ and some other negative literal. Hence there are at most $\binom{2^{k-1}\binom{n-1}{k-1}}{\leq k\zeta n^{k-1}}^{2i}$ choices for (6). \qedhere
\end{enumerate}
\end{proof*}

Combining \cref{claim:I2i-case1} and \cref{claim:I2i-case2}, we get that
\[
|\cI^*_{2,1}(n,\zeta,\beta_1)\setminus \cI_{2,i}^*(n, \zeta,\beta_1,\beta_i)|\leq \exp_2\sqb{(1-c){n\choose k}}+\exp_2\sqb{i(1-c){n\choose k-1}}|\cI(n-i)|.\qedhere
\]
\end{proof}

\subsection{$k$-wise boolean combinations of monotone neighborhoods are typically not too small} It remains to upper bound the number of formulae in $\cI_2^*(n, \zeta,\vec\beta)$.  
In this section, we make some preparations for giving this bound.

For every $k$ variables $u_1,\dots,u_k\in X$,  define
\begin{align*}
    X_{-\mathbf u}&=X\setminus \{u_1,\dots,u_k\},\\
    N_{u_j}^{-\mathbf u}&=\{\{v_1, \ldots v_{k-1}\} \subseteq X_{-\mathbf u} : u_j v_1\cdots v_{k-1} \in G\},\\
    \ol{N_{u_j}^{-\mathbf u}}&=\binom{X_{-\mathbf u}}{k-1}\setminus \{\{v_1, \ldots v_{k-1}\} \subseteq X_{-\mathbf u} : u_jv_1\cdots v_{k-1} \in G\}.
\end{align*}

We say that vectors $\vec a, \vec b$ have $\vec a \succ \vec b$, if $a$'s entries are component-wise larger than $b$'s entries.

\begin{definition}\label{def:I3}
For   $\vec \beta \succ 0$, let
\begin{align*}
\cI_3^*(n, \zeta, \vec\beta) &=\bigg\{G \in \cI_2^*(n, \zeta, \vec\beta) : \text{ for all } u_1, u_2, \ldots, u_k \in X,\\
&\hspace{1cm}M_{u_1}\in\{N_{u_1}^{-\mathbf u},\ol{N_{u_1}^{-\mathbf u}}\},\dots, M_{u_k}\in\{N_{u_1}^{-\mathbf u},\ol{N_{u_k}^{-\mathbf u}}\},\\
&\hspace{1.75cm} \text{we have } |M_{u_1}\cap\dots\cap M_{u_k}|\geq \frac{1}{2^{k+1}}\binom{n-k}{k-1} \bigg\}.
\end{align*}

\end{definition}

\begin{lemma} \label{l:step3}
For all $k \ge 2$ and all sufficiently small $\vec \beta \succ 0$, if $\zeta$ is sufficiently small relative to $\beta_1$, there exists $c > 0$ such that 
$$|\cI_2^*(n, \zeta, \vec\beta) \backslash \cI_3^*(n, \zeta, \vec\beta)| \le \exp_2\left[\frac{k}{2}(1-c){n\choose k-1}\right]|\cI(n-k)|.$$
\end{lemma}
\begin{proof}Fix $\vec \beta \succ 0$ sufficiently small and choose $0<c\lll \zeta\lll \beta_1$.  For every formula $G\in\cI_2^*(n, \zeta, \vec\beta) \backslash \cI_3^*(n, \zeta, \vec\beta)$, there exist variables $u_1,\dots,u_k\in X$ 
and $M_{u_1}\in\{N_{u_1}^{-\mathbf u},\ol{N_{u_1}^{-\mathbf u}}\},\dots, M_{u_k}\in\{N_{u_1}^{-\mathbf u},\ol{N_{u_k}^{-\mathbf u}}\}$
such that $|M_{u_1} \cap \cdots\cap M_{u_k}|< \frac{1}{2^{k+1}} {n-k \choose k-1}. $

Every formula $G\in\cI_2^*(n, \zeta, \vec\beta) \backslash \cI_3^*(n, \zeta, \vec\beta)$ is uniquely determined by sequentially specifying the following information.
\begin{enumerate}
    \item $u_1, \ldots, u_k$;
    \item The non-monotone clauses in $G$ that use at least one of the variable $u_1, \ldots, u_k$.
    \item The monotone clauses in $G$ that use at least one of the variable $u_1, \ldots, u_k$.
    \item The clauses in $G$ that do not use any of the variables $u_1, \ldots, u_k$.
\end{enumerate}
We upper bound $|\cI_2^*(n, \zeta, \vec\beta) \backslash \cI_3^*(n, \zeta, \vec\beta)|$ by upper bounding the number of choices for each item:
\begin{align*}
    |\cI_2^*(n, \zeta, \vec\beta) \backslash \cI_3^*(n, \zeta, \vec\beta)|&\leq \underbrace{ {n\choose k}}_{(1)} \underbrace{{(2^k - 1){ n - 1 \choose k-1} \choose \le \beta_1n^{k-1} + 2(k-1)\zeta n^{k-1}}^k}_{(2)} \underbrace{\exp_2\left[H\left(\frac{1}{2^{k}}\right){n-k \choose k-1}  \right]}_{(3)}\underbrace{|\cI(n-k)|}_{(4)} \\
    &\leq \exp_2\left[\frac{k}{2}(1-c){n\choose k-1}\right]|\cI(n-k)|
\end{align*}
(recall that $c\lll \zeta\lll\beta_1$). In particular, upper bounds for (2) and (3) arise from the following observations.
\begin{enumerate}
    \item[(2)] For every variable $u_j$, by \cref{lem:ci_1starproperties}(b)(c), $G$ has at most $2(k-1)\zeta n^{k-1}$ non-monotone clauses that use the variable $u_j$ and contain some negative literal other than $\ol{u_j}$. Moreover, since $G\in \cI_{2,1}^*(n,\zeta,\beta_1)$, $G$ has at most $\beta_1n^{k-1}$ clauses of the form $\ol{u_j}v_1\cdots v_{k-1}$. 
    \item[(3)] There are $2^{O(n^{k-2})}$ possibilities to choose the monotone clauses in $G$ that use at least two of the variables $u_1,\ldots, u_k$. The monotone clauses in $G$ that use exactly one of the variables $u_1,\dots,u_k$ is completely determined by
    \begin{align*}
        \{ M_{u_1}\cap \dots \cap M_{u_k}: M_{u_1}\in\{N_{u_1}^{-\mathbf u},\ol{N_{u_1}^{-\mathbf u}}\},\dots, M_{u_k}\in\{N_{u_1}^{-\mathbf u},\ol{N_{u_k}^{-\mathbf u}}\}\},
    \end{align*}
     which form an ordered partition of $\binom{X_{-\mathbf u}}{k-1}$ into $2^k$ parts. Since $G\notin \cI_3^*(n,\zeta,\vec\beta)$, at least one of these parts has size less than $\binom{n-k}{k-1}/2^{k+1}$. The number of such ordered partitions is at most 
\[
2^k\exp_2\left[H\left(\frac{1}{2^{k+1}}\right){n-k \choose k-1}\right]\binom{{n-k \choose k-1}+2^{k}-2}{2^{k}-2}
         \leq \exp_2\left[H\left(\frac{1}{2^{k}}\right){n-k \choose k-1}  \right]. \qedhere 
\]
\end{enumerate}
\end{proof}

\subsection{Most formulae in $\cI_3^*(n, \zeta, \vec \beta)$ are monotone}\label{s:endik}
In this step, we bound the number of formulae in $\cI_{3}^*(n,\zeta,\vec\beta)$ that are not monotone.

\begin{definition}Let
$$\cI_{4}^*(n,\zeta,\vec\beta)=\{G\in \cI_3^*(n, \zeta,\vec\beta):\text{$G$ is monotone}\}.$$
\end{definition}

We give an upper bound on $|\cI_{3}^*(n,\zeta,\vec\beta)\setminus \cI_{4}^*(n,\zeta,\vec\beta)|$. 

\begin{lemma}\label{l:step4-k=2}
For $k = 2$, for $n$ sufficiently large $\cI_{3}^*(n,\zeta,\vec\beta)\setminus \cI_{4}^*(n,\zeta,\vec\beta) = \emptyset$.
\end{lemma}
\begin{proof}
We let $n$ be sufficiently large.
Fix $G \in \mcI_3^*(n,\zeta,\vec\beta)$. For every pair of variables $u, v$, we know from \cref{def:I3} that $N_u \cap N_v \neq \emptyset$.

Suppose to the contrary that $G$ has a non-monotone clause. If some $\ol u v \in G$, then for any $w \in N_u \cap N_v$, $\{\ol u v, u w, v w\}$ certifies the non-minimality of $G$, a contradiction. 

Else, if some $\ol u \ol v \in G$, we can again select some $w \in N_u \cap N_v$. Since $\mcI_3^*(n,\zeta,\vec\beta) \subseteq \mcI_1^*(n,\zeta)$, $w$ is used by at least $n/20-2\zeta n$ clauses of the form $zw$ or $\ol z w$ in $G$. Thus  we can find some other variable $z \not \in \{u, v, w\}$ such that  $z w $ or $\ol zw$ lies in $G$. Since either $\{\ol u \ol v, u w, vw, zw\}$ or $\{\ol u \ol v, u w, vw, \ol zw\}$ certifies the non-minimality of $G$, we also get a contradiction.
\end{proof}
Hereafter, we may assume that $k \ge 3$ as $\cI_{3}^*(n,\zeta,\vec\beta)\setminus \cI_{4}^*(n,\zeta,\vec\beta) = \emptyset$ for $k = 2$.
The following observation about $\ell$-hypergraphs generalizes a theorem on graphs in \cite{Hak65} (see also \cite[Theorem 61.1]{Sch03}), and will be useful in our proof.

\begin{lemma}
\label{lem:orientation} For $\ell \ge 2$, every $\ell$-uniform hypergraph $H$ can be oriented (i.e., every edge directed at some vertex) such that for every $\ell-1$ vertices $v_1,\dots,v_{\ell-1}$, there are at most $$\ceil{\frac{(\ell-1)!}{(\ell!)^\frac{\ell-1}{\ell}}\cdot e(H)^\frac{1}{\ell}}$$ edges that contain all of $v_1,\dots,v_{\ell-1}$ and are directed at a vertex $w$ other than $v_1,\dots,v_{\ell-1}$.
\end{lemma}
\begin{proof}
Consider an $\ell$-graph $H$. Write
$$L=\ceil{\frac{(\ell-1)!}{(\ell!)^\frac{\ell-1}{\ell}}\cdot e(H)^\frac{1}{\ell}}.$$
Construct an associated bipartite graph $B_H$ in which one part is $E(H)$ (the edge set of $H$), and the other part is the multiset $\bigcup_{S\in {V(H)\choose \ell-1}}\{\!\!\{\underbrace{S,\dots,S}_{L\text{ times}}\}\!\!\}$. For every $e\in E(H)$ and $S\in {V(H)\choose \ell-1}$, we define that $e\sim S$ in $B_H$ if and only if $S\subseteq e$.

Observe that
\begin{align*}
    &\text{$B_H$ has an $E(H)$-matching} \\& \Longleftrightarrow \text{$H$ has an orientation in which every $\ell-1$ vertices are contained}\\
    &\quad\quad \text{in at most $L$ edges directed at another vertex.}
\end{align*}
Thus, it suffices to show that $B_H$ has an $E(H)$-matching.
By Hall's theorem, it is equivalent to show that 
$|W|\leq|N_{B_H}(W)|$ for every $W\subseteq E(H)$.

For every $W\subseteq E(H)$, 
$N_{B_H}(W)$ consists of all $S\in {V(H)\choose \ell-1}$ (each repeated $L$ times) that are contained in some $e\in W$. By the Kruskal--Katona theorem (\cref{thm:kruskalkatona}), every collection $\cS\subseteq {V(H)\choose \ell-1}$ covers at most $((\ell-1)!\abs{\cS})^\frac{\ell}{\ell-1}/\ell!$ edges in $H$, so every $W\subseteq E(H)$ satisfies
\begin{align*}
    |W|&\leq \left((\ell-1)! |N_{B_H}(W)|/L\right)^\frac{\ell}{\ell-1}/\ell!.
\end{align*}
Rearranging gives
\begin{align*}
    |N_{B_H}(W)|&\geq\frac{(\ell!)^\frac{\ell-1}{\ell}|W|^\frac{\ell-1}{\ell}L}{(\ell-1)!}\geq  \frac{(\ell!)^\frac{\ell-1}{\ell}|W|^\frac{\ell-1}{\ell}}{(\ell-1)!}\cdot \frac{(\ell-1)!}{(\ell!)^\frac{\ell-1}{\ell}}\cdot e(H)^\frac{1}{\ell}
    = |W|^\frac{\ell-1}{\ell}\cdot e(H)^\frac{1}{\ell}\geq |W|.\qedhere
\end{align*}
\end{proof}

We will also use the following version of Shearer's inequality \cite{CGFS86} (see also \cite[Lemma 6.5]{IK12}). 
\begin{lemma}\label{lem:shearer}
Let $W$ be a set and $\cF$ be a family of subsets of $W$. Let $\cH$ be a hypergraph with  $V(\cH)=W$ and $\deg_{\cH}(v)\geq k$ for every $v\in W$. Then
$$\log_2|\cF|\leq\frac{1}{k}\sum_{e\in E(\cH)}\log_2|\Tr(\cF,e)|,$$
where  $\Tr(\cF,e)=\{F\cap e:F\in\cF\}$.
\end{lemma}

Next, we verify the following technical result.

\begin{lemma}\label{claim:step4final}
Suppose $\vec \beta, \zeta$ are all sufficiently small constants (relative to $k$). Then, there exists $c > 0$ such that
for all $0\leq t_1,\dots,t_{k-2}\leq \zeta n^k$, we have
\begin{align*}
    & \left|\left\{G\in \cI_3^*(n, \zeta,\vec\beta)\setminus \cI_4^*(n, \zeta,\vec\beta):\text{$G$ contains no clause with exactly $k$ negative literals, }\right.\right.\\
    &\left.\left.\qquad\qquad\qquad\qquad\qquad\qquad\qquad |G_1|=t_1,\dots, |G_{k-2}|=t_{k-2} \right\}\right|\\
    & <\exp_2\left[{n\choose k}-2cn\right],
\end{align*}
where $G_i$ is the set of clauses in $G$ with exactly $i$ negative literals (so $G$ is the disjoint union of $G_0,G_1,\dots,G_{k-2}$).
\end{lemma}
The proof of this lemma is rather involved, so below we give a brief informal outline of the general argument:
\begin{itemize}
    \item We consider some $k$-SAT formula $G$ satisfying the lemma conditions for a particular choice of parameters $(t_1, \ldots, t_{k-2}).$
    \item We partition the clauses of $G$ based on the number of negative literals $i$, and for each $i$ we construct an auxiliary oriented hypergraph $H_i$.
    \item Each edge of our hypergraph $H_i$ encodes a constraint on which clauses can and cannot appear in $G$.
    \item Thus, our primary proof strategy is to enumerate the number of such $G$ by counting (a) the number of possible hypergraph sequences $\{H_i\}_{i \in [k-2]}$ and (b) the number of possible formulae that correspond to a given hypergraph sequence.
\end{itemize}

\begin{proof}[Proof of~\cref{claim:step4final}]
Fix $0\leq t_1,\dots,t_{k-2}\leq \zeta n^k$. Consider some $G\in \cI_3^*(n, \zeta,\vec\beta)\setminus \cI_4^*(n, \zeta,\vec\beta)$ such that $G$ does not contain clauses with exactly $k$ negative literals, and $|G_1|=t_1$,\dots, $|G_{k-2}|=t_{k-2}$.

For every $i\in[k-2]$, we construct an auxiliary collection $G_i'\subseteq {X\choose i}\times {X\choose k-i-1}\times X$ based on $G_i$, by doing the following: for every clause $\ol u_1\cdots\ol u_iv_1\cdots v_{k-i}\in G_i$, we pick one $j\in\{1,\dots,k-i\}$, and put the element $(\{u_1,\dots,u_i\},\{v_1,\dots,v_{j-1},v_{j+1},\dots, v_{k-i}\},v_j)$ in $G_i'$. In other words, we obtain $G_i'$ from $G_i$ by making the set of positive literals in every clause ``directed'' at some $v_j$. 

There are clearly many ways to build $G_i'$. Using \cref{lem:orientation}, we show that there exists some $G_i'$ such that for every $u_1,\dots,u_i\in X$ and $v_1,\dots,v_{k-i-1}\in X\setminus\{u_1,\dots,u_i\}$, there are at most $$\ceil{\frac{(k-i-1)!}{((k-i)!)^\frac{k-i-1}{k-i}}\cdot \beta_i^\frac{1}{k-i}\cdot n}$$ elements of the form $(\{u_1,\dots,u_i\},\{v_1,\dots,v_{k-i-1}\},w)$. For every $u_1,\dots,u_i\in X$, define the $(k-i)$-hypergraph $H_{\mathbf u}$ by
\begin{align*}
    V(H_{\mathbf u})&=X\setminus\{u_1,\dots,u_i\},\\
    E(H_{\mathbf u})&=\{\{v_1,\dots,v_{k-i}\}:\ol u_1\cdots\ol u_iv_1\cdots v_{k-i}\in G_i\}.
\end{align*}
Since $G$ contains at most $\beta_in^{k-i}$ clauses of the form $\ol u_1\dots\ol u_iv_1\dots v_{k-i}$, we have $e(H_{\mathbf u})\leq \beta_in^{k-i}$.  By~\cref{lem:orientation},  $H_{\mathbf u}$ can be directed such that every $k-i-1$ vertices in $X\setminus\{u_1,\dots,u_i\}$ are contained in at most $\ceil{\frac{(k-i-1)!}{((k-i)!)^\frac{k-i-1}{k-i}}\cdot \beta_i^\frac{1}{k-i}\cdot n}$ edges as the undirected part. Combining orientations over all $X_{\mathbf u}$ gives us $G_i'$ as wanted.

Now  fix the above $G_i'$. For every $i\in\{1,\dots,k-2\}$, define the $(i+1)$-multigraph $H_i$ by
\begin{align*}
   V(H_i)&={X\choose k-i},\\
   E(H_i)&=\{\{\{u_1,v_1,\dots,v_{k-i-1}\},\dots, \{u_i,v_1,\dots,v_{k-i-1}\},\{w,v_1,\dots,v_{k-i-1}\}\}:
   \\&\hspace{6cm} (\{u_1,\dots,u_i\},\{v_1,\dots,v_{k-i-1}\},w)\in G_i'\},
\end{align*}
so that $|E(H_i)|=t_i$. We emphasize that $G_i$ is uniquely determined by $H_i$. Moreover, due to the orientation of $G_i'$, $H_i$ has small maximum degree.

\begin{claim}For every $i\in\{1,\dots,k-2\}$, we have
$$\Delta_1(H_i)\leq (k-i){n\choose i-1}\left(\ceil{\frac{(k-i-1)!}{((k-i)!)^\frac{k-i-1}{k-i}}\cdot \beta_i^\frac{1}{k-i}\cdot n}+\zeta n\right).$$
\end{claim}
\begin{proof*}
Fix some $\mathbf v\in\binom{X}{k-i}$.  Every edge  $e\in E(H_i)$ in which $\mathbf v$ takes the role of $\{u_1,v_1,\dots,v_{k-i-1}\}$ (resp. $\{u_2,v_1,\dots,v_{k-i-1}\}$, \dots, $\{u_i,v_1,\dots,v_{k-i-1}\}$) is uniquely determined by the following:
\begin{enumerate}
    \item[(i)] $u_1\in \mathbf v$;
    \item[(ii)] $u_2,\dots,u_i\in X$;
    \item[(iii)] $w$.
\end{enumerate}
Unlike in some earlier arguments in this section that look superficially similar, here, the underlying hypergraph $H_i$ is fixed, as is $G$. Thus, the number of such edges is just at most
\begin{align*}
   \underbrace{(k-i)}_{(i)}\underbrace{\binom{n}{i-1}}_{(ii)}\underbrace{\ceil{\frac{(k-i-1)!}{((k-i)!)^\frac{k-i-1}{k-i}}\cdot \beta_i^\frac{1}{k-i}\cdot n}}_{(iii)},
\end{align*}
where (iii) arises from the fact that given $u_1,\dots,u_i,v_1,\dots,v_{k-1-i}$, there are $\leq\ceil{\frac{(k-i-1)!}{((k-i)!)^\frac{k-i-1}{k-i}}\cdot \beta_i^\frac{1}{k-i}\cdot n}$ elements of the form $(\{u_1,\dots,u_i\},\{v_1,\dots,v_{k-i-1}\},w)$ in $G_i'$.
Meanwhile, every edge  $e\in E(H_i)$ in which $\mathbf v$ takes the role of $\{w,v_1,\dots,v_{k-i-1}\}$ is uniquely determined by  the following:
\begin{enumerate}
    \item[(i)] $w\in \mathbf v$;
    \item[(ii)] $u_1,\dots,u_{i-1}\in X$;
    \item[(iii)] $u_i$.
\end{enumerate}
The number of  such edges is at most
\begin{align*}
  \underbrace{(k-i)}_{(i)}\underbrace{\binom{n}{i-1}}_{(ii)}\underbrace{\zeta n}_{(iii)},
\end{align*}
where (iii) follows from the  fact that given $u_1,\dots,u_{i-1},v_1\dots v_{k-i-1},w$, by \cref{lem:ci_1starproperties}(a), there are at most $\zeta n$ variables $u_i$ in $G$ such that $\ol{u_1}\cdots \ol{u_i}v_1\cdots v_{k-1-i}w\in G$.
\end{proof*}

Let $\nu_i$ and $\tau_i$ denote the matching and vertex cover numbers of $H_i$. Suppose $0\lll \beta_i \lll \rho_i\lll 1$ and $0 \lll \zeta \lll 1$, so that we have the degree bound
$$\Delta_1(H_i)\leq (k-i){n\choose i-1}\left(\ceil{\frac{(k-i-1)!}{((k-i)!)^\frac{k-i-1}{k-i}}\cdot \beta_i^\frac{1}{k-i}\cdot n}+\zeta n\right)\leq\rho_in^i.$$
We thus have the  inequality (recall that $|E(H_i)|=t_i$) $$(i+1)\nu_i\geq\tau_i\geq \frac{t_i}{\rho_i n^i}.$$

We now give an upper bound on the number of choices for $H_i$ given $t_i$, $\nu_i$ and $\tau_i$.
\begin{claim}\label{claim:step4final-2}
    For every $i\in\{1,\dots, k-2\}$, given $t_i$, $\nu_i$ and $\tau_i$, the number of choices for $H_i$ is at most
    \[
\binom{\binom{n}{k-i}}{\tau_i}\binom{\tau_i\binom{n}{i}\binom{k}{i+1}}{t_i}.
\]
\end{claim}
\begin{proof*}Recall that $|V(H_i)|=|\binom{X}{k-i}|=|\binom{n}{k-i}|$. Observe that $H_i$ is uniquely determined by  the following:
\begin{itemize}
    \item[(1)] A vertex cover $\mathcal T_i$ of $H_i$;
    \item[(2)] The edge set $E(H_i)$ of size $t_i$, such that every edge contains at least one vertex in $\mathcal T_i$.
\end{itemize}
Given $t_i$, $\nu_i$ and $\tau_i$, there are $\binom{|V(H_i)|}{\tau_i}=\binom{\binom{n}{k-i}}{\tau_i}$ possible vertex covers of $H_i$. Moreover, any vertex cover of $H_i$ gives at most $\tau_i\binom{n}{i}\binom{k}{i+1}$ edges that could possibly lie in $H_i$; thus, it gives at most $\binom{\tau_i\binom{n}{i}\binom{k}{i+1}}{t_i}$ possible choices for $H_i$. Hence given $t_i$, $\nu_i$ and $\tau_i$, the number of possible choices for $H_i$ is at most
\[
\binom{\binom{n}{k-i}}{\tau_i}\binom{\tau_i\binom{n}{i}\binom{k}{i+1}}{t_i}.
\]
\end{proof*}

We now upper bound the number of choices for $G_0$ given $H_1,\dots,H_{k-2}$.
Observe that for every
$$e=\{\{u_1,v_1,\dots,v_{k-i-1}\},\dots, \{u_i,v_1,\dots,v_{k-i-1}\},\{w,v_1,\dots,v_{k-i-1}\}\}\in E(H_i)$$ 
and any other $i$ variables $a_1,\dots,a_i\notin\{u_1,\dots, u_i,v_1,\dots,v_{k-i-1},w\}$, since $\ol u_1\cdots \ol u_{i}v_1\cdots v_{k-i-1}w\in G_i$, to avoid the non-minimal subformula
\begin{align*}
    &\big\{v_1\cdots v_{k-i-1}a_1\cdots a_iw, \\
    &\hspace{0.3cm} v_{1}\cdots v_{k-i-1}a_1\cdots a_iu_1,\\
    &\hspace{0.3cm} \dots, \\
    &\hspace{0.6cm} v_1\cdots v_{k-i-1}a_1\cdots a_iu_i,\\
    &\hspace{0.7cm} \ol u_1\cdots \ol u_{i}v_1\cdots v_{k-i-1}w\big\}
\end{align*}
(as it is impossible to satisfy only the first clause), we have
$$|G_0\cap F_{e,a_1,\dots,a_i}|<i+1,$$
with $F_{e,a_1,\dots,a_i}$ defined by
$$F_{e,a_1,\dots,a_i}=\{v_1\cdots v_{k-i-1}a_1\cdots a_iw, \, v_{1}\cdots v_{k-i-1}a_1\cdots a_iu_1, \, \dots, \, v_1\cdots v_{k-i-1}a_1\cdots a_iu_i\}.$$

\begin{claim}\label{claim:step4final-3}
    Suppose $H_1,\dots,H_{k-2}$ have been determined. Then the number of possible choices for $G_0$ is at most
\[\exp_2\left[{n\choose k}-\sum_{i=1}^{k-2}\delta_i{n-k\choose i}\tau_i\right],\]
where
$$\delta_i= \frac{i+1-\log_2(2^{i+1}-1)}{(2^k-k-2)(i+1){k-1\choose i}}.$$
\end{claim}
\begin{proof*}
Let $\mathcal M_i$ be a  maximum matching of $H_i$, so that $|\mathcal M_i|=\nu_i$.
Define the family of monotone formulae $\mathcal J_i$ by
\begin{align*}
    \mathcal J_i&=\bigcup_{e\in \mathcal M_i}\{F_{e,a_1,\dots,a_i}:a_1,\dots,a_i\notin\{u_1,\dots, u_i,v_1,\dots,v_{k-i-1},w\}\\
    & \qquad\qquad\text{ with }e=\{\{u_1,v_1,\dots,v_{k-i-1}\},\dots, \{u_i,v_1,\dots,v_{k-i-1}\},\{w,v_1,\dots,v_{k-i-1}\}\}\}.
\end{align*}
For every possible  $G_0$ and every $F\in \mathcal J_i$, by the above argument, we have $|G_0\cap F|< i+1.$
Moreover, since every $e\in \mathcal M_i$ gives $\binom{n-k}{i}$ different edges in $\mathcal J_i$, and every   $F\in\mathcal J_i$ can arise from at most ${k-1\choose i}$ different  $e\in\mathcal M_i$, we have
$$|\mathcal J_i|\geq\frac{\nu_i{n-k\choose i}}{{k-1\choose i}}.$$

For every monotone clause $C\in\binom{X}{k}$, let $d_i(C)$ denote the number of elements $F\in\mathcal J_i$ that contain $C$. Note that we always have $d_i(C)\leq\binom{k}{i}$, or there would be two edges $e\in \mathcal M_i$ sharing the same vertex, contradicting the fact that $\mathcal M_i$ is a matching. Let $\mathcal H$ be a (non-uniform) multigraph on the set of monotone clauses, defined by
\begin{align*}
    V(\mcH)={X\choose k},\qquad
    E(\mcH)=(\bigcup_{i=1}^{k-2}\mathcal J_i)\cup(\bigcup_{\substack{C\in {X\choose k}}}\,\{\!\!\{\underbrace{\{C\},\dots,\{C\}}_{\sum_{i=1}^{k-2}(\binom{k}{i}-d_i(C))\text{ times}}\}\!\!\}\, ),
\end{align*}
so that every vertex in $\cH$ has degree  $\sum_{i=1}^{k-2}\binom{k}{i}=2^k-k-2$. Let $\mathcal G_{0}$ be the collection of all possible choices for $G_0$. By~\cref{lem:shearer}, we have
\begin{align*}
    (2^k-k-2)\log_2|\mathcal G_{0}|&\leq  \sum_{F\in E(\cH)} \log_2|\Tr(\cG_0,F)|.
\end{align*}
For all $i\in[k-2]$ and $F\in \mathcal J_i$, since $|G_0\cap F|< i+1$ for all $G_0\in \mathcal G_{0}$, we have $|\text{Tr}(\mathcal G_{0},F)|\leq 2^{i+1}-1$. Meanwhile, since  all $\text{Tr}\left(\mathcal G_0,\{C\}\right)\leq 2$, it follows that
\begin{align*}
    (2^k-k-2)\log_2|\mathcal G_{0}|
    &\leq \sum_{C\in\binom{X}{k}}(2^k-k-2-d_i(C))+\sum_{i=1}^{k-2}|\mathcal J_i|\log_2(2^{i+1}-1).
\end{align*}
Since $\sum_{C\in {V\choose k}}d_i(C)=(i+1)|\mathcal J_i|$, we then have
\begin{align*}
(2^k-k-2)\log_2|\mathcal G_{0}|&\leq (2^k-k-2){n\choose k}-\sum_{i=1}^{k-2}(i+1-\log_2(2^{i+1}-1))|\mathcal J_i| \\
    &\le (2^k-k-2){n\choose k}-\sum_{i=1}^{k-2}\frac{(i+1-\log_2(2^{i+1}-1)){n-k\choose i}\nu_i}{{k-1\choose i}} \tag{$|\mathcal J_i|\geq\frac{\nu_i{n-k\choose i}}{{k-1\choose i}}$}\\
    &\le (2^k-k-2){n\choose k}-\sum_{i=1}^{k-2}\frac{(i+1-\log_2(2^{i+1}-1)){n-k\choose i}\tau_i}{(i+1){k-1\choose i}}. \tag{$\frac{\tau_i}{i+1}\leq\nu_i$}
\end{align*}
Therefore, given $H_1,\dots,H_{k-2}$, we have
\begin{align*}
    |\mathcal G_0|&\leq\exp_2\left[{n\choose k}-\sum_{i=1}^{k-2}\delta_i{n-k\choose i}\tau_i\right]. \qedhere
\end{align*}
\end{proof*}

To summarize, every $G$ in \cref{claim:step4final} is uniquely determined by sequentially specifying the following information:
\begin{enumerate}
    \item[(1)] $0\leq \nu_i\leq \frac{\binom{n}{k-i}}{i+1}$ for $i=1,\dots,k-2$;
    \item[(2)] $0\leq \tau_i\leq \binom{n}{k-i}$ for $i=1,\dots,k-2$;
    \item[(3)] $H_i$ for $i=1,\dots,k-2$;
    \item[(4)] $G_0$.
\end{enumerate}
There are at most 
\[
\prod_{i=1}^{k-2}\left(\frac{\binom{n}{k-i}}{i+1}+1\right)\left(\binom{n}{k-i}+1\right)\leq n^{2k^2}
\]
choices for (1)(2). For every fixed sequence $\tau_1,\dots,\tau_{k-2}$, by \cref{claim:step4final-2} and \cref{claim:step4final-3}, the number of choices for (3)(4) is at most
\begin{align*}
    &\exp_2\left[{n\choose k}-\sum_{i=1}^{k-2}\delta_i{n-k\choose i}\tau_i\right] \prod_{i=1}^{k-2}{n^{k-i}\choose\tau_i}{{k\choose i}\tau_i{n\choose i}\choose t_i}\\
    &\le \exp_2\left[{n\choose k}+\sum_{i=1}^{k-2}\left(O(\log n)+H\left(\frac{ i!\rho_i}{{k\choose i}}\right)\binom{k}{i}n^i-\delta_i{n-k\choose i}\right)\tau_i\right] \tag{$\tau_i\geq\frac{t_i}{\rho_i n^i}$} \\
    &\leq \exp_2\left[{n\choose k}-3cn\right]
\end{align*}
(recall that $\rho_i\lll\delta_i$).
Hence given $t_1,\dots,t_{k-2}$, the number of choices for  $G\in\cI_4^*(n, \zeta,\vec\beta)$ in \cref{claim:step4final} is at most 
\[n^{2k^2}\exp_2\left[{n\choose k}-3cn\right]\leq \exp_2\left[{n\choose k}-2cn\right]. \qedhere\] 
\end{proof}

We now can show that most formulae in $\cI_3^*(n, \zeta,\vec\beta)$ also lie in $\cI_4^*(n, \zeta,\vec\beta)$.
\begin{lemma}\label{l:step4}
Suppose $\vec \beta, \zeta$ are all sufficiently small constants (relative to $k$). Then, there exists $c > 0$ such that
$$|\cI_3^*(n, \zeta,\vec\beta)\setminus \cI_{4}^*(n,\zeta,\vec\beta)|< \exp_2\left[{n\choose k}-cn\right].$$
\end{lemma}
\begin{proof}
Let each of the entries of $\vec \beta$ and $\zeta$ be sufficiently small positive constants (as a function of $k$), and let $c > 0$ be sufficiently small relative to all of $\vec \beta, \zeta.$

Observe first that every $G\in \cI_3^*(n, \zeta,\vec\beta)$ cannot contain any clause with exactly $k-1$ negative literals. Indeed, if $\ol u_1\cdots\ol u_{k-1}v\in G$, then since $N_{u_1}\cap\dots\cap N_{u_{k-1}}\cap N_v\neq\emptyset$ (\cref{def:I3}), $G$ would have a non-minimal suformula of the form
$$\{w_1\cdots w_{k-1}v, w_1\cdots w_{k-1}u_1,\dots, w_1\cdots w_{k-1}u_{k-1},  \ol u_1\cdots\ol u_{k-1}v\}.$$

Next, we upper bound the number  formulae $G\in \cI_3^*(n, \zeta)$ that contain at least one clause with exactly $k$ negative literals. Every such $G$ is uniquely determined  by sequentially specifying the following information:
\begin{enumerate}
    \item The lexicographically smallest variable set $\{u_1,\dots,u_k\}\subseteq X$ such that $\ol{u_1}\cdots\ol{u_k}\in G$;
    \item $\bigcap_{i=1}^kN_{u_i}\subseteq\binom{X}{k-1}$;
    \item The non-monotone clauses in $G$;
    \item The monotone clauses in $G$.
\end{enumerate}
The number of such formulae is at most
\begin{align*}
    \underbrace{\binom{n}{k}}_{(1)}\underbrace{2^{\binom{n}{k-1}}}_{(2)}\underbrace{\binom{2^k\binom{n}{k}}{\leq \zeta n^k}}_{(3)}\underbrace{\exp_2\sqb{{n\choose k}-\frac{n-2k+1}{k\cdot 2^{k+1}}{n-k\choose k-1}}}_{(4)}\leq \exp_2\sqb{(1-c){n\choose k}}
\end{align*}
(recall that $c\lll\zeta$). In particular, bounds (3) and (4) arise from the following observations:
\begin{enumerate}
    \item[(3)] $G \in \cI_1^*(n, \zeta)$ has at most $\zeta n^k$ non-monotone clauses.
    \item[(4)] Since $\ol u_1\cdots\ol u_k\in G$, for all $\{w_1,\dots, w_{k-1}\}\in\bigcap_{i=1}^kN_{u_i}$, to avoid the non-minimal subformula
\begin{align*}
    &\{w_1\cdots w_{k-1}v,\, w_1\cdots w_{k-1}u_1,\,\dots,\,w_1\cdots w_{k-1}u_k,\,\ol u_1\cdots\ol u_k\},
\end{align*}
there cannot be any $v\in X$ such that $w_1\cdots w_{k-1}v\in G$. Thus (1)(2) have determined there are at least 
$$\frac{1}{k}\cdot (n-2k+1)\left|\bigcap_{i=1}^kN_{u_i}\right|\geq \frac{n-2k+1}{k\cdot 2^{k+1}}{n-k\choose k-1}$$
monotone clauses that cannot belong to $G$.
\end{enumerate}

Finally, we upper bound the number of formulae $G\in \cI_3^*(n, \zeta,\vec\beta)\setminus \cI_4^*(n, \zeta,\vec\beta)$ that do not contain  clauses with exactly $k$ negative literals. For every such $G$, we write $$G=G_0\cup G_1\cup \dots\cup G_{k-2},$$ where $G_i$ is the set of clauses in $G$ with exactly $i$ negative literals. Since $G$ has at most $\zeta n^k$ non-monotone clauses total, we have $|G_1|,\dots,|G_{k-2}|\leq \zeta n^k$.

Since the number of choices of $t_1,\dots,t_{k-2}$ is at most $(\zeta n^k+1)^{k-2}$, by~\cref{claim:step4final}, there are at most $(\zeta n^k+1)^{k-2} \exp_2\left[{n\choose k}-2cn\right]$ formulae in $|\cI_3^*(n, \zeta,\vec\beta)\setminus \cI_4^*(n, \zeta,\vec\beta)|$ that do not contain clauses with exactly $k$ negative literals. Combining with the previous two observations, we get the bound (since $c\lll\zeta$)
\begin{align*}
    |\cI_3^*(n, \zeta,\vec\beta)\setminus \cI_4^*(n, \zeta,\vec\beta)|&\leq (\zeta n^k+1)^{k-2} \exp_2\left[{n\choose k}-2cn\right]+\exp_2\left[(1-c){n\choose k}\right]\\
&    \leq \exp_2\left[{n\choose k}-cn\right]. \qedhere
\end{align*}
\end{proof}

\subsection{Bound on $|\cI_1^*(n, \zeta)|$}

\begin{proof}[Proof of \cref{t:boundcistar}] 
For $k\geq 4$, suppose $0<c\lll \zeta\lll \beta_1\lll\beta_{i}\lll 1/k$ for all $i \in \{2, \ldots, k-2\}.$
We have
\begin{align*}
    |\cI_1^*(n, \zeta, \beta_1)\setminus\cI_2^*(n, \zeta,\vec\beta)| &\le (k-2)\exp_2\sqb{(1-c){n\choose k}}+\exp_2\sqb{(1-c){n\choose k-1}}|\cI(n-1)|\\
    &\quad+\sum_{i=2}^{k-2}\exp_2\sqb{i(1-c){n\choose k-1}}|\cI(n-i)|, \tag{\cref{l:step2.1,l:step2.2}}
\end{align*}
\begin{align*}
    |\cI_2^*(n, \zeta,\vec\beta)\setminus\cI_3^*(n, \zeta,\vec\beta)|   &\le\exp_2\sqb{\frac{k}{2}(1-c){n\choose k-1}}|\cI(n-k)|, \tag{\cref{l:step3}}
\end{align*}
\begin{align*}
     |\cI_3^*(n, \zeta,\vec\beta)\setminus\cI_4^*(n, \zeta,\vec\beta)| &\le   \exp_2\sqb{{n\choose k}-cn},
\tag{\cref{l:step4}}
\end{align*}
\begin{align*}
     |\cI_4^*(n, \zeta,\vec\beta)| &\le 2^{{n\choose k}}.
\end{align*}
 Combining the above inequalities gives the upper bound
 \begin{flalign*}
  |\cI_1^*(n, \zeta)|   &\le 2^{n\choose k}+(k-2)\exp_2\sqb{(1-c){n\choose k}}+\sum_{i=1}^{k-2}\exp_2\sqb{i(1-c){n\choose k-1}}|\cI(n-i)|\\
    &\qquad+
    \exp_2\sqb{\frac{k}{2}(1-c){n\choose k-1}} |\cI(n-k)|+\exp_2\sqb{{n\choose k}-cn} ,
\end{flalign*}
which yields the desired bound by replacing $c$ by a constant multiple.

For $k=3$, suppose $0<c\lll \zeta\lll \beta_1\lll 1/k$. Since
\begin{align*}
    |\cI_1^*(n, \zeta, \beta_1)\setminus\cI_2^*(n, \zeta,\vec\beta)|\le \exp_2\sqb{(1-c){n\choose k}}+\exp_2\sqb{(1-c){n\choose k-1}}|\cI(n-1)| \tag{\cref{l:step2.1}},
\end{align*}
we get that
\begin{align*}
    |\cI_1^*(n, \zeta)|   &\le 2^{n\choose k}+\exp_2\sqb{(1-c){n\choose k}}+\exp_2\sqb{(1-c){n\choose k-1}}|\cI(n-1)|\\
    &\qquad+
    \exp_2\sqb{\frac{k}{2}(1-c){n\choose k-1}} |\cI(n-k)|+\exp_2\sqb{{n\choose k}-cn}.
\end{align*}

For $k=2$, suppose $0<c\lll \zeta\lll \beta_1\lll 1/k$. Since
\begin{align*}
    &|\cI_1^*(n, \zeta, \beta_1)\setminus\cI_2^*(n, \zeta,\vec\beta)|\le 2\exp_2\sqb{(1-c){n\choose k}} \tag{\cref{l:step2.1-k=2}},\\
    &|\cI_3^*(n, \zeta, \vec\beta)\setminus\cI_4^*(n, \zeta,\vec\beta)|=0 \tag{\cref{l:step4-k=2}},
\end{align*}
we get that
\begin{align*}
    |\cI_1^*(n, \zeta)|   &\le 2^{n\choose k}+2\exp_2\sqb{(1-c){n\choose k}}+
    \exp_2\sqb{\frac{k}{2}(1-c){n\choose k-1}} |\cI(n-k)|.
\end{align*}
\end{proof}

\end{document}